
\magnification=1200
\baselineskip=12pt 


\catcode`@=11 


\font\ninerm=cmr9

\font\eightrm=cmr8
\font\sixrm=cmr6

\font\ninei=cmmi9
\font\eighti=cmmi8
\font\sixi=cmmi6
\skewchar\ninei='177 \skewchar\eighti='177 \skewchar\sixi='177

\font\ninesy=cmsy9
\font\eightsy=cmsy8
\font\sixsy=cmsy6
\skewchar\ninesy='60 \skewchar\eightsy='60 \skewchar\sixsy='60

\font\eightss=cmssq8

\font\ninebf=cmbx9
\font\eightbf=cmbx8
\font\sixbf=cmbx6

\font\ninett=cmtt9
\font\eighttt=cmtt8

\font\ninesl=cmsl9
\font\eightsl=cmsl8

\font\nineit=cmti9
\font\eightit=cmti8


\newskip\ttglue
\def\tenpoint{\def\rm{\fam0\tenrm}%
  \textfont0=\tenrm \scriptfont0=\sevenrm 
\scriptscriptfont0=\fiverm
  \textfont1=\teni \scriptfont1=\seveni \scriptscriptfont1=\fivei
  \textfont2=\tensy \scriptfont2=\sevensy \scriptscriptfont2=\fivesy
  \textfont3=\tenex \scriptfont3=\tenex \scriptscriptfont3=\tenex
  \def\it{\fam\itfam\tenit}%
  \textfont\itfam=\tenit
  \def\sl{\fam\slfam\tensl}%
  \textfont\slfam=\tensl
  \def\bf{\fam\bffam\tenbf}%
  \textfont\bffam=\tenbf \scriptfont\bffam=\sevenbf
   \scriptscriptfont\bffam=\fivebf
  \def\tt{\fam\ttfam\tentt}%
  \textfont\ttfam=\tentt
  \tt \ttglue=.5em plus.25em minus.15em

  \normalbaselineskip=12pt
  \def\MF{{\manual META}\-{\manual FONT}}%
  \let\sc=\eightrm
  \let\big=\tenbig
  \setbox\strutbox=\hbox{\vrule height8.5pt depth3.5pt width\z@}%
  \normalbaselines\rm}

\def\ninepoint{\def\rm{\fam0\ninerm}%
  \textfont0=\ninerm \scriptfont0=\sixrm \scriptscriptfont0=\fiverm
  \textfont1=\ninei \scriptfont1=\sixi \scriptscriptfont1=\fivei
  \textfont2=\ninesy \scriptfont2=\sixsy \scriptscriptfont2=\fivesy
  \textfont3=\tenex \scriptfont3=\tenex \scriptscriptfont3=\tenex
  \def\it{\fam\itfam\nineit}%
  \textfont\itfam=\nineit
  \def\sl{\fam\slfam\ninesl}%
  \textfont\slfam=\ninesl
  \def\bf{\fam\bffam\ninebf}%
  \textfont\bffam=\ninebf \scriptfont\bffam=\sixbf
   \scriptscriptfont\bffam=\fivebf
  \def\tt{\fam\ttfam\ninett}%
  \textfont\ttfam=\ninett
  \tt \ttglue=.5em plus.25em minus.15em
  \normalbaselineskip=11pt
  \def\MF{{\manual hijk}\-{\manual lmnj}}%
  \let\sc=\sevenrm
  \let\big=\ninebig
  \setbox\strutbox=\hbox{\vrule height8pt depth3pt width\z@}%
  \normalbaselines\rm}

\def\eightpoint{\def\rm{\fam0\eightrm}%
  \textfont0=\eightrm \scriptfont0=\sixrm 
\scriptscriptfont0=\fiverm
  \textfont1=\eighti \scriptfont1=\sixi \scriptscriptfont1=\fivei
  \textfont2=\eightsy \scriptfont2=\sixsy \scriptscriptfont2=\fivesy
  \textfont3=\tenex \scriptfont3=\tenex \scriptscriptfont3=\tenex
  \def\it{\fam\itfam\eightit}%
  \textfont\itfam=\eightit
  \def\sl{\fam\slfam\eightsl}%
  \textfont\slfam=\eightsl
  \def\bf{\fam\bffam\eightbf}%
  \textfont\bffam=\eightbf \scriptfont\bffam=\sixbf
   \scriptscriptfont\bffam=\fivebf
  \def\tt{\fam\ttfam\eighttt}%
  \textfont\ttfam=\eighttt
  \tt \ttglue=.5em plus.25em minus.15em
  \normalbaselineskip=9pt
  \def\MF{{\manual opqr}\-{\manual stuq}}%
  \let\sc=\sixrm
  \let\big=\eightbig
  \setbox\strutbox=\hbox{\vrule height7pt depth2pt width\z@}%
  \normalbaselines\rm}

\def\tenbig#1{{\hbox{$\left#1\vbox to8.5pt{}\right.\n@space$}}}
\def\ninebig#1{{\hbox{$\textfont0=\tenrm\textfont2=\tensy
  \left#1\vbox to7.25pt{}\right.\n@space$}}}
\def\eightbig#1{{\hbox{$\textfont0=\ninerm\textfont2=\ninesy
  \left#1\vbox to6.5pt{}\right.\n@space$}}}


\def\fourteenpointbold{\def\rm{\fam0\frtnb}%
  \textfont0=\frtnb \scriptfont0=\tenbf \scriptscriptfont0=\sevenbf
  \textfont1=\frtnbmi \scriptfont1=\tenbmi 
\scriptscriptfont1=\sevenbmi
  \textfont2=\frtnbsy \scriptfont2=\tenbsy 
\scriptscriptfont2=\sevenbsy
  \textfont3=\tenex \scriptfont3=\tenex \scriptscriptfont3=\tenex
  \def\it{\fam\itfam\frtnbit}%
  \textfont\itfam=\frtnbit
  \def\sl{\fam\slfam\frtnbsl}%
  \textfont\slfam=\frtnbsl
  \def\bf{\fam\bffam\frtnrm}
  \textfont\bffam=\frtnrm \scriptfont\bffam=\tenrm
   \scriptscriptfont\bffam=\sevenrm
  \def\tt{\fam\ttfam\frtntt}%
  \textfont\ttfam=\frtntt
  \tt \ttglue=.5em plus.25em minus.15em

  \normalbaselineskip=17pt
  \def\MF{{\manual META}\-{\manual FONT}}%
  \let\sc=\frtnsc
  \let\big=\tenbig
  \setbox\strutbox=\hbox{\vrule height8.5pt depth3.5pt width\z@}%
  \normalbaselines\rm}


\font\frtnb=cmbx12 scaled\magstep1
\font\frtnbmi=cmmib10 scaled\magstep2
\font\frtnbit=cmbxti10 scaled\magstep2
\font\frtnbsl=cmbxsl10 scaled\magstep2
\font\frtnbsy=cmbsy9 scaled\magstep2
\font\frtnrm=cmr12 scaled\magstep1
\font\frtntt=cmtt10 scaled\magstep2
\font\frtnsc=cmcsc10 scaled\magstep2


\font\tenbmi=cmmib10
\font\sevenbmi=cmmib7
\font\tenbsy=cmbsy9
\font\sevenbsy=cmbsy7


\def\ftnote#1{\edef\@sf{\spacefactor\the\spacefactor}#1\@sf
      \insert\footins\bgroup\eightpoint
      \interlinepenalty100 \let\par=\endgraf
        \leftskip=\z@skip \rightskip=\z@skip
        \splittopskip=10pt plus 1pt minus 1pt \floatingpenalty=20000
        \smallskip\item{#1}\bgroup\strut\aftergroup\@foot\let\next}
\skip\footins=12pt plus 2pt minus 4pt 
\dimen\footins=30pc 


\newcount\ftnoteno
\def\Fnote{\advance\ftnoteno by 1 \ftnote{$^{\number\ftnoteno}$}}

\font\sc=cmcsc10



%
\font\ss=cmss10
\font\sixss=cmss8.263   


\def\head#1#2{\headline{\eightss {\ifnum\pageno=1 
\underbar{\raise2pt
 \line {#1 \hfill #2}}\else\hfil \fi}}}

\def\title#1{\centerline{\fourteenpointbold #1}}

\def\author#1{\bigskip\centerline{\sc By #1}}

\def\abstract#1{\vskip.6in\begingroup\ninepoint\narrower\narrower

 \noindent{\bf Abstract.} #1\bigskip\endgroup}

\outer\def\proclaim #1. #2\par{\medbreak \noindent {\sc 
#1.\enspace }
 \begingroup\it #2 \endgroup
 \par \ifdim \lastskip <\medskipamount \removelastskip \penalty
 55\medskip \fi}  

\def\section#1. #2.{\vskip0pt plus.3\vsize \penalty -150 \vskip0pt
 plus-.3\vsize \bigskip\bigskip \vskip \parskip
 \centerline {\bf \S#1. #2.}\nobreak \medskip \noindent}



\def\qed{\quad \vrule height7.5pt width4.17pt depth0pt}
\def\Qed{\qed\ifmmode \relax \else \medbreak \fi}
\long\def\remark{\noindent{\sc Remark.\enspace}}
\long\def\Remark#1\par{\medbreak\remark #1 \medbreak}
\def\example{\noindent{\sc Example:\enspace}}
\def\Example#1\par{\medbreak\example #1 \medbreak}


\def\ref#1 (#2){\par\hangindent=.8cm\hangafter=1\noindent {\sc 
#1}\ (#2).}
\def\rfr#1 (#2){\par\hangindent=.8cm\hangafter=1\noindent {\sc 
#1}\ (#2).}
\def\and{{\rm and }}




\def\cbuldot{{\,\raise.25ex\hbox{$\scriptscriptstyle\bullet$}}}



\def\la#1{\mathop{#1}\limits^{\lower 6pt\hbox{$\leftarrow$}}}
\def\ra#1{\mathop{#1}\limits^{\lower 6pt\hbox{$\rightarrow$}}}

\def\L{{\cal E}_L}


\overfullrule=0pt


\def\bigb{\bigbreak\noindent}         
\def\med{\medbreak\noindent}


\input amssym.def
\input amssym
\hsize 175 truemm
\hoffset=-.2truein
\vsize=9.0 truein

\def\P {{\Bbb P}}  
\def\R {{\Bbb R}} \def\N {{\Bbb N}} \def\Z {{\Bbb Z}}
\def\B {{\Bbb B}}  
 \def\sB {{\cal B}} \def\sC {{\cal C}}
 \def\sE {{\cal E}} 
  
 \def\sK {{\cal K}}

\def\sS {{\cal S}}   
 \def\sW {{\cal W}}  
  \def\Kappa{\sK}
\def\F {{\cal F}}  \def\S {{\cal S}}  \def\G  {{\cal G}}
\def\Y {{\vec Y}} \def\W {{\bf W}} \def\V{{\bf V}}
   
\def\CC{{\mathchoice{{\hbox{\ss Clust}}}{{\hbox{\ss 
Clust}}}{{\hbox{\eightss
 Clust}}}{{\hbox{\sixss Clust}}}}}   
\def\sib{{\mathchoice{{\hbox{\ss sib}}}{{\hbox{\ss 
sib}}}{{\hbox{\eightss
 sib}}}{{\hbox{\sixss sib}}}}}   
   
\def\minver{{\mathchoice{{\hbox{\ss First}}}{{\hbox{\ss 
First}}}{{\hbox{\eightss
 First}}}{{\hbox{\sixss First}}}}}   
\def\pare{{\mathchoice{{\hbox{\ss par}}}{{\hbox{\ss 
par}}}{{\hbox{\eightss
 par}}}{{\hbox{\sixss par}}}}}   
\def\Spine{{\mathchoice{{\hbox{\ss Spine}}}{{\hbox{\ss 
Spine}}}{{\hbox{\eightss
 Spine}}}{{\hbox{\sixss Spine}}}}}   
\def\spine{{\mathchoice{{\hbox{\ss spine}}}{{\hbox{\ss 
spine}}}{{\hbox{\eightss
 spine}}}{{\hbox{\sixss spine}}}}}   
   
\def\reg{R}
  
\def\grow{{\mathchoice{{\hbox{\ss card}}}{{\hbox{\ss 
card}}}{{\hbox{\eightss
 card}}}{{\hbox{\sixss card}}}}}   
\def\minsize{{\mathchoice{{\hbox{\ss minsize}}}{{\hbox{\ss 
minsize}}}
 {{\hbox{\eightss minsize}}}{{\hbox{\sixss minsize}}}}}   
\def\cutsize{{\mathchoice{{\hbox{\ss size}}}{{\hbox{\ss size}}}
 {{\hbox{\eightss size}}}{{\hbox{\sixss size}}}}}   
\def\maxsize{{\mathchoice{{\hbox{\ss maxsize}}}{{\hbox{\ss 
maxsize}}}
 {{\hbox{\eightss maxsize}}}{{\hbox{\sixss maxsize}}}}}   
\def\rate{{\mathchoice{{\hbox{\ss rate}}}{{\hbox{\ss 
rate}}}{{\hbox{\eightss
 rate}}}{{\hbox{\sixss rate}}}}}   

\def\couple{{\mathchoice{{\hbox{\ss couple}}}{{\hbox{\ss couple}}}
 {{\hbox{\eightss couple}}}{{\hbox{\sixss couple}}}}}   
\def\tcpl{\tau_{\couple}}
\def\squid{{\mathchoice{{\hbox{\ss not}}}{{\hbox{\ss 
not}}}{{\hbox{\eightss
 not}}}{{\hbox{\sixss not}}}}}   
\def\TT{T_{\squid}}

\def\var {\varepsilon}

\def\q{\quad} \def\qq{\qquad}

\def\bB{{\Bbb B}}

\def \half {{\scriptstyle {1\over 2}}}
\def \thalf {\textstyle {1 \over 2}}
\def \third {\textstyle {1\over 3}}

\def \threequarters {{\textstyle {3\over 4}}}

\def \qed {\hfill$\square$\par}
\def \llq {\lq\lq}

\parskip=3pt
\def\no{\noindent}
\def\pain{\displaystyle {\sum{\kern-2pt 
\raise8pt\hbox{$\prime$}}}}

\def\UT{{\tilde{U}}}
\def\ee{\varepsilon}
\def\ab{\allowbreak}

\def\wconv{ \,\, {\buildrel w \over \longrightarrow} \,\,} 
\def\L1conv{ \,\, {\buildrel L^1 \over \rightarrow} \,\,}

\overfullrule=0pt
\def\dst{\displaystyle}
\pageno=0

\centerline {{ \bf DIFFUSION-LIMITED AGGREGATION ON A TREE}}
\bigb
\centerline {by}
\bigskip
\centerline {Martin T. Barlow 
\footnote*{The author's research was supported by an NSERC 
operating grant.}
 \hskip 0.95in Robin Pemantle
\footnote\dag{The author's research was supported by NSF grants 
DMS-9300191
and DMS-9353149 and by a Sloan Foundation Fellowship.}
\hskip 0.95in Edwin A. Perkins\ * }
\centerline { U. British Columbia \hskip.88in U. Wisconsin 
\hskip1.16in \ U. 
British Columbia}\vskip.25in

\centerline {Summary}
\bigb
We study the following growth model on a regular $d$-ary tree.  
Points at distance $n$ adjacent to the existing subtree are
added with probabilities proportional to $\alpha^{-n}$, where
$\alpha < 1$ is a positive real parameter.  The heights of
these clusters are shown to increase linearly with their total size;
this complements known results that show the height increases 
only logarithmically when $\alpha \geq 1$.  Results are obtained
using stochastic monotonicity and regeneration results which
may be of independent interest.  Our motivation comes from 
two other ways in which the model may be viewed: as a problem 
in first-passage percolation, and as a version of 
diffusion-limited aggregation (DLA), adjusted so that 
``fingering'' occurs.  
\vskip1in
1991 AMS Subject Classifications:  Primary. 60K40, 60K30, 
60K99\smallskip
\hskip 65 truemm Secondary. 60F05, 60F15, 60K35

\bigskip
Key Words:  diffusion-limited aggregation, trees, random clusters,
                           
\hskip 23 truemm  first-passage percolation,  fingering, regeneration.

\vfill\eject

\noindent {\bf {1. Introduction}}
\medskip

Consider the following dynamical method for growing a subtree of the regular 
$d$-ary tree $\B$.  Initially the subtree consists of only the root.  Vertices 
are then added one by one from among those neighbouring the current subtree.  
The choice of which vertices to add is random, with vertices in generation $n$ 
chosen with probabilities proportional to $\alpha^{-n}$, where $\alpha>0$ is a 
fixed  parameter. 
 Let $A_n$ denote the subtree at step $n$ and let $h(A_n) = \max \{ |x| : x 
\in A_n \}$ denote the maximum height of a vertex in $A_n$.  We are 
interested in the form of the infinite cluster 
$A_\infty := \cup^\infty_{n=0} A_n$, as well as the 
behaviour of $h(A_n)$ and related quantities as $n \rightarrow \infty$.  
In this paper we treat the case $\alpha < 1$; the case 
$\alpha \geq 1$ has already been studied. Our main result, contained in 
Theorem~6.4 and Corollary~6.5, is as follows. 


\no{\bf Theorem 1.1 (Strong Law and CLT)}.  Let $\B$ be the regular $d$-ary 
tree and $0 < \alpha < 1$.  There exist constants $\mu_0 (\alpha , d) \in 
(0,1)$ and  $\sigma^2 =\sigma^2(\alpha,d) > 0$ such that 

\no (a) \quad $\lim_{n \rightarrow \infty} n^{-1} h(A_n)  = 
   \mu_0(\alpha,d)  \;$ a.s.

\no (b) \quad $ n^{-1/2}  (h(A_n) - n \mu_0(\alpha,d)) \wconv N(0 , \sigma^2)
   \; \hbox{ as } n \rightarrow \infty \, . $

\no This model has arisen in a number of different contexts,
and, as we will see below, can be interpreted both as diffusion limited 
aggregation and as first passage percolation on the tree $\B$. 

The case $\alpha = 1$ for binary trees arises in binary search algorithms in 
computer science and has been studied by numerous authors including Pittel 
(1984) and Devroye (1986).  The case $\alpha = 2$ (again for $d=2$) arises in 
an entropy estimation procedure of Ziv (1978) and has been studied by Pittel 
(1985) and (along with $\alpha >1$) by Aldous and Shields (1988).  A 
discussion of the cases $\alpha =1$ and $\alpha =2$ within the general context 
of random search trees in computer science may be found in Chapters 2 and 6, 
respectively, of Mahmoud (1992). 

If $\alpha \geq 1$ it is easy to see that $A_\infty := \cup^\infty_{n=0} A_n$ 
is almost surely the entire $d$-ary tree.  Simply look at the vacant node 
$x_0$ closest to the root, note that $\partial A_n$ has $(d-1) n+1$ points of 
which $x_0$ is the most likely to get filled, and use the fact that 
$\prod^\infty_{k=n} (1-((d-1) k+1)^{-1}) = 0$.  (This argument also works for 
any tree $\B$ of bounded degree.) In the case $\alpha \geq 1$ it is therefore 
natural to look at the growth rates of both $h(A_n)$ and $l (A_n) = \min \{ 
|x|: x\in A^c_n\}.$ 

\bigb {\bf Theorem 1.2}. Let $\B$ be a binary tree (i.e. $d = 2$). 

\item {(a)}  (Pittel (1984), Devroye (1986)) If $\alpha = 1$, then,
writing $\beta_1=4.311...$, $\beta_0 =0.373...$ for the roots of the
equation  $ \half \beta e^{(1-\beta)/\beta} =1,$
$$ \lim_{n\rightarrow \infty} h(A_n) (\log_2n)^{-1} = \beta_1, \,\, 
\hbox{ and } \lim_{n\rightarrow \infty} l(A_n) (\log_2n)^{-1} = 
\beta_0 \hbox { a.s.}$$ 

\item {(b)} (Aldous and Shields (1988))\ \ If $\alpha >1$, then 
$$\lim_{n\rightarrow \infty} h(A_n) (\log_2 n)^{-1} = 
\lim_{n\rightarrow \infty} l(A_n) (\log_2n)^{-1} = 1\ {\rm a.s.}$$ 

\med The case $\alpha = 2$ is also in Pittel (1985, Corollary 1). 
We remark that Aldous and Shields (1988) also study the form of the 
cluster between $l(A_n)$ and $h(A_n)$, and that while they do not state 
explicitly the result for $l(A_n)$ and $\alpha >1$, this is readily 
derived using their methods.  In addition, they mention
the growth dynamics for $\alpha<1$ as an interesting open problem.

Comparing Theorems 1.1 and 1.2 we see that there is a dramatic 
phase transition at $\alpha=1$. For 
$\alpha > 1$ the process $\{ A_n \}$ exhibits the 
same balanced growth as a deterministic procedure in which vertices are added 
in lexicographic order.  Thus the subtree is essentially as short as possible,
($h(A_n) \approx \log_2 n$) and $l(A_n) / h(A_n) \rightarrow 1$.  
  For $\alpha=1$ the results of Pittel and Devroye show that fluctuations on a 
logarithmic scale arise, while for $\alpha<1$ Theorem 1.1(a) implies that the 
subtree is essentially as long as possible $(h(A_n) \approx n)$. (In this case 
it is also easy to see that $A_\infty \not = \B$ -- see Theorem 2.1 and Remark 
2.2 below). 

It would be of interest to study this ``phase transition" at $\alpha =1$ more 
closely.  An understanding of the asymptotics of $\mu_0(\alpha ,d)$ as $\alpha 
\uparrow 1$ would be a step in this direction.  The bounds in~(6.9) show that 
$\mu_0(\alpha ,d) \geq c_1(d) \exp (- \lambda (d) / (1-\alpha ))$, (where
$\lambda(d)>0$), and we suspect that 
 $$ 0 < \lim_{\alpha \uparrow 1} 
     (1-\alpha ) \log \mu_0 (\alpha , d)^{-1} < \infty. \leqno(1.1)$$ 

Our original motivation for studying this model was that it is
a version of diffusion-limited aggregation (DLA) on a tree.
DLA on $\Z^d$ was introduced by 
Witten and Sander (1981) to model aggregates of a condensing metal vapour, and 
since then it has been attracted much interest as a model for various physical 
phenomena (see for example Vicsek (1989)).  DLA is a Markov chain taking 
values in the space of connected finite subsets of $\Z^d$.  Given the 
current configuration, particles diffuse in from infinity according to a 
random walk conditioned to hit the ``boundary" of the current cluster and 
attach themselves to the first point they hit which is adjacent to the current 
cluster.  Although the process has quite a simple description, there are very 
few rigorous mathematical results. 
 If $A_n$ is the cluster at step $n$, $A_0=\{0\}$ and $h(A_n) = \max \{ |x|: 
x\in A_n\}$, then one hard open problem is to prove the existence, and find 
the value of 
$$\beta_d = \lim_{n\rightarrow\infty} {\log h(A_n) \over \log n}. $$ 
Kesten (1987, 1990) has shown the $\lim \sup$ of the above ratio is at most 
$2/(d+1)$.  A conjectured value of $\beta_d$ which agrees quite well with 
numerical simulations is $(d+1) (d^2 +1)^{-1}$ (Lawler (1991, Sec. 2.6)). 
There is no rigorous lower bound for $\beta_d$ aside from the trivial 
$\beta_d\geq d^{-1}$, and hence no rigorous proof of the existence of the 
``fingering" $(\beta_d>d^{-1})$ which simulations suggest. For some further 
surveys of DLA from a mathematical perspective see Lawler (1991) and Barlow 
(1993). 

Analyzing DLA on $\Z^d$ is a hard problem, but the 
same process on a tree is tractable for two reasons: 
\item{(i)} There is a simple formula 
for harmonic measure on the boundary of a cluster 
(see Lemma~1.3 below). 
\item{(ii)} The absence of loops in the graph means that disjoint parts of the 
cluster evolve nearly independently. 

 While the usual heuristic is that a $d$-ary tree (or {\it Bethe Lattice}) 
will exhibit the limiting behaviour for $\Z^d$ as $d \rightarrow \infty$, we
do not know to what extent our model is relevant to DLA on $\Z^d$. However,
it may be interesting to note that, even though the harmonic measure of any
cluster is easy to calculate, the proof of Theorem 1.1(a) is still quite long 
and hard. 

To describe more precisely the connection between the model given above and
DLA on a tree, we need some more notation.
We begin by presenting the notation used to describe an abstract rooted 
tree with no leaves.  In order to be able to move pieces of trees around, we 
find it convenient to view an arbitrary rooted tree (with countably many 
vertices) as a subset of the universal rooted tree ${\bf T}= \cup^\infty_{n=0} 
\N^n$ where $\N^0 = \{0\}$ and $0$ denotes the root of each tree.  (Here 
$\N$ is taken not to include zero.)  
First, the notation for ${\bf T}$ itself is as follows.
If $x \in {\bf T}, |x|=n$ if and only 
if $x\in \N^n$ and $x|j = (x_1,\ldots , x_j)$ for $j\leq |x|$
 (we set $x|0=0)$.  If 
$x,y \in {\bf T}$, let $x\oplus y = (x_1, \ldots , x_{|x|}, y_1, \ldots , 
y_{|y|})$ $(0\oplus x = x)$ and $x\wedge y = (x_1, \ldots , x_j)$ 
where $j$ is 
the largest integer such that $x|j = y|j\ \ (x\wedge y = 0$ if $j=0)$.  
We write $x\leq 
y$ if  $x$ is an ancestor of $y$, i.e. $y = x \oplus z$ 
for some $z$, and we say $x$ is the {\sl parent} of 
$y$ if in addition, $|y|=|x|+1$, and write $x=\pare (y)$.  
Next, we view an arbitrary ordered (the children of each node come 
with an order) tree as a subtree of ${\bf T}$.  Given $m: {\bf T} 
\rightarrow  \N$, inductively 
define the associated locally finite rooted tree (with no leaves) 
$\B\subset {\bf T}$ by 
$$\B (0) = \{ 0\}, \q  \B(n+1) = \{ x \oplus i: x\in 
\B (n),\ i \in \{ 1,2, \ldots , m(x)\}\}, \q \B =  
\cup^\infty_{n=0} \B (n). \leqno(1.2)$$ 
Note that $m(x)$ is the number of children of $x$: the condition 
$m(x)\ge1$
is equivalent to the assertion that $\B$ has no leaves.  
If $m(x)=d$ for all $x$ we obtain the regular $d$-ary tree.  The values 
of $m(x)$ for $x \notin {\B}$ are of course irrelevant.
Intervals $\B(j,k)$, $\B(j,\infty)$, etc. are defined, for example, by 
$\B(j,k) = \cup_{j \leq n\leq k} \B (n)$.  

A subset $A$ of $\B$ is a {\sl rooted subtree} if $\pare (x) \in A$ for all 
$x \in A$. $S$ (respectively, $S_0$) denotes the set of all (respectively, all 
finite) such subtrees.  For $A \in S$, the (external) boundary of $A$ is 
 $$\partial A = \{ x \in A^c: \pare (x) \in A\},$$ 
and the height of $A$ is  $h(A) = \sup \{|x|: x\in A\}$.  

Assume, until we indicate otherwise, that $\B$ is the regular $d$-ary 
tree.  Fix $\alpha >0$.
Let $Q^\mu$ denote the law of the random walk 
$(Y_0 , Y_1 , \ldots)$ on $\B$
started from initial distribution $\mu$, with transition
probabilities given by
$p(x,y) = \alpha / (\alpha + d)$ when $y$ is the parent of $x$,
$p(x,y) = 1 / (\alpha + d)$ when $y$ is a child of $x$
and $x$ is not the root, and $p(0,y) = 1/d$ when $y$ is a child
of $0$. Thus $(Y_n,n\in \Z_+)$ is the random walk on $\B$ 
obtained by assigning conductances $\alpha^{-n}$ to each edge from 
generation $n$ to generation $n+1$. 
Write $\tau (A) = \min \{ n \in \Z_+ : Y_n
\in A \}$ for the hitting time of $A \subset \B$.


The following lemma describes harmonic measure on the boundary of
an arbitrary subtree of $\B$.

\med 
{\bf Lemma 1.3.} Let $A \in S_0$ be a non-empty rooted subtree of $\B$ and 
fix an $N>h(A)$; let $\mu_N = d^{-N} \sum_{x \in \B(N)} \delta_x$ be the 
uniform measure on $\B(N)$.  If $\alpha_0 = \alpha \wedge d$ then 
 $$ Q^{\mu_N} (Y_{\tau(\partial A)} = x | \tau (\partial A) < \infty )
= \alpha_0^{-|x|} (\sum_{y\in \partial A} 
\alpha_0^{-|y|})^{-1},\q x\in \partial A. \leqno (1.3) $$ 
\noindent
{\it Proof.} Note that $|Y_n|$ is a simple reflecting random walk on $\Z_+$ 
which moves to the right with probability $p=d\alpha^{-1} (1+ d\alpha^{-1})^{-
1}$ and to the left with probability $1-p$.  Therefore if $y\in \B(N)$ and 
$0\leq j \leq N$, then 
$$ Q^{\delta_y} (|Y_n| = j \hbox{ for some } n \geq 0)= 
((1-p)/p) \wedge 1)^{N-j} = (\alpha_0/d)^{N-j}.  $$ 
If $\tau (\partial A) < \infty , \ x\in \partial A$, and $|Y_0| = N> h (A)$, 
then $Y_{\tau (\partial A)}=x$ if and only if $Y_0$ is one of the $d^{N-|x|}$ 
descendants of $x$ in $\B(N)$.  Therefore if $x\in \partial A$, then 
$$ \eqalign{ Q^{\mu_N}( Y_{\tau (\partial A)} = x, T (\partial A)<\infty ) 
&= \sum_{y \in \B(N)} d^{-N} 1(x\leq y) Q^{\delta_y} (\tau (\partial 
A)<\infty)\cr 
&= d^{N-|x|} d^{-N} (\alpha_0/d)^{N-|x|} \cr
& = (\alpha_0/d)^N \alpha_0^{-|x|},\cr}$$ 
and the result follows. \qed

As the hitting distribution in (1.3) is independent of $N$ we arrive at 
the following definition. 

\no{\bf Definition.} If $A\in S_0$ and $A\not= \emptyset$, then 
{\sl harmonic measure on $\partial A$ with parameter $\alpha >0$} is given by 
$$ H^\alpha_{\partial A} (x) = \alpha^{-|x|}  (\sum_{y\in 
\partial A} \alpha^{-|y|})^{-1}, \q x\in \partial A.$$ 
Strictly speaking, in view of Lemma 1.3, we should restrict to 
$\alpha \leq d$,  
but allowing $\alpha >d$ is harmless, and in any case in this work we are only 
interested in $\alpha \in (0,1)$.

\med {\bf Definition}. {\sl DLA on $\B$ with parameter}
 $\alpha >0$ is the 
$S_0$-valued Markov chain $(A_n, n\in \N)$ such that $A_1 = \{ 0\} $ and 
$A_{n+1} = A_n \cup \{ D_{n+1}\}$, where 
$$ P(D_{n+1} = x |A_n ) = 
H^\alpha_{\partial A_n} (x), \q x\in \partial A_n. \leqno (1.4)$$ 

It is clear from (1.4) that (for $0< \alpha \leq d$) the growth model 
described at the beginning of the Introduction is exactly DLA on $\B$.

For $\alpha <1$, $H^\alpha_{\partial A_n}$ favours large $|x|$ values in 
$\partial A_n$, for $\alpha >1$ it favours small $|x|$ values and for $\alpha 
= 1$ we have the uniform law on $\partial A_n$ (the ``Eden model").  
 From the perspective of classical DLA the case $\alpha <1$ is of greater 
interest as it is here that we obtain ``fingering", i.e. $h(A_n) \gg 
(\log_2n)^\beta$ for some $\beta >1$. 
This model has been studied in the physics literature -- see Vannimenus et al. 
(1984). As well as some calculations in the case $\alpha=1$, this paper 
describes computer simulations which suggested that $h(A_n)\sim c(\alpha ) n$ 
when $d=2$ and $\alpha <1$. 

A second motivation is that our model is equivalent to {\it
first-passage percolation}.  In first-passage percolation, each edge
is assigned a positive random variable, called a passage time, and
thought of as the time it takes for water (or information, etc.) to
pass from one endpoint to the other. To each vertex $x$ one associates
a time $T(x)$, which is the minimum over paths from the root to $x$ of
the sum of the transit times.  This represents the time before $x$
gets wet when the root is a source of water from time 0 onward. The
model has been widely studied on $\Z^d$, and some quite nontrivial
results have been obtained on trees as well (see Bramson~(1978) and
Pemantle and Peres~(1994) for two examples).  On a tree, it is natural
to rescale edges so that ones more distant from the root are shorter.
For example, the tributaries of a river system behave this way; see
also the limit trees of Aldous~(1991), whose edge lengths decay
exponentially in the distance from the root.  While the edge random
variables in first passage percolation may, of course, have any
distribution on $[0, \infty)$, the simplest case is that of
exponentials. In this case the lack of memory property of the
exponential distribution implies that the cluster $C(t)$ of vertices
which are wet at time $t$ will be a continuous time Markov process on
the space of rooted subtrees $S$.

In Section~2, we prove that DLA on $\B$ arises from first-passage
 percolation when the passage times are independent exponential random
 variables with mean $\alpha^n$ (for an edge from $\B(n)$ to
 $\B(n+1)$).  Throughout the paper, we will use these various
 viewpoints interchangeably; thus we usually refer to the subtree as a
 cluster or as the DLA, but keep the notion of passage times in the
 foreground as well, and in fact, most of our analysis takes place in
 the continuous-time setting of first-passage percolation.

The organization of the paper is as follows.  In Section 2 we show
that DLA can be embedded in first-passage percolation, and give a
number of general results on the model.  
Theorem 2.1 proves that either $A_\infty=\B$ a.s. or else 
the cluster $A_\infty$ has unique infinite line of descent (the
``backbone'').
In the latter case it is possible to decompose $A_\infty$ into the backbone 
plus a sequence of finite clusters 
attached to successive points of the backbone.  These clusters are i.i.d. 
given the ``percolation times'' along the backbone (Theorem~2.5).  
If $\alpha <1/d$,  a simple law of large numbers argument, based on estimates
of cluster sizes in Lemma 2.4, 
shows that $  {\lim\inf}_{n\rightarrow \infty}n^{-1} h(A_n) \geq c(\alpha, d) 
>0$ a.s. (Theorem 2.6).  This result exhibits our basic approach while 
avoiding the technical problems involved when considering the case when 
$\alpha$ is close to 1. 
Finally, Theorem 2.8 gives a general strong Markov property. 

 
The hard work is in Sections 3 and 4. Section 4 contains
the key estimates on
the sizes of the finite clusters, giving bounds in $L^1$ (Theorem~4.4)
and in $L^2$ (Theorem~4.6).  To handle the dependence which arises in these
proofs it is necessary along the way to prove stochastic monotonicity
results for the conditional distribution of the clusters
given the backbone times, and this is done in Section 3 (Lemmas~3.1 and~3.3).  
 Section~5 shows that the percolation times along the backbone form a Markov 
chain, and that this chain converges exponentially fast to its stationary 
measure (Theorem~5.2 and Corollary~5.7).  This paves the way for a Strong Law 
of Large Numbers (Theorem~6.1) and Central Limit Theorem (Theorem~6.2) holding 
for a general class of functionals of the finite clusters.  
 In Section~7 we find a sequence of {\it regeneration times} for cluster, which 
allow it to be decomposed into i.i.d. pieces. Combining these results with the 
theorems in Section 6 we complete the proof of Theorem 1.1. Finally, Section~8 
adds some remarks on the asymptotics of the growth dynamics as $\alpha 
\uparrow 1$, and on a related particle system. 

In most of this work we will only be concerned with the regular $d$-ary tree, 
and the DLA model described above. However, we may note that most of the 
results of Sections 2 and 3 hold for a more general model. First, we may 
consider a  general locally finite rooted tree $\B$, defined by (1.2). 
Secondly, we can fix a function $f: \B \rightarrow (0,\infty)$, and take 
the passage time  between $\pare(x)$ and $x$ to be exponential with mean
$f(x)$. We can then consider a 
 process $A_n$ which has growth probabilities given by
$$ P(A_{n+1}=A_n \cup \{x\} |A_n) = 
  f (x)^{-1} \bigl( \sum_{y\in \partial A_n} f(y)^{-1} \bigl)^{-1}.$$ 

While these extensions involve no new ideas, describing the tree and the
process in this more general setup does require some quite cumbersome
notation. Apart from Theorem 2.1, and some simple estimates on cluster
size in Lemma 2.4, we will therefore restrict our proofs to the case
of a regular tree and $f(x)=\alpha^{|x|}$.
An earlier version of this paper, which treats the general case in Sections 2 
and 3, is available by anonymous ftp from {\tt ftp.math.ubc.ca} (directory
{\tt pub/barlow}.)

A general notational convention is that $c_{i.j}$ denotes a globally defined 
constant introduced in Section $i$, whereas $c$, $c^\prime$, $c(\alpha) , 
\ldots$ may represent different values in different lines.  Dependence of 
$c_{i.j}$ on parameters such as $(\alpha , d)$ will at times be suppressed if 
there is no ambiguity. The integral of a function $\varphi$ with respect to a 
measure $\mu$ (or $\P$) is written $\mu (\varphi )$ (or $\P(\varphi ))$. 
                 
\bigbreak \no {\bf 2. The Continuous Time Model for General Trees}

Let $\B$ be a general locally finite rooted tree with no leaves and 
let 
$f:\B \rightarrow (0, \infty )$.  Consider the $S_0$-valued Markov 
chain 
$\{ A_n\}$ such that $A_1 = \{0\}$ and $A_{n+1} = A_n \cup \{ 
D_{n+1}\}$ where 

$$ \P(D_{n+1} = x |A_n) =\cases{ f(x)^{-1} (\sum_{y\in \partial A_n} 
f(y)^{-1})^{-1}& \hbox{ for }  $x\in \partial A_n$\cr 
0 &\hbox{ for } $x\notin \partial A_n$.\cr} \leqno (2.1)$$ 

\no Note that $f(x) = \alpha^{|x|}$ on the regular $d$-ary tree gives 
the 
cluster dynamics of the previous section. Set $A_\infty = 
\cup^\infty_{n=0} 
A_n$.  We now embed $(A_n, \, n\geq 0)$ in a continuous time 
process.

On some complete $(\Omega, {\cal F} , \P)$ let
 $\{ U_x : x\in \bB\}$ be i.i.d. exponential random variables with 
mean one, and define $T(x) = \sum_{y\leq x} f(y) U_y$. (We set
$T(\pare(0))=0$). In terms of first passage percolation,
if $f(y) U_y$ is the time for liquid to percolate from $\pare (y)$ to 
$y$ 
and $f(0)U_0$ is the time it takes for $0$ to get wet, then $T(x)$ is 
the time 
it takes to percolate to $x$.  Let $C(t) =\{ x: T(x) \leq t\}$ be the 
nodes 
which are wet at time $t$ and define $a(t):C(t)\rightarrow [0,\infty 
)$ by 
$a(t) (x) = t-T(x)$.  
Let $ {\cal S} = \{ (C,a): C\in S, \, a: C\rightarrow [0,\infty )\}$ 
and let ${\cal S}_0$ denote the same set with $S_0$ in place of $S$.  
If $\Delta$ is added to $\R$ as a discrete point, define $\Pi_x: 
{\cal S} \rightarrow\R \cup\{\Delta\}$ for $x$ in $\bB$ by $\Pi_x 
(C,a) = a(x) 1(x\in C) + \Delta 1(x\notin C)$.  Give ${\cal S}$ 
(respectively ${\cal S}_0$) the smallest $\sigma$-field ${\cal F}
({\cal S})$ (respectively, ${\cal F}({\cal S}_0))$ generated by 
the maps $\{ \Pi_x: x\in 
\bB\}$ (respectively their restrictions to ${\cal S}_0$).  
The process $Y(t) = (C(t), a(t))$ is an ${\cal S}$-valued process.  
${\cal G}^0_t = \sigma (T(x) \wedge t, x\in \bB$) and ${\cal G}_t 
= {\cal G}^0_{t+}$ are filtrations on $(\Omega , {\cal F})$ and $T(x)$ 
is a $({\cal G}_t)$ - stopping time for all $x$ in $\bB$ because 
$\{ T(x) < t\} = \{ T(x) \wedge t < t\} \in {\cal G}^0_t$.  
Clearly $Y$ is $({\cal G}_t)$-adapted because $\Pi_x(Y(t))$ takes on 
the value
$t-T(x)$ on $\{T(x)\leq t\}$ and the value $\Delta$ elsewhere.  

Let $\grow_n = \inf \{ t: \#C(t) = n\}$: thus $\{\grow_n\}$ is a sequence of 
a.s. finite $({\cal G}_t)$-stopping times.  Clearly 
$C(\grow_{n+1}) = C(\grow_n) \cup \{ D_{n+1}\}$ 
 where $D_{n+1} \in \partial C(\grow_n)$.  Using the lack of memory property 
of the exponential, it is easy to see, as in Section~1 of Aldous and Shields 
(1988), that the process $A_n = C(\grow_n)$ satisfies (2.1). 
 %
 We therefore may, and shall, take $A_n = C(\grow_n)$ throughout this work. 

\no{\bf Notation.}  Let 
$T (\infty) = \lim_{n \rightarrow \infty} \grow_n \leq \infty$ be the 
time to percolate to $\infty$. It is clear that 
$A_\infty =  C(T (\infty) -) := \{ x: T(x) < T (\infty) \}$. 
Set $T(x,y) = T(y) - T(x)$ if $x\leq y$ and $T(x,y)=\infty$ otherwise.
Let
$$\eqalign{
 T(x, \infty) &= \lim_{n\rightarrow \infty} \min \{ T(x,y) : y\geq 
x, y\in \bB (n)\},\cr 
 T(x-,\infty) &= \lim_{n \rightarrow \infty} \min \{ T(\pare (x),y) : 
y\geq x, y\in \bB(n)\} = f(x)U_x +T(x,\infty).}$$  
Thus $T(x,\infty)$ is the time to percolate from $x$ to $\infty$, while 
$T(x-,\infty)$ is the time to percolate from $\pare(x)$ to $\infty$ 
through
$x$.  If $x \in \B$, let $[x,\infty)$ denote $\{ y \in \B : y \geq x \}$
and similarly define $(x,\infty) , [0,x]$, etc.  For $G \subset \B$
let $\F_G = \sigma (U_x : x \in G)$ and define $\F_n = 
\sigma (U_x : |x| \leq n)$.

 The next result is due to Brennan and Durrett (1986, Sec. 3), but 
we include a proof because the settings are a little different. 

\no{\bf Theorem 2.1}. (a) $\P(T (\infty)<\infty )=0$ or $1$.
\item {(b)} If $\P(T (\infty) = \infty)=1$ then $A_\infty = \bB$ a.s.
\item {(c)}  If $\P(T (\infty) < \infty)=1$ then there is a.s.\ a 
unique infinite line of descent in $A_\infty$, i.e., there is a unique
sequence $\{ \spine_n, n \geq 1\}$ such that
$$ \spine_n \in \B(n)\cap A_\infty, \qquad 
 \spine_n =\pare(\spine_{n+1}), \quad\hbox{for all } n\geq 1.\leqno(2.2)$$


\med {\bf Remark 2.2}. If $\bB$ is a regular $d$-ary tree and $f(x) 
= g(|x|)$ then by comparing $T(\infty)$ with the time to percolate to
$\infty$ along a fixed path we see that 
 $$ \sum_{n=0}^\infty g(n) < \infty \quad \hbox{implies } \quad
   \P(T (\infty) < \infty ) = 1. $$
Rather surprisingly, if $g$ is monotone then the converse is also true:
$\sum g(n) = \infty $ implies  $\P(T (\infty) < \infty ) = 0$ -- see
Pemantle and Peres (1994).
 
\med {\it Proof.} (a) The event
  $  \{ T (\infty) = \infty \} = \bigcap_{x\in \bB(n)} 
\{T(x,\infty )=\infty \}$ is clearly in $  {\cal T} =  
\bigcap^\infty_{n=1}  \F_{\B(n, \infty )}$.  
The latter is a $0-1$ $\sigma$-field by the Kolmogorov $0-1$ theorem. 

\no (b) Clearly $  T (\infty) = \infty$ a.s.  implies $A_\infty = \{ x: 
T(x) < T (\infty) \} = \bB$ a.s.

\no (c) Fix $j\in \N$.  
To show that for each $j$ there is a unique choice of
$\spine_j$ it is enough to show that
$$ \hbox{ There exists a unique } X^j \in \bB (j) \hbox{ such that } 
A_\infty \cap [X^j, \infty) \hbox{ is infinite,}\leqno (2.3)$$ 
Assume (2.3) fails.  
Then there are distinct vertices $x_1, x_2$ in $\bB (j)$ such that 
with positive probability
$A_\infty \cap [x_i, \infty )$ is infinite for both $i=1,2$.
It follows that, with positive probability,
$$ T(\pare (x_1)) + T(x_1-, \infty) = T (\infty) = 
T(\pare (x_2)) + T(x_2-, \infty ).$$ 
Therefore conditional on 
${\cal F}_{j-1},\ T(x_1-,\infty) - T(x_2-,\infty)$ has an atom at the 
${\cal F}_{j-1}$-measurable point $T(\pare 
(x_2)) - T(\pare (x_1))\in \R$.  On the other hand $T(x_1-,\infty)$ 
and $T(x_2-,\infty)$ are independent random variables with densities 
(because $T(x_i-,\infty) = U_{x_i} f(x_i) + T(x_i , \infty)$ with 
$T(x_i,\infty )$ independent of the exponential random variable
$U_{x_i}$) and are jointly independent of ${\cal F}_{j-1}$.  Therefore 
$T(x_1-,\infty) - T(x_2-,\infty)$ has a conditional 
density given ${\cal F}_{j-1}$.  This contradiction completes the 
proof. \qed 

\medskip Similar arguments later will require the following quantitative 
estimate on densities of sums of the variables $U_x$.  The elementary proof is 
omitted. 

\no{\bf Notation}.  Let $s_p (\alpha ) = \prod^p_{i=1} (1 -\alpha^i)^{-1} $, 
$\alpha \in (0,1)$,  $p\in \Z_+ \cup \{ \infty \}$, and write $s(\alpha ) = 
s_\infty (\alpha )$. 

\med 
{\bf Lemma 2.3}.   If $\{ U_i,\ i \in \Z_+\}$ are i.i.d. 
exponential r.v. with mean $1$, and $\alpha \in (0,1)$, then for 
$p\in \Z_+ \cup\{ \infty \}$, $S_p = \sum^p_{i=0} \alpha^i U_i$ 
has a density $u_p (t) \leq s_p (\alpha ) e^{-t}$.   
\med

\no{\bf Notation.} If $x \in \B$, then 
$\B_x = \{ y : x \oplus y \in \B \}$ is the locally finite rooted tree 
of descendants of $x$ (properly translated) and 
$\B_x(n) = \{ y \in \B_x : |y| = n \}$ is the $n^{th}$ generation of the 
nodes in $\B_x$.  If $x \neq 0$ and $y \in \B_x$, let 
$$T^{(x)} (y) = f(x)^{-1} T (\pare (x) , x \oplus y)$$
be the {\it rescaled} percolation times for the tree $\B_x$.
Let 
$$ \eqalign{
T(n) &= \inf \{ T(x) : x \in \B (n) \} \;\; \hbox{ and }, \cr
T^{(x)}(n) &= \inf \{ T^{(x)} (y) : y \in \B_x (n) \} }$$
be the percolation
and rescaled percolation times, respectively, to the $n^{th}$
generation in $\B$ and $\B_x$.  Finally, let 
$$T^{(x)}(\infty) = T(x- , \infty) f(x)^{-1}$$
be the rescaled time to percolate to
infinity in $\B_x$.  Note that if $\B$ is a regular $d$-ary
tree and $f(x) = \alpha^{|x|}$, then $\B_x = \B$ for all $x$,
and $\{ T^{(x)} (y) : y \in \B \}$ has the same distribution as
$\{ T (y) : y \in \B \}$. 

We now derive upper bounds on the $L^1$ and $L^2$ norms of the cluster size at 
a fixed time, conditioned on being finite.  These bounds are crude but their 
proofs are fairly simple,
 and the bounds are good enough to enable us to prove that 
$\liminf n^{-1}h(A_n) >0$ 
in the case when $f(x) = \alpha^{|x|}$ 
and $\alpha \sup_x m(x) < 1$ (c.f.\ Theorem~2.6 and Remark~2.7).  
In Sections~3 and~4 we will have to work much harder to obtain better bounds
(e.g.\ Theorem~4.4) which lead to the linear growth of $h(A_n)$ 
for a regular $d$-ary tree and all $\alpha \in (0,1)$.  

\no{\bf Lemma 2.4}.  Let $f(x) = \alpha^{|x|}$ for some $\alpha \in (0,1)$.
\item {(a)}  $  \P(\# C(t) |T (\infty) >t) \leq c (\alpha ) 
\sum^\infty_{n=0}
\#\bB(n)\alpha^n \qquad  \hbox{ for all } t\geq 0$.
\item {(b)}  $  \P(\# C(t)^2 |T (\infty)  > t) \leq c (\alpha ) 
\sum_{z\in \bB} \alpha^{|z|}   \Big( \sum_{m=0}^\infty \alpha^m \# 
\bB_z (m) \dst \Big) ^2 \quad \hbox{ for all } t\geq 0.$ 
\med
\no{\it Proof}.  (a) If $x\in \bB, t\geq 0$ and $S_\infty$ is as in 
Lemma
2.3, then
$$  \eqalign{
\P(T(x)\leq t<T (\infty)  )
&\leq \P(\P(t-T (x, \infty ) < T(x) \leq t |T(x,\infty )))\cr
&\leq s(\alpha ) \P(e^{-(t-T (x,\infty ))} T(x,\infty ))\qquad \ \hbox{ 
(Lemma 
\ 2.3)}\cr 
&\leq s(\alpha ) e^{-t} \P(\exp(\alpha^{|x| +1} S_\infty) \alpha^{|x| +1} 
S_\infty)\cr 
&\leq (s(\alpha )/(1-\alpha ))^2 e^{-t} \alpha^{|x|+1} \qquad \hbox{ 
(Lemma 
 2.3).}}\leqno (2.4)$$ 
Using the fact that $\P(T (\infty)  >t )\geq \P(U_0>t)=e^{-t}$, we 
conclude that
$$  \eqalign{
\P(\# C(t)|T (\infty)  >t) &= \sum_{x\in \bB} \P(T(x) \leq t<T (\infty)
 )
\P(T (\infty)  >t)^{-1}\cr
&\leq (s(\alpha)/(1-\alpha ))^2 \sum^\infty_{n=0} \# \bB(n) 
\alpha^{n+1}.\cr}$$

\no (b)  Let $x_i \in \bB$, $i=1,2$ be both distinct from 
$x_1 \wedge x_2$.  Let $p= |x_1 \wedge x_2|$,  $x_i^\prime = 
x_i|(p+1)$, 
and write $x_i= x_i^\prime \oplus y_i$. Then
%
$$ \eqalign{ \P( &T(x_1) \vee T(x_2) \leq t <T (\infty)  ) \cr
&\leq \P \bigl( T(x_1 \wedge x_2)+\alpha^{p+1} T^{(x_i^\prime )}(y_i 
) 
\leq t < T(x_1 \wedge x_2) +\alpha^{p+1} T^{(x_i^\prime)}(\infty), 
i=1,2
\bigr) \cr 
 &= \P \big( \prod_{i=1}^2 \P(T^{(x_i^\prime )} (y_i) 
\leq (t-T(x_1 \wedge x_2))\alpha^{-p-1} < T^{(x_i^\prime )}(\infty) 
| {\cal F}_p) \big)\cr 
 &\leq c(\alpha ) \P(1(T(x_1 \wedge x_2)<t ) \exp (-2 \alpha^{-p-1} 
(t-T(x_1 \wedge x_2)))\alpha^{|x_1|+|x_2|-2p}) \q \hbox{ (by (2.4))} 
\cr
&\leq c(\alpha ) e^{-t} \alpha^{|x_1|+|x_2|-p}.\cr } \leqno (2.5)$$

\no In the last line we again used Lemma 2.3 with $S_p = T(x_1 
\wedge x_2)$.  
If either $x_1$ or $x_2$ equals $x=x_1 \wedge x_2$ the above 
inequality is 
clear from (2.4). Decomposing the sum over $x_1,x_2$ in $\bB$ 
according to 
the value of $z=x_1\wedge x_2$, we obtain 
 $$ \eqalignno{
\P( \# C(t)^2 1(T (\infty)  >t))
&=\sum_{x_1, x_2 \in \bB} \P(T(x_1) \vee T(x_2) \leq t <T (\infty) 
)\cr
&\leq \sum_{z\in \bB} \sum^\infty_{m=0} \sum^\infty_{n=0} \#
\bB_z (m) \#\bB_z (n) c(\alpha )e^{-t} \alpha^{m+n+|z|} \q 
\hbox{ (by (2.5))}\cr 
&\leq c (\alpha ) \P(T (\infty)  > t) \Big(\sum_{z\in \B} 
\alpha^{|z|} \big( \sum^\infty_{m=0} \# \bB_z (m) \alpha^ m 
\dst \big)^2\dst \Big).  &\square \cr}$$ 

These bounds are far from optimal.  The very first inequality in the proof 
of (a) ignores the critical fact that $T(x) < T (\infty) $ only holds
for a small proportion of vertices $x$ in $\bB(n)$.  This leads to the more 
restrictive conditions on $\alpha$ in Theorem~2.6 below. 

{\bf In the remainder of the paper we will assume $\B$ is a regular
$d$-ary tree and $f(x) = \alpha^{|x|}$ with $\alpha \in (0,1)$.}
By Remark 2.2 we have $T (\infty)  < \infty$ a.s.  The unique infinite
line of descent in $A_\infty$, defined by (2.2),
is called $\Spine = \{\spine_n, n \in \N\}$, or the ``backbone''
of the cluster $A_\infty$.  The cluster $A_\infty$ may be partitioned 
into $\Spine$ and a collection of disjoint finite clusters which branch
off $\Spine$.  

\no{\bf Notation.} For $n \in \Z_+$ and $x \in \B (n+1)$, let $\{ \sib_j 
(x) :
j < d \}$ denote the siblings of $x$, i.e., the points in
$\{ (x|n) \oplus i : i \neq x_{n+1} , i \leq d \}$ in increasing
order of $i$.  We write $e_{n,j}$ for 
$\sib_j (\spine_{n+1})$ 
and for $n \in \Z_+$ and $1 \leq j < d$, set
$$\CC_{n,j} = \{ x \in \B : e_{n,j} \oplus x \in A_\infty \} \, .$$
Thus $\{ \CC_{n,j} : j < d \}$ are the (possibly empty) clusters
which branch off the backbone in generation $n$.  
Define $a_{n,j} : \CC_{n,j} \rightarrow [0,\infty)$ by
$$ a_{n,j} (x) = (T (\infty)  - T(e_{n,j} \oplus x)) \alpha^{-n-1},$$
and let $Y_{n,j} = (\CC_{n,j} , a_{n,j})$, which is almost surely 
in $\S_0$ by Theorem~2.1.  
Let $W_n = T(\spine_n , \infty) \alpha^{-n-1}$
denote the normalized time to percolate along the backbone 
 from generation $n$ to infinity.  Let $\W$ denote the sequence
$\{ W_n : n = 0 , 1 , 2 , \ldots \}$.  For each $t \geq 0$,
define a law $\nu_t$ on $\S_0$ by
$$\nu_t (\cdot) = \P (Y(t) \in \cdot | T (\infty)  > t) \, .$$

\no{\bf Theorem 2.5}.  Conditional on $(\Spine , \W)$, 
the collection $\{ Y_{n,j} : n \in \Z_+ , j < d \}$ is independent
as $n$ and $j$ vary, and the joint conditional distribution of
each $Y_{n,j}$ is given by
$$\P (Y_{n,j} \in \cdot | \Spine , \W)(\omega) = \nu_{W_n (\omega)} 
\, .$$
\med
The intuitive explanation of this is that 
the only information passed to the subtree beneath $e_{n,j}$ by conditioning 
on the backbone is that the time to infinity inside this subtree has to
be greater than the time along the backbone.  

\no{\it Proof.} Choose $N \in \N$ and $x \in \B (N)$.  Consider 
an event $D = \bigcap_{0 \leq n < N , 1 \leq j < d} D_{n,j}$,
where $D_{n,j}$ is of the form 
$$D_{n,j} = \{ Y_{n,j} = (b_{n,j} , a) \hbox{ for some } a 
   \in F_{n,j} \} ,$$
where $b_{n,j} \in S_0$ and $F_{n,j}$ is the set of
nonnegative functions $\varphi$ on $b_{n,j}$ such that
$\varphi (y) \in F^{n,j,y}$ for each $y$, where $\{ F^{n,j,y} \}$
is a specified collection of measurable sets.  This class of events
generates $\sigma (Y_{n,j} : n \in \Z_+ , j < d)$ and
is closed under finite intersection, so it suffices to
show that
$$\P (D | \Spine ,\W) = \prod \nu_{W_n} (D_{n,j}) \, .$$

If $V_n = U_{\spine_n}$, then clearly $\sigma (\W) = \sigma 
(V_n : n \in \N)$.  For $n < N$, let $R_n = T(x|n , x) + T(x , \infty)$ be
the time to percolate from $x|n$ to infinity through  $x$.  Note that
$\spine_N = x$ if and only if for each $0 \leq n \leq N-1$, the fastest
route from $x|n$ to infinity is through $x$, that is,
$$\{ \spine_N = x \} = \bigcap_{0 \leq n < N} \bigcap_{1 \leq j < d} 
   \{ R_n < T(\sib_j (x|(n+1))- , \infty) \} \, . \leqno(2.6)$$
Note also that on $\{ \spine_N = x \}$ we have 
$$\CC_{n,j} = \{ y' : T^{(\sib_j (x|n+1))} (y') < R_n \alpha^{-n-1} \} $$
and for $y \in \CC_{n,j}$,
$$ \eqalign{
a_{n,j} (y) &= (T (\infty)  - T (\sib_j (x|n+1) \oplus y)) \alpha^{-n-1} 
\cr
&= (T(x|n) + R_n - T(\sib_j(x|n+1) \oplus y)) \alpha^{-n-1} \cr
&= (R_n - T(x|n , \sib_j (x|n+1) \oplus y)) \alpha^{-n-1} \cr
&= R_n \alpha^{-n-1} - T^{(\sib_j (x|n+1))} (y) \, .    } $$
If $B_n$ are measurable subsets of the positive reals, then 
using~(2.6) and the above, we have
$$\eqalign{
\P ( V_n &\in B_n \hbox{ for } 1 \leq n \leq N , \spine_N = x , Y_{n,j} 
\in
   D_{n,j} \hbox{ for } 0 \leq n < N , 1 \leq j < d ) \cr
 = &\P \Big( 1(U_{x|n} \in B_n \hbox{ for } 1 \leq n \leq N)
   \prod_{n=0}^{N-1} \prod_{j=1}^{d-1} \big\{ 
      1(R_n < T(\sib_j (x|n+1)- , \infty)) \, \cr 
   &\times 
      1( \{ y : T^{(\sib_j(x|n+1))} (y) < R_n \alpha^{-n-1} \} 
          = b_{n,j}) \cr 
   &\times
      1(R_n \alpha^{-n-1} - T^{(\sib_j (x|n+1))} (y) \in F^{n,j,y}
         \hbox{ for all } y \in b_{n,j})
   \big\} \Big) \, .} $$
Let $\G(x) = \F_{[0,x]} \vee \F_{[x,\infty)}$.  Observe that
$R_n$ is $\G(x)$-measurable and that $T^{(\sib_j (x|n+1))} (y)$ and
$T(\sib_j (x|n+1)- , \infty)$ are $\F_{[\sib_j (x|n+1) , \infty)}$-
measurable.
The collection of $\sigma$-fields $\F_{[\sib_j (x|n+1) , \infty)}$ for $j < d$ 
and $0 \leq n < N$, together with $\G(x)$, are all mutually independent.  
Condition the above integrand with respect to $\G(x)$ to see that it equals 
$$\eqalign{ 
  \int 1(U_{x|n} (\omega) &\in B_n , 1 \leq n \leq N) 
     \prod_{n=0}^{N-1} \prod_{j=1}^{d-1} 
\Big\{  \P \big( R_n (\omega) < T(\sib_j(x|n+1)- , \infty) \big) \cr
   \times \P \Big( \{ &y : T^{(\sib_j(x|n+1))} (y) < 
      R_n (\omega) \alpha^{-n-1} \} =    b_{n,j} ,\cr 
 R_n (&\omega) \alpha^{-n-1} - T^{(\sib_j (x|n+1))} (y) 
   \in  F^{n,j,y} \hbox{ for all } y \in b_{n,j}  \cr
  \Big| &R_n (\omega) \alpha^{-n-1} 
  < T^{(\sib_j (x|n+1))}(\infty) \Big)
      \Big\} \, d\P (\omega) \,  . }  $$
By~(2.6) the product of the first factors in the curly braces equals 
$\P (\spine_N = x | \G(x)) (\omega)$.  On $\{\spine_N = x \}$ we
have $R_n \alpha^{-n-1} = W_n$ for $n < N$, and so the above leads to
 $$\eqalign{
\P ( V_n &\in B_n \hbox{ for } 1 \leq n \leq N, \spine_N = x , D ) \cr
 &= \int 1(\spine_N = x , U_{\spine_n} \in B_n \hbox{ for } 1 \leq n 
\leq N) \cr
  & \qquad \times \prod_{n=0}^{N-1} \prod_{j=1}^{d-1} \P \big( \{ y : 
T^{(\sib_j(x|n+1))} (y)
     < W_n (\omega) \} = b_{n,j} , \cr
  &\qquad W_n (\omega) - T^{(\sib_j(x|n+1))} (y) \in F^{n,j,y} \hbox{ 
for all } 
     y \in b_{n,j} \cr
  &\qquad | T^{(\sib_j(x|n+1))}(\infty) > W_n (\omega) \big) \, 
d\P(\omega) \cr
= &\int 1(V_n \in B_n \hbox{ for } 1 \leq n \leq N , \spine_N = x )
   \prod_{n=0}^{N-1} \prod_{j=1}^{d-1} \nu_{W_n (\omega)} (D_{n,j})
   \, d\P (\omega) \, \, . }$$
We have used the equivalence in law of $\{ T^{(x)} (y) : y \in \B \}$
and $\{ T(y) : y \in \B \}$ in the last line.  \qed

The above decomposition and the $L^1$ and $L^2$ bounds in Lemma~2.4 allow us 
to use the law of large numbers to establish linear growth of $h(A_n)$ for 
sufficiently small $\alpha$.  The proof illustrates the basic approach we will 
take in Section~6 to obtain the result for all $\alpha < 1$.

\no{\bf Theorem 2.6}. Assume $\alpha \in (0, d^{-1})$. Then 

$$  \liminf_{n\rightarrow \infty}n^{-1}  h(A_n) \geq 
c_{2.1} (\alpha ,d)> 0\ \ \ \hbox{ a.s.}$$ 

\no{\it Proof.} If $\sigma (n) = \inf \{k:h(A_k)=n\}$, the result is 
equivalent 
to 
$$  \limsup_{n\rightarrow\infty}n^{-1} \sigma (n)  \leq c
< \infty \ \ \ \hbox{a.s.} \leqno (2.7)$$

\no Decompose $C(T(n))$ into the backbone vertices in $C(T(n))$ 
(there 
are at most $n+1$) and the portions of the clusters $\{e_{k,i} \oplus
 \CC_{k,i} ,\ 
k< n,\ i< d \}$ which are contained in $C(T(n))$.  This shows that 
 $${\sigma (n) \over n} = {\# C(T(n) ) \over n} \leq {n+1 \over n} + 
n^{-1} \sum^{n-1}_{k=0} \Bigg( \sum^{d-1}_{i=1} \# \CC_{k,i} 
\dst\Bigg) . \leqno (2.8)$$ 
Conditional on $(\W,\Spine)$, $\{ \# \CC_{k,i}: k\geq 0,\ i<d\}$ are 
independent random variables (Theorem~2.5) such that 
 $$  \eqalign{
\P(\# \CC_{k,i}^2 | \Spine, \W) (\omega )
&= \nu_{W_n (\omega)} (( \# C)^2 ) \cr 
& \leq c (\alpha ) \sum_{z\in \bB} \alpha^{|z|} 
(\sum_{m=0}^\infty \alpha^m d^m)^2\cr 
 &\leq c(\alpha ) (1-\alpha d)^{-3}.}\eqno (\hbox{Lemma}\  2.4)$$
Therefore $\mu_{k,i}(\omega) = \P(\# \CC_{k,i} |\Spine,\W)(\omega)$ is also 
uniformly bounded by $c(\alpha ,d)$ say, and the strong law of large numbers 
(applied conditionally) implies that 
 $$  \limsup_{n\rightarrow \infty} n^{-1} \sum^{n-1}_{k=0}\ \  
\sum^{d-1}_{i=1} (\# \CC_{k,i} -\mu_{k,i} ) = 0\ \ \hbox{ a.s.}$$ 
Use this in (2.8) to see that
 $ \limsup_{n\rightarrow\infty} \sigma (n) n^{-1} \leq 1+dc(\alpha , 
d)$ a.s., thus proving (2.7). \qed

\no{\bf Remark 2.7.} In the more general setting of Theorem~2.1 
(general
$\B$ and $f$), it is just as easy to decompose $A_\infty$ into 
$\Spine$ (the backbone) and clusters $\{ \CC_{n,j} : j < m(\spine_n) , 
n \in \Z_+ \}$
which branch off the backbone in generation $n$.  With only 
notational
changes in the proof, it is then possible to derive an analogue of 
Theorem~2.5.  The lack of scaling means the conditional law of 
$Y_{n,j}$
will also depend on the tree $\B_{e_{n,j} (\omega)}$ and the 
appropriately
shifted and rescaled version of $f$.  One can then show that 
Theorem~2.6
remains valid if $\B$ is a rooted tree with no leaves such that $m(x)
\leq d$ for all $x$ and $f(x) = \alpha^{|x|}$ for some $\alpha \in 
(0,d^{-1})$.
The proof is the same.

\no {\bf Notation.} The wide sense past up to $x$ is defined by 
$\sE_x = \F_{(x,\infty)^c}$ and we let $\sE_{x-} = \F_{[x,\infty)^c}$.  

The following strong Markov 
property will be used in Section~6.  It states that if you stop at a
stopping time when $\pare(x)$ has been reached, but $x$ has not 
been reached,
then the remaining times to hit vertices from the subtree rooted at
$x$, rescaled, are equal in law to the original system of hitting times.  
As might be expected, the main difficulty in obtaining this strong Markov
property is in getting the statement right.  

\no{\bf Theorem 2.8.} If $\sigma$ is an a.s.-finite $(\G_t)$-stopping 
time
and $x \in \B$, then for each measurable $B \in [0,\infty)^\B$,
$$\P \Big( ((T (x \oplus y) - \sigma)\alpha^{-|x|} : y \in \B) 
   \in B | \G_{\sigma} \vee \sE_{x-} \Big) 
   \; = \; \P \big( (T(y) : y \in \B) \in B \big) \leqno(2.9)$$
almost surely on $\{ x \in \partial C(\sigma) \}$.

\no{\it Proof.} Note that $\{ x \in \partial C(\sigma) \} = \{ T(\pare 
(x)) \leq \sigma < T(x) \}$ is in $\G_{\sigma}$. 
(Recall that $T(\pare(0)) = 0$).  
Assume first that $\sigma \geq 0$
is constant and consider~(2.9) with $\G_t^o$ in place of $\G_t$.  
Let $\nu$ be the exponential law with mean 1 and let $\nu^F$
denote product measure on $[0,\infty)^F$.  Define $\UT_0 =
(T(x) - \sigma) \alpha^{-|x|}$ and $\UT_y = U_{x \oplus y}$ for
$y \in \B \setminus \{ 0 \}$.  We claim that
$$\P (( \UT_y : y \in \B) \in \cdot | \G_\sigma^o \vee \sE_{x-} )
   = \nu^\B ( \cdot ) \hbox{  a.s. on the event } \{ x \in 
   \partial C(\sigma) \} \, . \leqno(2.10)$$
Indeed, $\G_\sigma^o \vee \sE_{x-} = \sigma (T(y) \wedge \sigma : y 
\geq x) \vee
\sE_{x-} \vee \sigma (1(T(x) > \sigma))$, and $T(y) \wedge \sigma = 
\sigma$
for all $y \geq x$ on $\{ x \in \partial C(\sigma) \}$; the latter
event is measurable with respect to $\sE_{x-} \vee \sigma (1(T(x) > 
\sigma))$,
and so~(2.10) is equivalent to
$$\P (( \UT_y : y \in \B) \in \cdot | \sE_{x-}  \vee
   \sigma (1(T(x) > \sigma)))  
= \nu^\B ( \cdot ) \hbox{  a.s. on the event } \{ x \in 
   \partial C(\sigma) \} \, . \leqno(2.11)$$
Let $F_1$ and $F_2$ be finite subsets of $\B \setminus \{ 0 \}$ and
$[x,\infty)^c$ respectively, let $B_0$ be a measurable set of
nonnegative reals, and let $B_j$ be Borel subsets of $[0,\infty)^{F_j}$
for $j = 1 , 2$.  If $G_1 = \{ (U_{x \oplus y} : y \in F_1) \in B_1 \}$
and $G_2 = \{ (U_y : y \in F_2) \in B_2 \}$, then
$$\eqalign{
\P ( \UT_0 \in &B_0 , G_1 , G_2 , x \in \partial C(\sigma) ) \cr
   &= \P \big( 1( (T(x) - \sigma) \alpha^{-|x|} \in B_0 , G_2 , 
T(\pare(x))
      \leq \sigma < T(x)) \P (G_1 | \sE_x) \big) \cr
   &= \nu^{F_1} (B_1) \P \Big( 
      1(G_2 , T(\pare(x)) \leq \sigma) \cr  
     & \qquad \times  \P (U_x - (\sigma - T(\pare(x))) \alpha^{-|x|}
\in B_0 ,  U_x > (\sigma - T(\pare (x))) 
   \alpha^{-|x|} | \sE_{x-} ) \Big) \, .}$$
Since $T(\pare(x)) \in \sE_{x-}$ and $U_x$ is independent of 
$\sE_{x-}$,
the lack of memory property of the exponential shows that the
conditional expectation term in the last line above is equal to
$$\nu (B_0) \exp [ - \alpha^{-|x|} (\sigma - T(\pare(x))) ] 
   = \nu (B_0) \P (U_x > \alpha^{-|x|} (\sigma - T(\pare(x))) | \sE_{x-}) 
$$
on the event $\{ T(\pare(x)) \leq \sigma \}$.  Substitute this in the
previous equation to conclude that
$$\P (\UT_0 \in B_0 , (\UT_y : y \in F_1) \in B_1 , G_2 , x \in 
   \partial C(\sigma)) = \nu (B_0) \nu^{F_1} (B_1) 
   \P (G_2 , x \in \partial C(\sigma)) \, .$$
It is easy to see this implies~(2.11) and hence~(2.10).  Noting that
$$(T(x \oplus y) - \sigma) \alpha^{-|x|} = \sum_{z \in \B , z \leq y}
   \UT_z \alpha^{-|z|} \, , $$
one derives~(2.9) for $\sigma$ constant and $\G_\sigma^o$ in place 
of $\G_\sigma$.
The entire argument generalizes easily to the case where $\sigma$
is a $(\G_t^o)$-stopping time taking countably many values.

For a general $(\G_t)$-stopping time $\sigma$, choose $(\G_t^o)$-stopping 
times $\sigma_n > \sigma$ so that $\sigma_n$ is a multiple of $2^{-n}$ with 
$\sigma_n \downarrow \sigma$.  Note that $x \in \partial C(\sigma)$ implies 
that $x \in \partial C(\sigma_n)$ for sufficiently large $n$.  Taking limits 
in the above result, we arrive at~(2.9) with $\bigcap_n (\G_{\sigma_n}^o \vee 
\sE_{x-})$ in place of the smaller $\sigma$-field $\G_\sigma \vee \sE_{x-}$.  
The result follows.  \qed 

\no{\bf Remark 2.9.} This result justifies our earlier assertion 
that the process $A_n = C(\grow_n)$ is a Markov
chain satisfying~(2.1).  (Recall that $\grow_n$ is the time
$C(t)$ reaches size $n$.)  

\bigb {\bf 3.  Some Stochastic Monotonicity Lemmas}
\medskip
The stochastic monotonicity results derived in this section will play a
pivotal role in the proof of the key $L^1$ and $L^2$ estimates in Section 4.

There is an obvious isomorphism between $\R^{\B(0,n)}$ and 
$\R^{\{0\}} \times ( \R^{{\bB} (0,n-1)})^{\B (1)}$ which  
we denote $u \mapsto \bar u$ and we extend this isomorphism to 
functions
$\varphi : \R^{\B(0,n)}\rightarrow \R$ by defining 
$\bar\varphi (\bar u) = \varphi (u)$.  
We will use the same notation $(u\rightarrow \bar u)$ to denote the 
isomorphism between $\R^{\bB(n)}$ and $\prod_{y\in \bB(1)}  
\R^{\bB_y(n-1)}$ and hence define $\bar\varphi (\bar u) = \varphi 
(u)$ for 
$\varphi : \R^{\bB(n)} \rightarrow \R$.  The purpose of
this notation is to allow test functions $\varphi$ to be built
recursively, yielding inductive proofs of distributional inequalities.

If $(B, {\cal B}, \leq )$ is a partially ordered ordered measurable 
space and 
$\mu , {\nu}$ are probability laws on $(B, {\cal B})$, $\mu$ is 
{\sl stochastically smaller than} ${\nu}$ (write $\mu \prec {\nu}$) if 
and 
only if $\int \psi d \mu \leq \int \psi d{\nu}$ for every bounded 
measurable 
non-decreasing $\psi : B \rightarrow \R$.  If the law of $X$ is 
stochastically smaller than that  of $Y$, write $X\prec Y$.  If $B=A^T$ 
for 
$A\subset \R$ and $T$ a countable set, we always partially order 
$B$ by 
$f\leq g$ if and only if $f(t) \leq g (t)$ for all $t\in T$.  Partially 
order 
${\cal S}$ by $(C,a) \leq (C^\prime , a^\prime )$ if and only if 
$C\subset 
C^\prime$ and $a\leq a^\prime$ on $C$. 

If $\Gamma$ is a finite subset of $\bB$, $t\in \R$, and $n\in \Z_+$, 
define a probability ${\nu} (\Gamma , n, t)$ on $(-\infty , 
0)^\Gamma$ by 
$$  {\nu} (\Gamma , n,t)(\cdot ) = \P((t-T(x) : x\in \Gamma ) \in 
   \cdot \, |T(n) > t)$$ 
\no (for $t\leq 0$ the conditioning is trivial).  
 \medbreak
\no{\bf Lemma 3.1.}  
${\nu} (\B(0,n) , n,s) \prec {\nu} (\B(0,n), n, t)$  whenever 
$-\infty < s \leq t<\infty $ and $n\in \Z_+$.
\med
\no{\it Proof}.   If $n=0$, ${\nu} ({\bB}(0), 0, s) = 
{\nu} ( {\bB}(0), 0,0)$ for $s\geq 0$ by the lack of memory 
property of the exponential distribution and is trivially stochastically 
non-decreasing for $s\leq 0$. 

Assume for induction that the result holds for $n-1$.  Fix $t \in \R$ 
and
for $z \in \B (1)$, let $\varphi_z$ be bounded measurable functions
on $\R^{\B (0,n-1)}$.  Let $\varphi_0$ be bounded and measurable 
on $\R$
and define $\varphi$ on $\R^{\B (0,n)}$ by $\bar \varphi (\bar u) = 
 \varphi_0 (\bar u_0) \prod_{z \in \B (1)}  \varphi_z (\bar u 
(z))$.  
Also let $\varphi^{(\alpha)}
(u) = \varphi (\alpha u)$, and set $g(u) = \P (T(n-1)  > u)^d$.  Use
the independence of the vectors $V_z = (T^{(z)} (x) : x \in \B (0,n-1))$
and integrate out the exponential variable $U_0$ to see that
$$\eqalign{
\P &(T(n)  > t) \nu (\B(0,n) , n , t) (\varphi) \cr
   &= \int_0^\infty e^{-s} \varphi_0 (t-s) \cr
  & \quad \times \P \big( \prod_{z \in \B (1)} 
     1(s + \alpha T^{(z)} (n-1) > t) \varphi_z (t-s-\alpha T^{(z)} (x) :
     x \in \B (0,n-1)) \big) \, ds \cr
   &= \int_0^\infty e^{-s} g((t-s) \alpha^{-1}) \varphi_0^{(\alpha)}
     ((t-s)/ \alpha) \cr
 &\quad \times \prod_{z \in \B (1)} \P (\varphi_z^{(\alpha)} ((t-s)
     \alpha^{-1} - T^{(z)} (x) : x \in \B(0,n-1)) 
       | ~\; T^{(z)} (n-1) > {t-s \over \alpha} ) \, ds \cr
   &= \int_0^\infty e^{-s} g((t-s) \alpha^{-1}) \big[ \bigotimes_{j=1}^d
     \nu ( \B(0,n-1) , n-1 , (t-s) \alpha^{-1}) \big] (
     {\overline{\varphi^{(\alpha )}}} ((t-s) \alpha^{-1} , \cdot )) \, ds
\, . }$$
The restriction on $\varphi$ is then easily removed, to yield the 
above
equality for all bounded measurable $\varphi$ on $\R^{\B(0,n)}$.  If
$$\psi (v) =  \big[ \bigotimes_{j=1}^d \nu ( \B(0,n-1) , n-1 , v) \big]
    ({\overline{\varphi^{(\alpha )}}} (v , \cdot )) \, ,$$
a change of variables in the above integral leads to
$$\nu (\B (0,n) , n , t) (\varphi) = (e^t \P (T(n)  > t))^{-1}
   \alpha \int_{-\infty}^{t/\alpha} e^{\alpha v} g(v) \psi (v) \, dv \, 
.$$
Take $\varphi$ to be the constant function 1, hence $\psi \equiv 1$ 
also, 
to see that
$$e^t \P (T(n)  > t) = \alpha \int_{-\infty}^{t/\alpha} e^{\alpha v} 
   g(v) \, dv \, ,$$
and therefore conclude
$$ \nu (\B(0,n) , n , t) (\varphi) = \int_{-\infty}^{t/\alpha} 
   e^{\alpha v} g(v) \psi (v) \, dv \big( \int_{-\infty}^{t/\alpha}
   e^{\alpha v} g(v) \, dv \big)^{-1} \, . \leqno(3.1)$$
If $\varphi$ is non-decreasing, then $\psi$ is nondecreasing by
the induction hypothesis.  Since~(3.1) expresses $\nu (\B (0,n) , n , t)
(\varphi)$ as a weighted average of values of $\psi$ with respect
to a weighting measure that stochastically increases in $t$, it follows 
that $\nu (\B (0,n) , n , t) (\varphi)$ is nondecreasing in $t$ and
the induction is complete.   \qed

\med \no{\bf Corollary 3.2}. ${\nu}_s \prec {\nu}_t$ whenever $0\leq 
s\leq t$.

\no{\it Proof}. Define $\Gamma : [0,\infty) \rightarrow \R^{\bB} $ 
by 
$\Gamma (t) (x) = t-T(x)$ and $\psi : \R^{\bB} \rightarrow {\cal S}$ 
by 
 $$ \psi(\gamma ) = (C,a),\q  C=\{ x: \gamma (x) \geq 0\}, \q 
a(x) = \gamma (x)\hbox{  for } x\in C:$$
here $\gamma \in \R^\B$.
Then $\psi$ is non-decreasing and $Y(t) = \psi (\Gamma (t))$.  
It therefore suffices to show that $t\rightarrow \P(\Gamma (t) 
\in \cdot \, | T(\infty) > t)$ is stochastically non-decreasing, and for 
this it 
suffices to fix $m \in \N$ and show that if $\Gamma_m (t) = 
\Gamma (t) 
|_{ {\bB}(0,m)},$ then $t\rightarrow \P(\Gamma_m (t) \in \cdot \, |
T(\infty) 
> t)$ is stochastically non-decreasing (see Kamae, Krengel and O'Brien 
(1977, 
Proposition 2)).  By taking limits, one reduces this in turn to proving
$$  t\rightarrow \P(\Gamma_m (t) \in \cdot \, |T(n) > t) \hbox{ is 
stochastically non-decreasing for all } n\geq m.\leqno (3.2)$$ 
\no It suffices to consider (3.2) with $m=n$ since decreasing $m$ 
only 
weakens the conclusion.  But (3.2) with $m=n$ is precisely the 
conclusion of 
Lemma 3.1.\qed 

The inductive arguments from here on require 
a second set of percolation times, defined analogously to the
first but not including the percolation time at the root of
each subtree.  We apologize for doubling the notation but
promise not to do it again.  

\no{\bf Notation.} Let
$$\eqalign{
\TT(n) &= \inf \{ T (0 , x) : x \in \B(n) \} = T(n)  - U_0 \cr
\TT^{(x)}(y) &= T(x , x \oplus y) \alpha^{-|x|} = T^{(x)} (y) - U_x \cr
\TT^{(x)}(n) &= \inf \{ \TT^{(x)} (y) : y \in \B (n) \} = T^{(x)}(n) - U_x 
   \, . }$$
In short, the $T^{(x)}$'s include a contribution from $U_x$ while
the $\TT^{(x)}$'s do not;  times with superscripts are rescaled.
It is evident that $\{ \TT^{(x)} (y) : y \in \B \}$
is equal in law to $\{ \TT^{(0)} (y) : y \in \B \}$.
For $|z|\geq 1$, let $\minver_n (z)$ denote the a.s. unique vertex 
in $\bB(n)$ such that $T^{(z)}(n) = T^{(z)} (\minver_n (z))$, so
that $z \oplus \minver_n (z)$ is the first descendant of $z$ in
generation $|z|+n$ to be reached. (Of course, this is not necessarily the
one through which infinity is reached from $z$).  
Let $\minver_n$ be the a.s.\ unique vertex in 
$\bB(n)$ such that $T(n) =T(\minver_n)$.  For $z \in \bB(n)$ let 
$\P_z(\cdot ) = \P(\cdot \, |\minver_n = z)$ and 
$\{ \mu (z,t) : t\geq 0\}$ be a set of regular 
conditional probabilities on $(-\infty , 0]^{\bB(n)}$ for 
$$\P_z ((\TT(n) - T(0,x) : x\in \bB(n))\in \cdot \, |\TT(n)= t).$$  

With the available 
symmetry, we could have defined $\mu (n , t)$ instead of
$\mu (z,t)$, but in this case keeping greater generality 
also reduces confusion of types.  
\med  
{\bf Lemma 3.3.}  For any $z\in \bB \setminus \{ 0 \}$, there 
is a version of the set $\{ \mu (z,t) : t \geq 0 \}$ such that 
$$\mu ( z,s) \prec \mu (z, t) \hbox{ for } 0\leq s\leq t.$$

\no{\bf Question:} Is there a useful description of 
the increasing limit, $\mu (z, \infty)$ of the measures
$\mu (z,t)$?  How about the measure $\mu (z,0)$?  In either
case, sending $|z|$ to infinity and rescaling by
$\alpha^{-|z|}$ should then yield a locally finite point process
on $(-\infty , 0]$.  

\no{\it Proof}.  We start by establishing a pair of auxiliary results,
(3.3) and (3.4), whose proofs could be omitted on a first reading.
Let $z = (z_1 , z_2 , \ldots , z_n) = z_1 \oplus z' \in \B (n)$, $n \geq 
1$.
 The first result is that, conditional on $\minver_n = z$ and on the 
percolation
times from the vertex $z_1$, the vectors $\{ T^{(x)} (y) : y \in \B \}$
are i.i.d. as $x$ ranges over the other first generation vertices, 
and distributed as $\{ T (y) : y \in \B \}$ conditioned
on $T(n-1) > \TT(n)(\omega) \alpha^{-1}$.  More precisely, letting 
$T$ denote the vector $(T(x) : x \in \B)$, we show that
$$ \P_z \big( (T^{(x)} : |x| = 1 , x \neq z_1 ) \in \cdot | 
\F_{[z_1,\infty)} 
  \big) (\omega) = \bigotimes_{j=1}^{d-1} \P ( T \in \cdot | 
  T(n-1) > \TT(n)(\omega) \alpha^{-1} )\, ,  \leqno (3.3)$$
i.e., the RHS of~(3.3) defines a regular conditional probability for the
left side.
 
To prove~(3.3), note that 
$$ \{ \minver_n = z \} = \{ T^{(z_1)}(n-1) = \alpha^{-1} T(0,z) \} \cap
   \bigcap_{|x| = 1, x \neq z} \{ T^{(x)}(n-1) > \alpha^{-1} T(0,z) \}$$
almost surely.  This shows that if $\varphi_x : \R^\B \rightarrow \R$ 
are
bounded and measurable for $x \in \B (1)$, then
$$\eqalign{
\P \big( 1 (\minver_n = z) \prod_{|x| = 1} \varphi_x (T^{(x)} ) \big) 
   = \P \Big( &\varphi_{z_1} (T^{(z_1)}) 1(T^{(z_1)}(n-1) = \alpha^{-1}
     T(0,z)) \cr
    &\times \prod_{x \in \B(1) \setminus \{ z_1 \} } 1(T^{(x)}(n-1) > 
     \alpha^{-1} T(0,z)) \varphi_x (T^{(x)}) \Big) \, .}$$  
The term in front of the product is $\F_{[z_1,\infty)}$-measurable 
(since $\F_{[z_1,\infty)}$ is just $\sigma (T^{(z_1)})$) and 
conditional on $\F_{[z_1,\infty)}$, the vectors $T^{(x)}$ for
$x \neq z_1$ are i.i.d. copies of $T$.  Thus
$$\eqalign{
\P_z &\big( \prod_{x \in \B (1)} \varphi_x (T^{(x)}) \big) \P 
(\minver_n = z) 
 \cr &= \int \varphi_{z_1} (T^{(z_1)} (\omega)) 
       \prod_{x \in \B(1)
        \setminus \{ z_1 \} } \P (\varphi_x (T) | T(n-1) > \alpha^{-1} 
        T(0,z) (\omega)) \cr
     &\times \big[ 1(T^{(z_1)}(n-1) (\omega) = \alpha^{-1} T(0,z) 
(\omega))
        \prod_{x \in \B(1) \setminus \{ z_1 \} } \P (T^{(x)}(n-1) > 
        \alpha^{-1} T(0,z) (\omega)) \big] \, d\P (\omega) \, \, .}$$ 
The term in square brackets equals $\P (\minver_n = z | \F_{[z_1 , 
\infty)})
(\omega)$ and so we may conclude that
$$\P_z \big( \prod_{x \in \B (1)} \varphi_x (T^{(x)}) \big) = \int 
   \varphi_{z_1} (T^{(z_1)}) \prod_{x \in \B(1) \setminus \{ z_1 \} } 
\P 
   (\varphi_x (T) | T(n-1) > \alpha^{-1} T(0,z) (\omega)) \, d\P_z 
(\omega)\, .$$
 From this,~(3.3) follows immediately upon noting that $T(0,z) = 
\TT(n)$
a.s.\ with respect to $\P_z$.

The second result is that if the joint distribution of
the times $\TT^{(z_1)} (\cdot)$ is conditioned on $\TT^{(z_1)}(n-1)$, 
on
$U_{z_1}$ and on $\minver_n = z$, then the value of $U_{z_1}$ is 
irrelevant,
and the information $\minver_n = z$ may be replaced by the weaker
$\minver_{n-1} (z_1) = z'$.  Formally, if $\P_{z'}' (\cdot)$ denotes
$\P (\cdot | \minver_{n-1} (z_1) = z')$ and $\varphi : \R^\B 
\rightarrow \R$
is bounded and measurable, then, letting $\TT^{z_1}$ denote the 
vector
$(\TT^{z_1}(x):x \in\B)$, we claim
$$\P_z (\varphi (\TT^{(z_1)}) | \TT^{(z_1)}(n-1) , U_{z_1}) = \P_{z'}'
   (\varphi (\TT^{(z_1)}) | \TT^{(z_1)}(n-1)) \hbox{  a.s.} \leqno 
(3.4)$$
Note that $\P_z$ is absolutely continuous with respect to $\P_{z'}'$,
so that in~(3.4) we are asserting that the right side defines a version
of the left side.

To prove~(3.4), first drop the conditioning on $U_{z_1}$ from the LHS
by observing $U_{z_1}$ to be independent of $\sigma (\TT^{(z_1)} 
(\cdot))$.
Next, write
$$ \{ \minver_n = z \} = \{ \minver_{n-1} (z_1) = z' \} \cap 
\bigcap_{x \in B(1)
   \setminus \{ z_1 \}} \{ U_{z_1} + \TT^{(z_1)}(n-1) < T^{(x)}(n-1) \} 
.$$
The independence of $\TT^{(z_1)}$ and $(U_{z_1} , T^{(x)} :
x \in B(1) \setminus \{ z_1 \})$ shows that
$$\P (\minver_n = z | \TT^{(z_1)}) = \P (\minver_n = z | 
1(\minver_{n-1} (z_1) = z') ,
   \TT^{(z_1)}(n-1)) \hbox{ a.s.} $$
If $\psi$ is a bounded, measurable real function, then 
$$\eqalign{
\P (\minver_n &= z) \P_z (\varphi (\TT^{(z_1)}) \psi (\TT^{(z_1)}( n-
1))) \cr
   = &\P \big( \varphi (\TT^{(z_1)}) \psi (\TT^{(z_1)}(n-1)) \P 
(\minver_n = z | 
      1(\minver_{n-1} (z_1) = z'), \TT^{(z_1)}(n-1)) \big) \cr
   = &\P \big( \P (\varphi (\TT^{(z_1)}) | 1(\minver_{n-1} (z_1) = 
z'),\TT^{(z_1)
     }(n-1)) \psi (\TT^{(z_1)}(n-1)) 1(\minver_n = z) \big) \, .}$$
A separate consequence of the definition of conditional expectation is 
that almost surely on $\{ \minver_{n-1} (z_1) = z' \}$ (which 
contains 
$\{ \minver_n = z \}$)
$$\P (\varphi (\TT^{(z_1)}) | 1(\minver_{n-1} (z_1) = z') , 
\TT^{(z_1)}(n-1))
   = \P_{z'}' (\varphi (\TT^{(z_1)}) | \TT^{(z_1)}(n-1)) \, .$$
Combining this with the previous identity shows that
$$\eqalign{
\P_z (\varphi (&\TT^{(z_1)}) \psi (\TT^{(z_1)}( n-1))) \cr 
= &\P (\P_{z'}' (\varphi
   (\TT^{(z_1)}) | \TT^{(z_1)}(n-1)) \psi (\TT^{(z_1)}(n-1)) 
   1(\minver_n = z)) 
   \P (\minver_n = z)^{-1} \cr
= &\P_z (\P_{z'}' (\varphi (\TT^{(z_1)}) | \TT^{(z_1)}(n-1)) \psi 
   (\TT^{(z_1)}(n-1))) }$$ 
and~(3.4) follows.   

The lemma is now proved by induction on
$n$ as follows.  For $n = 1$, let $Q_z = \delta_0$ and for $x \in \B (1) 
\setminus \{ z \}$,
let $Q_x$ be the law of $-\alpha$ times an exponential of mean 1.
The lack of memory property for the exponential law shows that
$\bigotimes_{x \in \B (1)} Q_x$ is a version of $\mu (z,t)$ for
all $t \geq 0$, so the result holds with equality of all the laws.

Assume the result now for $|z| < n$, retaining the above notation.
Let $\varphi_x : \R^{\B (n-1)} \rightarrow \R$ be bounded and 
measurable,
with $\varphi : \R^{\B (n)} \rightarrow \R$ given by $\bar
\varphi (\bar u) = \prod_{x \in \B (1)} \varphi_x (\bar u (x))$, i.e., 
by $\varphi ( u (\cdot)) = \prod_{x \in \B (1)} \varphi_x ( 
u (x \oplus y) : y \in \B (n-1))$.  Note that if $x \in \B (n-1)$
then almost surely with respect to $\P_z$,
$$\eqalign{ 
\TT(n)&= \alpha (U_{z_1} + \TT^{(z_1)}(n-1)) \cr
\TT(n)- \alpha T^{(z_1)} (x) &= \alpha (U_{z_1} + \TT^{(z_1)}(n-1)) 
   - \alpha (U_{z_1} + \TT^{(z_1)} (x)) \cr
&= \alpha (\TT^{(z_1)}(n-1) - \TT^{(z_1)} (x))
   \, . } \leqno(3.5) $$
Condition on $\F_{[z_1 , \infty)}$ and use~(3.3) and~(3.5) to see that
$$\eqalign{
\P_z &\big( \varphi (\TT(n)- T(0,x) : x \in \B (n)) | \TT^{(z_1)}(n-1) ,
   U_{z_1} \big) (\omega) \cr
&= \P_z \big( \prod_{y \in \B (1)} \varphi_y  (\TT(n)- \alpha T^{(y)} 
(x) : 
   x \in \B (n-1)) | \TT^{(z_1)}(n-1) , U_{z_1} \big) (\omega) \cr
&= \prod_{y \in B(1) \setminus \{ z_1 \}} \P \big( \varphi_y (\alpha
   (\TT(n)(\omega) \alpha^{-1} - T(x)) : x \in \B (n-1)) | T(n-1) >
   \TT(n)(\omega) \alpha^{-1} \big) \cr
&\qquad \times \P_z \big( \varphi_{z_1} (\alpha (\TT^{(z_1)}(n-1) - 
\TT^{(z_1)} (x)) : 
   x \in \B (n-1)) | \TT^{(z_1)}(n-1) , U_{z_1} \big) (\omega) \, .}$$
Recall the notation $\varphi^{(\alpha)} (u) = \varphi (\alpha u)$.
The previous line now becomes
$$\eqalign{
\P_z \big( \varphi (\TT(n)&- T(0,x) : x \in \B (n)) | \TT^{(z_1)}(n-1) ,
   U_{z_1} \big) (\omega) \cr
 = \Big\{ &\prod_{y \in \B(1) \setminus \{ z_1 \}} \big[ \nu ( \B (n-
1)
   , n-1, \TT(n)(\omega) \alpha^{-1}) (\varphi_y^{(\alpha)}) \big] 
\Big\} \cr
&\times \P_z \big( \varphi_{z_1}^{(\alpha)} (\TT^{(z_1)}(n-1) - 
\TT^{(z_1)} (x) :
   x \in \B (n-1)) | \TT^{(z_1)}(n-1) , U_{z_1} \big) (\omega) \; .}$$ 
Now use~(3.4) to equate the above to
$$\eqalign{
\Big\{ \prod_{y \in \B(1) \setminus \{ z_1 \}} &\big[ \nu ( \B (n-1)
   , n-1, \TT(n)(\omega) \alpha^{-1}) (\varphi_y^{(\alpha)}) \big] 
\Big\} \cr
&\times \P_{z'}' \big( \varphi_{z_1}^{(\alpha)} (\TT^{(z_1)}(n-1) - 
\TT^{(z_1)} (x) : 
   x \in \B (n-1)) | \TT^{(z_1)}(n-1) \big) (\omega) \, .}$$ 
Use the identity $\P_{z'}' ( (\TT^{(z_1)} (x) : x \in \B) \in 
\cdot) = \P_{z'} ((T(0,x) : x \in \B) \in \cdot)$ (in words:
looking in the subtree from $z_1$ conditioned on $z_1 \oplus z'$ 
being the first of its generation reached among the subtree is 
the same as looking in the whole tree conditioned on $z'$ being
the first in its generation) to get
$$\eqalign{
\P_z &\big( \varphi (\TT(n)- T(0,x) : x \in \B (n)) | \TT^{(z_1)}(n-1) ,
   U_{z_1} \big) (\omega) \cr
&= \mu (z' , \TT^{(z_1)}(n-1) (\omega)) (\varphi_{z_1}^{(\alpha)}) 
    \prod_{y \in \B(1) \setminus \{ z_1 \}} \big[ \nu ( \B (n-1)
   , n-1, \TT(n)(\omega) \alpha^{-1}) (\varphi_y^{(\alpha)}) \big] \, 
.}$$
Condition both sides of the above with respect to $\sigma (\TT(n))$
which is contained in $\sigma (\TT^{(z_1)} , U_{z_1}) \vee \{ \P_z 
\hbox{-null sets} \}$ and use~(3.5) together with the independence 
of $U_{z_1}$ and $\TT^{(z_1)}(n-1)$ to conclude that
$$\eqalign{
\P_z \big( &\varphi (\TT(n)- T(0,x) : x \in \B (n)) | \TT(n)\big) 
(\omega) \cr
= &\int_0^{\TT(n)(\omega)/\alpha} \Big[ \mu (z',s) \times \kern-2pt 
   \bigotimes_{ y \in \B (1) \setminus \{ z_1 \}} \nu ( \B (n-1) , n-1 , 
   \TT(n)(\omega) \alpha^{-1}) \Big] \cr
&\qquad \times  (\overline {\varphi^{(\alpha)}}) 
 \exp (- ((\TT(n)(\omega) / \alpha) - s)) \P_z (\TT^{(z_1)}(n-1) \in ds) \cr
&\qquad \times   
  \big( \int_0^{\TT(n)(\omega)/\alpha}  \exp (- ((\TT(n)(\omega) / 
\alpha) - s))
  \P_z (\TT^{(z_1)}(n-1) \in ds) \big)^{-1} \cr
= &\int_0^{\TT(n)(\omega)/\alpha} \Big[ \mu (z',s) \times 
\bigotimes_{
   y \in \B (1) \setminus \{ z_1 \}} \nu ( \B (n-1) , n-1 , 
   \TT(n)(\omega) \alpha^{-1}) \Big]  \cr
&\qquad \times  (\overline {\varphi^{(\alpha)}}) 
   e^s \P_z (\TT^{(z_1)}(n-1) \in ds) 
 \big( \int_0^{\TT(n)(\omega)/\alpha}  e^s
  \P_z (\TT^{(z_1)}(n-1) \in ds) \big)^{-1} .}$$
Let $\mu (z,t)(\varphi )$ be defined by the above expression with 
$t$ in 
place of $\TT(n)(\omega )$.  The above shows that $\{ \mu (z,t):t \ge 
0 \}$
are regular conditional probabilities for the required conditional 
distributions.  
By induction, there is a version of the measures $\mu (z' , s)$ that is
stochastically nondecreasing in $s$, while Lemma~3.1 shows that 
$\bigotimes_{y \in \B (1) \setminus \{ z_1 \}}
\nu (\B (n-1) , n-1 , t \alpha^{-1})$ is stochastically
nondecreasing in $t$.  Thus for each bounded, nondecreasing 
function,
$\varphi$ on $\R^{\B (n)}$,
$\mu (z,t) (\varphi)$ may be written as
a weighted average on $[0,t/\alpha]$ of a nondecreasing function:
$$\int_0^{t/\alpha} \psi (s,t) \gamma (ds) \gamma ([0,t/\alpha])^{-1} 
\, ,$$
where $\psi$ is nondecreasing in each variable (combining the effects of 
$\mu (z' , s)$ and the measures \break 
$\nu (\B (n-1) , n-1 , t / \alpha)$), and 
$\gamma$ is the locally finite measure given by the
$e^s$ term times the density of $\TT^{(z_1)}(n-1)$.  Such
a weighted average is clearly nondecreasing in $t$. \qed
\med
{\bf Corollary~3.4}.  If $\theta \geq 0$ and $n \in \N$, there is
a nondecreasing version of 
$$\Kappa_n (t) = \P \big( \sum_{x \in \B (n)} \exp ( -\theta 
   (T(0,x) - \TT(n))) | \TT(n) = t \big) .$$

\no{\bf Proof}. For $z \in \B (n)$, let 
$$\Kappa_{n,z} (t) = \P \big( \sum_{x \in \B (n)} \exp ( -\theta 
   (T(0,x) - \TT(n))) | \minver_n = z , \TT(n) = t \big) .$$
Since the sum is increasing in each $\TT(n) - T(0,x)$, Lemma~3.3
implies that each $\Kappa_{n,z} (t)$ has a nondecreasing 
version.  But by symmetry, $\Kappa_{n,z}$ is independent
of $z$, so is equal to $\Kappa_n$ almost surely.  \qed

The following elementary result is proved by an integration by parts 
and is 
stated for future reference. 

\no{\bf Lemma 3.5}.  Assume $p, q: \R\rightarrow [0,\infty )$ with 
$p$ non-decreasing and $q$ non-increasing.  Then for any random 
variable $X$,
$ \P(p(X) q(X)) \leq \P(p(X))\ \P(q(X)).$    \qed
\med

\bigb {\bf 4. $L^2$ Bounds for the Clusters}

As we saw in the proof of Theorem 2.6, the key in establishing the linear 
growth of $h(A_n)$ is a good bound on the size of each cluster off the 
backbone. These clusters are governed by the laws $({\nu}_t, t\geq 0)$ of 
$C(t)$ conditioned to still be finite at time $t$ (Theorem 2.5).  The main 
results of this section are Theorems 4.4 and 4.6, which give 
uniform bounds on ${\nu}_t (\# (C))$ and ${\nu}_t (\# (C)^2)$ for $t\geq 0$.

\no {\bf Notation.} Set $H(t) = \P(T( \infty) > t)$ and  $G(t) = e^t H(t)$.


\no{\bf Lemma 4.1.} (a) The functions $G$ and $H$ satisfy
$$ \leqalignno{ H(t) &= e^{-t} \bigr( 1+\alpha \int^{t/\alpha}_0 
e^{\alpha u} H (u)^d du \bigl), \q t \geq 0, &(4.1) \cr
 G(t) &= 1+ \alpha \int^{t/\alpha}_0 
e^{-u(d-\alpha )} G(u)^d du, \q t \geq 0. &(4.2) \cr}$$
{(b)} $G(t)$ increases to a finite limit $c_{4.1} 
(\alpha , d)$ as $t \uparrow \infty$ and
$$ 1\leq c_{4.1} (\alpha, d) \leq 
\exp \big( \log (d/(d-\alpha ))/(1-\alpha ) \big). $$

\no{\it Proof}.  Condition on $U_0$ and use scaling to conclude that
$$ H(t) = e^{-t} + \int^t_0 e^{-s} H((t-s)/\alpha )^d ds.$$
\no
Setting $u = (t-s) / \alpha$ gives (4.1), and (4.2) is then immediate.  
Equation~(4.2) shows 
that $G(t)$ increases to a (possibly infinite) limit $c_{4.1} (\alpha , d) 
\geq 1$ as $t\rightarrow \infty $.  Lemma 2.3 shows that 
$H(t) \leq s (\alpha ) e^{-t}$ and therefore 
$H(t) \leq e^{{-(t-t_0)^+}}$ where $t_0 = \log (s (\alpha ))$.  
Hence from~(4.1),
$$ \eqalign{
H(t) &\leq e^{-t} + \alpha e^{-t} 
\int^{t_0} _0 e^{\alpha u} du + 1(t/\alpha > t_0)\alpha e^{-t} 
\int^{t/\alpha}_{t_0} e^{\alpha u} e^{-d(u-t_0)} du\cr
& \leq e^{\alpha t_0-t} + \alpha (d-\alpha )^{-1} 
e^{-t + dt_0-(d-\alpha )t_0} = e^{-(t-t_1)},}$$
\no where $t_1 = h(t_0) := \alpha t_0 + \log (d/(d-\alpha ))$. 
Iterating this procedure, we obtain $H(t) \leq e^{-(t-t_n)^+}$ for 
$n \geq 1$, where $t_{n+1} = h(t_n)$ and 
$\lim_{n \rightarrow \infty} t_n = \log (d/(d-\alpha) ) (1- \alpha )^{-
1}$.
Let $n\rightarrow\infty$ in $G(t) \leq e^{t_n}$ to complete the proof. 
\qed

\no {\bf Remark.}
(4.1) and (4.2) are rather nasty equations since they are non-linear and 
(worst of all) anticipative.  (See Athreya (1985) for some similar equations,  
arising from the distribution function of the random variable $\sup_{x \in 
\bB} T(x)$.) 
 Solutions to these equations are not unique because
$H\equiv 1$ also satisfies (4.1).  It is, however, not hard to show that
 $\P(T(\infty) >t)$ is the unique non-increasing solution $H$ to (4.1) for 
which $e^t H(t)$ is bounded.  
 Although it seems difficult to get sharp estimates from these 
 equations, in the next section some closely related equations will 
help us analyze the process $(W_n, n\in \Z_+)$, and in Section 8 we will 
derive some asymptotic results as $\alpha \uparrow 1$. In particular, the 
upper bound on $c_{4.1}$ in $(b)$ is by no means optimal -- see Remark~8.3. 

\no{\bf Notation}. 
If $n$ is a non-negative integer, define a 
Laplace transform with respect to the real variables $\theta \in
[0,1]$ and $\gamma \geq 1$ by
$$r(n,\theta , \gamma ) =  \P(\sum_{x\in \bB(n)} 
\exp (\theta \TT(n) - \gamma \alpha^{-n-1} (\TT^{(0)} (x) - 
\TT(n)))), \, n \in 
{\Z} _+ .$$
Let 
$$d_0 = \thalf d(d-1)-1.$$
For $\gamma > \alpha$, $n \in \N$ and $0\leq \delta \leq t$, let
$$ S (n, \gamma , \delta , t) = \sum_{x\in \bB(n)} 
1(t- \delta \leq T(x) \leq t) \exp (-\gamma
\alpha^{-n-1} (t-T(n))).$$

To get a feeling for $S$, set $\delta = t$, and integrate over $t$
to get
$$ \int_0^\infty S(n , \gamma , t , t) dt = (\alpha^{n+1} / \gamma)
   \sum_{x \in \B (n)} \exp (- \gamma \alpha^{-n-1} (T(x) - T(n))) .$$
This is $O(\alpha^{n+1})$, provided that not too many times
$T(x) - T(n)$ are near 0 on a scale of $\alpha^n$.  
The next result gives the critical technical estimate on $r$, which
is then used to show that $S$ is indeed $O(\alpha^n)$ when $\gamma 
\geq \alpha + 1$.  This will in turn lead to Theorems~4.4 and~4.6
on the respective $L^1$ and $L^2$ behaviours of the cluster size.
Recall that $s(\alpha ) = \prod^\infty_{k=1} (1-\alpha^k)^{-1}$.

\med {\bf Lemma 4.2.} For all $n \in \Z_+$, $\gamma \geq 1$,
$$ r(n,1,\gamma ) \leq c_{4.2} (\alpha , \gamma , d ) :=
d^2(d-1)^{-1} (1-\alpha ) s (\alpha )^2 \exp (d_0 /(\gamma (1-\alpha 
))) .$$

\no{\it Proof}. The result is trivial if $n=0$, thus 
assume $n \geq 1$ and fix $\gamma \geq 1$.  For the induction, assume 
the result for all $n' < n$ and all $\theta'$ (what is actually 
needed is $\theta' = \alpha^j \theta$).   
Fix $z_1 \in \{ 1 , \ldots , d \}$ and
set $\tau = \min \{ T^{(y)}(n-1) : y \in \B (1)
\setminus \{ z_1 \} \}$.  Using the symmetry of the tree when
computing $r(n, \theta , \gamma)$, we may sum 
over only those $x$ with $x_1 = z_1$ and then multiply by $d$.
This leads to 
$$\eqalign{
 r(n , \theta , \gamma) 
 &= d \cdot \P \Big( \sum_{x \in \B (n-1)} \exp \big(
   \theta \TT(n) - \gamma \alpha^{-n-1} (\alpha U_{z_1} + \alpha 
\TT^{(z_1)}
   (x) - \alpha T^{(z_1)}(n-1) ) \cr
&\qquad \qquad - \gamma \alpha^{-n-1} (\alpha T^{(z_1)}(n-1) - 
\TT(n) 
   ) \big) \Big) \cr
&= d \cdot \P \Big( \kern-7pt \sum_{x \in \B (n-1)} 
   \exp \big( \theta \TT(n) - \gamma
   \alpha^{-n} (\TT^{(z_1)} (x) - \TT^{(z_1)}(n-1)) \cr
 &\qquad \qquad - \gamma \alpha^{-n}
   (T^{(z_1)}(n-1) - \alpha^{-1} \TT(n) ) \big) \Big).}$$
Divide the above expectation into two terms corresponding to
the events $\{ T^{(z_1)}(n-1) > \tau\}$ and $\{ T^{(z_1)}(n-1) 
\leq \tau\}$.  Use the fact that on the latter event 
$\TT^{(z_1)}(n-1) = \alpha^{-1} \TT(n)$ while on the former event
$\tau = \alpha^{-1} \TT(n)$ to see that
$$ r(n, \theta , \gamma) = I_1 + I_2, \qquad \hbox{where}$$
$$\eqalign{
I_1 &= d \cdot \P \Big( 1(T^{(z_1)}(n-1) > \tau) \exp (\theta \alpha 
 \tau -  \gamma \alpha^{-n} (T^{(z_1)}(n-1) - \tau)) \cr
&\qquad \times \sum_{x \in \B (n-1)}
   \exp ( - \gamma \alpha^{-n} (\TT^{(z_1)} (x) - \TT^{(z_1)}(n-1))) 
\Big), \cr
I_2 &= d \cdot \P \Big( 1(T^{(z_1)}(n-1) \leq \tau) \exp (\theta 
 \alpha T^{(z_1)}(n-1)) \cr
& \qquad \times \sum_{x \in \B (n-1)} \exp ( - \gamma \alpha^{-n} 
   (\TT^{(z_1)} (x) - \TT^{(z_1)}(n-1))) \Big) \; . } \leqno(4.3)$$
The plan is to bound these terms by constant multiples of
$r(n-1 , \alpha \theta , \gamma)$ and then apply the induction 
hypothesis. (In fact $I_1$ will be of smaller order).  

Consider $I_1$ first.  The term $1(T^{(z_1)}(n-1) > \tau) \exp (\theta
\alpha \tau - \gamma \alpha^{-n} (T^{(z_1)}(n-1) - \tau))$ will
make this term relatively small as $n$ becomes large;
heuristically, the typical difference between $T(n)$ and $T(x)$
will be of order 1 when $x_1 \neq z_1$, and so the factor of
$\alpha^{-n}$ in the exponent makes these terms small.  To verify 
this,
fix $z_2 \in B(1) \setminus \{ z_1 \}$ and set $R = \TT^{(z_1)}(n-1) 
- \TT^{(z_2)}(n-1)$.  Focus on the expression in front of the sum in 
$I_1$
and integrate over the pair $(U_{z_1} , U_{z_2})$ to see that
$$\eqalign{
\P &\big( 1(T^{(z_1)}(n-1) > \tau) \exp (\theta \alpha \tau - 
   \gamma \alpha^{-n} (T^{(z_1)}(n-1) - \tau)) | \F_{(z_1 , \infty)} 
   \big) \cr
&= (d-1) \P \big( 1(U_{z_1} + \TT^{(z_1)}(n-1) > \tau = U_{z_2} + 
   \TT^{(z_2)}(n-1)) \cr
&\qquad  \times \exp (\theta \alpha (U_{z_2} + \TT^{(z_2)}(n-1)) 
   - \gamma \alpha^{-n} (U_{z_1} - U_{z_2} + R)) | \F_{(z_1 , \infty)}
   \big) \cr
&\leq (d-1) \P \Big( \big[ \int_0^\infty \int_0^\infty 1 (u_2
   - u_1 \leq R) \exp (\theta \alpha u_2 - \gamma 
   \alpha^{-n} (u_1 - u_2) ) e^{-u_1} \, du_1 \, e^{-u_2} \, du_2 \big] 
\cr
&\qquad \times \exp (\theta \alpha \TT^{(z_2)}(n-1) - \gamma 
\alpha^{-n}
   R) | \F_{(z_1 , \infty)} \Big) \; . }$$
First integrate $u_2$ over $(0 , u_1 + R)$ and then integrate
$u_1$ over $((-R)^+ , \infty)$ to bound the above by 
$$\eqalign{
(d-1) (\alpha \theta + \gamma \alpha^{-n} &- 1)^{-1} (2 - \alpha
   \theta)^{-1} \cr
     &\times \P \big( \exp (\theta \alpha \TT^{(z_2)}(n-1) 
  - (2 - \alpha \theta) (-R)^+ + (\alpha \theta - 1)R) |
   \F_{(z_1 , \infty)} \big) \cr
 = \; (d-1) \alpha^n (&\gamma + \alpha^{n+1} \theta - \alpha^n)^{-1}
   (2 - \alpha \theta)^{-1} \exp (\theta \alpha \TT^{(z_1)}(n-1)) \cr
     &\times \P \big( \exp ( -(2-\alpha \theta) (-R)^+ - R) | 
   \F_{(z_1 , \infty)} \big) \; . }$$
Note that $-(2-\alpha \theta) (-R)^+ - R \leq 0$ for $\alpha ,
\theta \in [0,1]$, and so this is bounded by 
$$(d-1) \alpha^n (\gamma - \alpha^n)^{-1} (2 - \alpha \theta)^{-1}
   \exp (\theta \alpha \TT^{(z_1)}(n-1)) \, .$$
Substitute this bound into the expression for $I_1$ (the summation
being $\F_{(z_1 , \infty)}$-measurable) to conclude that $I_1$ is at 
most
$$\eqalign{
     d (d-1) (2-&\alpha \theta)^{-1} \alpha^n (\gamma - \alpha^n)^{-1}\cr
    &\times P \big( \sum_{x \in \B (n-1)} \exp (\theta \alpha \TT^{(z_1)}(n-1) 
 - \gamma \alpha^{-n} (\TT^{(z_1)} (x) - \TT^{(z_1)}(n-1))) \big) \cr
 &= (d_0 + 1) (1 - \alpha \theta / 2)^{-1} \alpha^n (\gamma -
   \alpha^n)^{-1} r (n-1 , \alpha \theta , \gamma) \, .} \leqno(4.4) $$

Consider now $I_2$.  By symmetry, as $\theta$ becomes small the terms in front 
of the summation in the expression for $I_2$ should have mean close to 1, and 
so $I_2$ should be close to $r (n-1 , \alpha \theta , \gamma)$.  This is true, 
but to make the argument rigorous we must control the possible correlations 
between the summation and the remaining terms in the integrand. This makes use 
of the monotonicity results from Section~3, and in particular Corollary~3.4.  

Let $q(t) = \P ( \exp (\theta \alpha U_{z_1}) 1(U_{z_1} + t \leq \tau))$
and note that $q$ is decreasing in $t$.  Use the independence of 
$\F_{(z_1 , \infty)}$ and $\sE_{z_1} \supset \sigma (U_{z_1} , \tau)$
to see that
$$\eqalign{
 \P \big( 1(&T^{(z_1)}(n-1) \leq \tau) \exp (\theta \alpha T^{(z_1)}(n-
   1)) | \F_{(z_1 , \infty)} \big) \cr
&= \P \big( 1(U_{z_1} + \TT^{(z_1)}(n-1) \leq \tau) \exp ( \theta 
\alpha 
   (U_{z_1} + \TT^{(z_1)}(n-1))) | \TT^{(z_1)}(n-1) \big) \cr
&= \exp (\theta \alpha \TT^{(z_1)}(n-1)) q(\TT^{(z_1)}(n-1)) \, .} 
\leqno(4.5)$$
The joint independence of $\tau , U_{z_1}$ and $\TT^{(z_1)}(n-1)$
shows that $${\tilde q} (u) = \P (U_{z_1} + \TT^{(z_1)}(n-1) \leq \tau 
| U_{z_1} = u)$$ is decreasing in $u$.  The same independence and 
Lemma~3.5 give 
$$\eqalign{
\P (q(\TT^{(z_1)}(n-1))) &= \P (\exp (\theta \alpha U_{z_1}) {\tilde q}
   (U_{z_1})) \cr
&\leq \P (\exp (\theta \alpha U_{z_1})) \P ({\tilde q} (U_{z_1})) 
 = (1 - \theta \alpha)^{-1} d^{-1} } \leqno(4.6)$$
by symmetry.  Using~(4.5) in the expression for $I_2$ we get
$$  \eqalign{ I_2 = d \P \Big( q (\TT^{(z_1)}(&n-1)) \exp (\theta \alpha 
\TT^{(z_1)}(n-1))  \cr
  &\times \P \big( \kern-5pt \sum_{x \in \B (n-1)} \exp (-\gamma \alpha^{-n}
 (\TT^{(z_1)} (x) - \TT^{(z_1)}(n-1))) | \TT^{(z_1)}(n-1) \big) \Big) \,.}$$
Corollary~3.4 shows we may assume the conditional expectation is
a nondecreasing function of $\TT^{(z_1)}(n-1)$.  Recalling that
$q$ is nonincreasing, we again use Lemma~3.5 to conclude that
$$\eqalign{
I_2 &\leq d \P (q(\TT^{(z_1)}(n-1))) \P \big( \exp (\theta \alpha 
   \TT^{(z_1)}(n-1)) \cr
  & \qquad \times \sum_{x \in \B (n-1)} \exp (-\gamma \alpha^{-n}
   (\TT^{(z_1)} (x) - \TT^{(z_1)}(n-1))) \big) \cr
&\leq (1 - \theta \alpha)^{-1} r (n-1 , \theta \alpha , \gamma),}$$
by~(4.6).

Combine the above with (4.3) and (4.4) to see that
$$ \eqalign{
r(n,\theta ,\gamma ) 
  &\leq (1-\alpha\theta)^{-1}\left(1+(d_0+1)\alpha^n(\gamma
-\alpha^n)^{-1}\right)
     r(n-1,\alpha\theta ,\gamma )\cr
  &= (1-\alpha\theta )^{-1}{\left(
1-\alpha^n\gamma^{-1}\right)}^{-1}(1+d_0\gamma^{-1}
     \alpha^n)r(n-1,\alpha\theta ,\gamma ).\cr}$$
We now use induction, recalling that $\gamma \geq 1$, to conclude that
$$ \eqalign{
r(n,1,\gamma )
 &\leq \big[ \prod^{n-1}_{k=1} (1-\alpha^k)^{-1} \prod^n_{k=2} 
(1-\alpha^k\gamma^{-1})^{-1} \prod^n_{k=2} (1+d_0\gamma^{-
1}\alpha^k) \big] 
r(1,\alpha^{n-1}, \gamma )\cr
 & \leq s(\alpha )^2 (1-\alpha ) \exp (d_0 \gamma^{-1}(1-\alpha )^{-1})
 r(1,\alpha^{n-1}, \gamma )\ .}\leqno (4.7)$$
\no Note that
$$ 
r(1, \alpha^{n-1} , \gamma )  \leq d\P(\exp (\alpha^{n-1} \TT(1) ))
  \leq d \P(\exp (\min_{x \in \B(1)} U_x)) = d^2/(d-1),$$
and use this in (4.7) to complete the argument.\qed
\med
{\bf Lemma 4.3}. For $\gamma \geq 1 + \alpha$, $n\in \Z_+$, 
$0 \leq \delta \leq t$,
$$\P(S(n,\gamma , \delta, t)) \leq e^{-t} c_{4.2} 
(\alpha , \gamma - \alpha , d ) (\alpha^{n+1} \wedge \delta).$$

\no{\it Proof.}  Let $\gamma , \delta , t$ be as above, and let
$n \in \Z_+$; for $x \in \B (n)$ we evaluate the corresponding
summand in $S(n , \gamma , \delta , t)$.  First, integrate out
$U_0$ to see that
$$\eqalign{
\P \big( 1(t &-\delta \leq T(x) \leq t) \exp (-\gamma \alpha^{-n-1}
   (t - T(n))) \big) \cr
&= \, \P \big( 1(\TT^{(0)} (x) \leq t) \int_{(t - \delta - \TT^{(0)} 
(x))^+}^{t 
   - \TT^{(0)} (x)} \exp (-u -\gamma (t - u - \TT^{(0)} (x)) \alpha^{-n-
1}) \, du \cr
&\qquad \times \exp (-\gamma (\TT^{(0)} (x) - \TT(n)) \alpha^{-n-
1}) \big) \; .}$$
Change variables to $s = t - u - \TT^{(0)} (x)$ and bound the 
Lebesgue integral from above on the event $\{ \TT^{(0)} (x) \leq t \}$ by:
$$\eqalign{
 \int_0^{\delta \wedge (t - \TT^{(0)}(x))} &\exp (- (\gamma \alpha^{-
n-1} - 1) s) \, ds \exp (-(t - \TT^{(0)} (x))) \cr
\leq &\alpha^{n+1} (\gamma - \alpha^{n+1})^{-1} (1 - \exp (-(\gamma 
   \alpha^{-n-1} - 1) \delta)) \exp (-t + \TT^{(0)} (x)) \cr
\leq &e^{-t} \alpha^{n+1} (\gamma - \alpha^{n+1})^{-1} \min \{
   1 , (\gamma \alpha^{-n-1} - 1) \delta \} \exp (\TT^{(0)} (x)) \cr
\leq & e^{-t} \min \{ \alpha^{n+1} (\gamma - \alpha)^{-1} , \delta \}
   \exp (\TT^{(0)} (x)) \; .}$$
Thus
$$\eqalign{
\P \big( &1(t-\delta \leq T(x) \leq t) \exp (-\gamma \alpha^{-n-1}
   (t - T(n))) \big) \cr
&\leq e^{-t} \min \{ \alpha^{n+1} (\gamma - \alpha)^{-1} , \delta \}
   \P \big( \exp (\TT(n) - (\gamma \alpha^{-n-1} - 1) (\TT^{(0)} (x) - 
\TT(n)) ) \big) \cr
&\leq e^{-t} \min \{ \alpha^{n+1} (\gamma - \alpha)^{-1} , \delta \}
   \P \big( \exp (\TT(n) - (\gamma - \alpha) \alpha^{-n-1} (\TT^{(0)} 
(x) - \TT(n))) \big) \, . }$$

\no 
Summing over $x$ in $\bB(n)$ and recalling that $\gamma \geq 1 + \alpha$, 
 $$ \eqalignno{
\P(S(n,\gamma , \delta , t)) 
&\leq e^{-t} \min (\alpha^{n+1} (\gamma - \alpha )^{-1} , \delta ) r(n, 
1, \gamma - \alpha )\cr
 & \leq e^{-t} c_{4.2} (\alpha , \gamma - \alpha , d) (\alpha^{n+1} 
\wedge \delta ) \q \hbox{ (by Lemma 4.2)}. &\square \cr}$$

\no {\bf Notation.}  Let $c_{4.3}(\alpha,d)=c_{4.1}(\alpha, d)^d
c_{4.2}(\alpha, d - \alpha, d)(1 - \alpha)^{-1}.$

In the following result we set $C(0-)= \emptyset.$

\bigb {\bf Theorem 4.4}. (a) For all $0\leq \delta \leq t$ and $n$ in 
$\Z_+$
$$ \P(\# ((C(t) - C((t-\delta )-)) \cap \bB(n))|T(\infty) > t)
 \leq c_{4.1} (\alpha , d)^d c_{4.2} (\alpha, d - \alpha, d)
 (\alpha^{n+1} \wedge \delta).$$

\no {(b)} For all $0\leq \delta \leq t$,
$$\eqalign{ \P (  \# ( C (t) - &C ((t-\delta )-)) | T(\infty) > t)\cr
 & \leq c_{4.1} (\alpha , d)^d c_{4.2} (\alpha , d-\alpha , d) 
(1-\alpha )^{-1} (\log^+ (1/\delta ) +1)
(\delta \wedge 1)\cr
 &\leq c_{4.3} (\alpha , d), }$$
and so in particular, $\P(\# C(t) |T(\infty) > t) \leq c_{4.3}(\alpha,d)$ for all 
$t\geq 0$.

\item {(c)} $\dst \P(\# (C(t\wedge T(\infty) ) \cap \bB (n)))\leq 
c_{4.3} 
(\alpha , d) (1-e^{-t} ) \hbox{ for all } t \geq 0, n\in \Z_+.$ 
\med
{\it Proof}. Clearly $T(\infty) \leq T(n) + T(\minver_n, \infty )$ and
 if $z\in \bB(n)$
then Lemma 4.1 shows
$$ \P(T(z, \infty ) > t) = \P(T(\infty) > t\alpha^{-n-1} )^d \leq c_{4.1} 
(\alpha, d)^d \exp (-dt\alpha ^{-n-1}).\leqno (4.8) $$ 

\no Therefore,
$$ \eqalign{
  \P \bigl( \# ((C(t) - C((t &- \delta )-))\cap \bB (n) ) 
1(T(\infty) >t)\bigr)\cr
&= \sum_{x\in \bB (n)} \P(t-\delta \leq T(x) \leq t< T(\infty))\cr
&\leq \sum_{x\in \bB(n)} \P \bigl( 1(t-\delta \leq T(x) \leq t)\ 
\P(T(\minver_n, \infty ) > t-T(n) | {\cal F}_n ) \bigr) \cr 
 & \leq c_{4.1} (\alpha , d)^d \P(S(n,d, \delta , t)) \q \hbox{ (by\ 
(4.8))}.}$$ 

\no Lemma 4.3 therefore shows that
$$ \eqalign{
\P\bigl(\# ((C(t) - C((t&-\delta ) - )) \cap \bB (n)) | 
   T(\infty) > t \bigr)\cr
&\leq c_{4.1} (\alpha , d)^d c_{4.2} (\alpha , d - \alpha , 
d)(\alpha^{n+1} 
\wedge \delta ) (e^t \P(T(\infty) > t))^{-1}\cr 
 &\leq c_{4.1} (\alpha , d)^d c_{4.2} (\alpha , d - \alpha , 
d)(\alpha^{n+1} 
\wedge \delta ).}\leqno (4.9) $$ 

\no This proves (a), and (b) follows upon summing over $n\in \Z_+$.

\no (c)  As in (a) we have
$$ \eqalign{
\P ( \# (C(&t \wedge T(\infty) ) \cap \bB(n) )) \leq \sum_{x\in 
\bB(n)} 
\P(1(T(x) \leq t) \P(T(\minver_n, \infty ) > T (x) - T(n) | {\cal 
F}_n))\cr 
& \leq \sum_{x\in \bB(n)} \P \big( 1 (T(x) \leq t) c_{4.1} 
   (\alpha , d)^d \exp (-d (T(x) - T(n) ) \alpha^{-n-1}) \big) \, 
   \hbox{ (as in (4.8)).}}$$
Integrate out $U_0$ to see that the last summation equals
$$\eqalignno{
c_{4.1} (\alpha , d)^d \sum_{x\in \bB(n)} \P \Big( \exp (-d 
(\TT^{(0)} (x) 
   - &\TT(n))\alpha^{-n-1}) 1(\TT^{(0)} (x) \leq t) \big[ 1-\exp 
   (-(t-\TT^{(0)} (x))) \big] \Big) \cr 
 &\leq c_{4.1} (\alpha , d)^d (1-e^{-t}) r(n,0,d)\cr
 &\leq c_{4.1} (\alpha , d)^d c_{4.2} (\alpha , d, d)(1-e^{-t})\ \ 
\hbox{ (by\ Lemma \ 4.2)}\cr
 &\leq c_{4.3} (\alpha , d)(1-e^{-t}).  &\square}$$

\no{\bf Remark 4.5.}  If $x\in (0,1)$ and 
$N= [(\log 1/x) / (\log 1/\alpha) ]$ then
$$ \eqalign{
s(\alpha )  &\leq \Big \{ \prod^N_1 (1-\alpha )^{-1} (1+\ldots + 
\alpha^{k-1})^{-1} \Big \} \exp ( \sum^\infty_{k=N+1}
\alpha^k (1- \alpha^k)^{-1})\cr
 &\leq (1-\alpha )^{-N} \prod^N_{k=1} (kx)^{-1} 
\exp (x(1-x)^{-1} (1-\alpha )^{-1}).\cr}$$

\no Use Stirling's Formula and optimize over $x$ to see there are 
constants, 
$c_{4.4}, c_{4.5} > 0$ such that 
$$ s(\alpha ) \leq c_{4.4} \exp (c_{4.5} 
(1-\alpha )^{-1}) \hbox{ for all } \alpha \in (0,1). \leqno (4.10)$$ 
\no Using (4.10) and Lemma 4.1(b) it follows that
$$ \leqalignno{ c_{4.3} (\alpha , d)  &=  c_{4.1} 
( \alpha , d)^d c_{4.2} (\alpha , d- \alpha , d)(1 - \alpha ) ^{-1} \cr 
 & \leq c_{4.4}^2 d^2 (d-1)^{-1} \exp \big( (2c_{4.5} + (d_0 / (d-1))+ 
d\log (d/(d-1)))/(1-\alpha ) \big)  &(4.11) \cr
 & = c_{4.6}(d) \exp ( c_{4.7}(d)/(1-\alpha)) .}$$ 

The Central Limit Theorem in Section 6 will require the following $L^2$ bounds 
on the conditioned clusters.  

\med {\bf Theorem 4.6}.  There is a $c_{4.8} (\alpha , d) > 0$ such 
that
$$ \P(\# C(t)^2 | T(\infty) > t) \leq c_{4.8}(\alpha,d) \qquad \hbox{ for all } 
t\geq 0.$$
\no {\it Proof}. 
(In this proof, we suppress the dependence of $c_{4.i}$ on $(\alpha , d)$). 
Let $\sum^\prime$ denote summation over $x_1, x_2$ in 
$\bB$ 
for which $x_1 \wedge x_2$ is distinct from both $x_1$ and $x_2$.  
Then 
 $$ \leqalignno{ \P(\# C(t)^2 1 (T(\infty) > t))  
\leq 2 \P(&\sum_{x_1} \sum_{x_2} 1(x_2 \geq x_1) 1(T(x_2)\leq t< 
T(\infty) ))
&(4.12)\cr 
&+ \P(\pain 1(T(x_1) \vee T(x_2) \leq t < T(\infty) )).}$$ 
\no If $y\not= 0$ let $\widehat T(y) = \inf \{ T(x) : |x| = |y| , 
x\not= y\}$ and let $\widehat Z (y)$ be the a.s. unique vertex in 
$\bB 
(|y|) - \{ y\}$ at which this minimum time is attained.  Define 
$\widehat 
T (0) = \infty$ and $\widehat Z (0) = 0$.  Note that $(\widehat T (y), 
\widehat Z (y))$ is ${\cal F}_{|y|}$-measurable and $T(\widehat Z (y), 
\infty)$ is $\sE_y$-measurable (recall from the beginning of Section 
2 
that $\sE_y$ is the``wide-sense" past up to $y$).  In the second 
summation 
in (4.12), write $x_1=v \oplus i \oplus y_1$ and $x_2 = v\oplus 
j\oplus y_2$ 
where $v=x_1 \wedge x_2$ and $i\not= j \in \{ 1, \ldots, d\}$ to get 
$$ \eqalign{ 
\P (\pain 1 (&T(x_1) \vee T (x_2) \leq t < T(\infty)))\cr
  \leq \P \Big\{ &\sum _v 1 (T(v) \leq t < \widehat T (v) + 
T(\widehat Z (v), 
     \infty ))\cr 
  &\times \P\Bigl[ \sum_{1 \leq i} \sum_{\not= j \leq d}
 \Bigl( \sum _{y_1} 1 (T (v) + T^{(v\oplus i)} (y_1) \alpha^{|v|+1} 
\leq t< T(v) + T^{(v\oplus i)}(\infty) \alpha^{|v|+1}) \Bigr)\cr 
 & \qquad \times \Bigl( \sum_{y_2} 1 (T(v) + T^{(v \oplus j)} (y_2) 
\alpha^{|v|+1} \leq t < T(v) + T^{(v \oplus j)}( \infty) 
\alpha^{|v|+1})\Bigr) \bigm| \sE_v \Bigl]  \Big\}. }$$
\no
$(T^{(v\oplus i)} (\cdot ), T^{(v\oplus j)} (\cdot ))$, $i\not= j$, are 
independent and are jointly independent of $\sE_v$.  The above 
therefore 
equals 
$$\eqalign{
d(d&-1) \int \sum_v 1 (T(v)(\omega) \leq t) \P \big( T (\widehat Z 
(v) ,
   \infty) > t - \widehat T (v) | \F_{|v|} \big) (\omega) \cr
&\qquad \times \big[ \P ( \# C((t - T(v)(\omega)) \alpha^{-|v| - 1}) 
   1(T(\infty) > (t - T(v)(\omega)) \alpha^{-|v| - 1})) \big]^2 \, 
   d\P (\omega) \cr
&\leq \; d(d-1) \int \sum_v 1(T(v)(\omega) \leq t) \P (T(0,\infty) >
   (t - \widehat T (v)(\omega)) \alpha^{-|v|}) \cr
&\qquad \times c_{4.3}^2 [ \P (T(\infty) > (t - T(v)(\omega)) 
\alpha^{-
   |v| - 1}) ]^2 \, d\P (\omega) \, .} \leqno (4.13)$$
The time $T(0,\infty)$ is equal in law to the minimum of $d$ 
independent copies of $\alpha T(\infty)$.  Lemma~4.1(c) and 
the fact 
that $T(|v|) = \min \{ \widehat T (v) , T(v) \}$ therefore show that 
if $T (v) (\omega) \leq t$, then
$$\eqalign{
\P (T(0,\infty) &> (t - \widehat T(v) (\omega)) \alpha^{-|v|})
   [\P (T(\infty) > (t - T(v)(\omega)) \alpha^{-|v| - 1})]^2 \cr
\leq &c_{4.1}^{d+2} \exp \big\{ -d \alpha^{-|v|-1} (t - \widehat T (v)
   (\omega))^+ - 2 \alpha^{-|v|-1} (t - T(v)(\omega)) \big\} \cr
\leq &c_{4.1}^{d+2} \exp \big\{ - 2 \alpha^{-|v|-1} (t - T(|v|)
   (\omega)) \big\} \, .}$$
Substitute this into the RHS of~(4.13), thus bounding the second term 
on
the RHS of~(4.12) by
$$d(d-1) c_{4.3}^2 c_{4.1}^{d+2} \int \sum_v 1( T(v)(\omega) \leq t)
   \exp \big\{ - 2 \alpha^{-|v|-1} (t - T(|v|) 
   (\omega)) \big\} d\P (\omega )$$
which is equal to 
$$d(d-1) c_{4.3}^2 c_{4.1}^{d+2} \P ( \sum_{n=0}^\infty S (n,2,t,t))$$
and hence at most $c(\alpha , d) e^{-t}$ by Lemma~4.3.

The first term on the right side of (4.12) equals
$$ \eqalign{ 2 \P(\sum_x |x| 1 (T(x) \leq t < T(\infty)))
&= 2\P (\sum^\infty_{n=1} n \# (C(t) \cap \bB(n)) |T(\infty) > t) 
   \P(T(\infty) > t)\cr 
&\leq c^\prime (\alpha, d) \P(T(\infty) > t)\ \qquad 
\hbox{ (by Theorem 4.4(a))}.}$$

\no Use this and the above bound in (4.12) to conclude
$$\eqalignno{ \P(\# C(t)^2 | T(\infty) > t) &\leq c^\prime (\alpha , d) + 
c(\alpha , d) (e^t \P(T(\infty) > t))^{-1}\cr
 &\leq c^\prime (\alpha , d) + c(\alpha , d)  \; . &\square}$$

\bigb { \bf 5. The Process $\W$}

Recall from Section 2 that 
$W_n = T(\spine_n , \infty ) \alpha^{-n-1}$ is the rescaled time to 
percolate along the backbone from generation $n$ to $\infty $, and 
$\W=\{W_n: n \in {\Z}_+ \}.$  In this 
section we will show $\W$ is an ergodic Markov chain which 
stochastically 
decreases to its unique invariant measure.  Moreover there is an 
exponentially 
fast coupling mechanism for the chain. 

\no{\bf Notation.}
For each $x\in \bB$, set
$$W (x) = T(x,\infty ) \alpha^{-|x|-1} = \alpha^{-1} \TT^{(x)}(\infty) \, 
.$$ 
Let $F$ denote the c.d.f.\ for $W_0$.
Then $W(x)$ is equal in law to $W_0$ for each $x \in \B$; 
since $W_0$ is the minimum of $d$ independent copies 
of $T(\infty)$, Lemma 4.1 shows that 
$$ 1-F(t) = \P(W_0 > t) \leq c_{4.1} (\alpha , d)^de^{-dt}.\leqno 
(5.1)$$ 
\no Hence we can define a law on $(0, \infty )$, which ``tilts'' $W_0$
to the right by $\alpha$:  
$$\Lambda (A) = \P ( e^{\alpha W_0} 1( W_0 \in A)) / \P(e^{\alpha
W_0} ) \, .$$
\no Write $\Lambda (t) = \Lambda ((0,t)).$

To define the filtration on which the sequence $\W$ is Markov, 
we introduce two more pieces of notation. 

\no{\bf Notation.}
Let $\bar{\cal F}_x = {\sE}_x \vee \sigma (W(x))$. 
Then $x\leq y$ implies $\bar{\cal F}_x \subset \bar{\cal F}_y$
and if $x\in \bB(n)$, then
$$ \{ \spine_n = x\} = \bigl\{ T(x) + \alpha^{n+1} W(x) < \min 
\{T(x^\prime )+
\alpha^{n+1} W(x^\prime ) : x^\prime \in \bB(n) -
\{x \}\} \bigr\} \in \bar{\cal F}_x.$$
\no Therefore $\spine_n$ is an ``$\bar{\cal F}_x$-stopping point" and 
we may 
define a filtration $\{ {\cal W}_n\}$ by
$$ {\cal W}_n = \bar{\cal F}_{\spine_n} = \{ A\in {\cal F} : 
A \cap \{\spine_n = x\} 
\in \bar {\cal F}_x\ \hbox{ for all } x \in \bB\}.$$
The verification that $\{ W_n \}$ is adapted to $\{ {\cal W}_n \}$ is 
immediate.
\bigb
\no{\bf Theorem 5.1}. (a)  $\W = \{ W_n : n\in \Z_+ \}$ is a 
$({\cal W }_n)$-Markov chain such that for Borel sets 
$B \subseteq (0,\infty)$,
$$ \P (W_{n+1} \in B | {\cal W}_n ) = p (B |W_n ) :=
\Lambda (B| (0, W_n / \alpha ))\ \, .$$
(b) $(W_n, n\in \Z_+)$ is stochastically non-increasing.
\med
\no{\it Proof}.  (a) Let $A\in {\cal W}_n, x\in \bB(n+1)$ and 
$B$ a Borel subset of $(0,\infty)$.  Then 
$$\leqalignno{  \P(W_{n+1} &\in B, A, \spine_{n+1} = x) &(5.2) \cr
  = \P \Bigl[ &1(A, \spine_n = x|n ) \P (W(x) \in B, \cr
  &U_x + \alpha W(x) 
< \min_{1 \leq i \leq d-1} U_{\sib_i(x)} + \alpha W(\sib_i(x))  
\bigm| \bar{\cal F}_{x|n} ) \Bigr].\cr} $$
Both $W(x|n )$ and the event in the conditional probability are
${\cal F}_{(x|n , \infty )}$-measurable.  Since ${\cal F}_{(x|n,\infty )}$ 
and $\sE_{x|n}$ are independent, while $\bar {\cal F}_{x|n} = 
\sE_{x|n} 
\vee \sigma (W(x |n ))$, it follows that~(5.2) equals
$$ \eqalign{
\P \Big[ 1(A, \spine_n & = x|n )  \P(W(x) \in B, \cr
 &U_x + \alpha W(x) 
< \min_{1 \leq i \leq d-1 } (U_{\sib_i(x)} + \alpha W(\sib_i(x)))
\bigm| W (x|n)) \Big] .} \leqno~(5.3)$$

\no Let $\bB(1) = \{z_1, \ldots , z_d\}$, and set
$$ \eqalign{ 
q(\cdot \, |w) &= d \cdot \P(W(z_d) \in \cdot , U_{z_d} + \alpha W(z_d) 
< \min_{1 \leq i \leq d-1 } (U_{z_i} + \alpha W(z_i))   \cr
&  \qquad \qquad \qquad \qquad \qquad \qquad \qquad \qquad
\bigm| \min_{1 \leq i \leq d } (U_{z_i} + \alpha W(z_i)) = w) \, .}$$

\no In other words, $q(\cdot \, | w)$ is a regular conditional 
probability for the right side.  The collections
$(W(x), U_x$, $((W(\sib_i (x))$, $U_{\sib_i(x)})$ : $i<d))$ and 
$(W(z_d), U_{z_d}, ((W(z_i), U_{z_i}) : i<d))$ are equal in law, and
$$W(x|n) = \min \{ U_x + \alpha W(x), U_{\sib_i(x)} + \alpha 
W(\sib_i(x)) : i< d\}.$$ 
 Therefore, using this in~(5.3), we deduce that 
$$ \eqalign{
\P(W_{n+1} \in B, A, \spine_{n+1} = x) &= d^{-1} 
\P \bigl(1(A , \spine_n =x|n) q(B|W(x|n))\bigr)\cr
 & =d^{-1} \P\bigl(1(A, \spine_n = x|n) q(B|W_n)\bigr).}$$
\no Sum over $x \in \bB(n+1)$ to conclude
$$\P(W_{n+1} \in B \, | \, \sW_n) = q(B| W_n)\ \hbox{\ 
a.s.}\leqno(5.4)$$
 %

If $h(y) = \P( \min_{1\leq i \leq d-1}(U_{z_i} + \alpha W(z_i)) > y)$ then 
clearly for $A\in \sB((0,\infty))$,
$$ \eqalign{
d\P\bigl(W (z_d) &\in B, 
\min_{1\leq i \leq d}(U_{z_i} + \alpha W(z_i))
\in A, U_{z_d}  + \alpha W(z_d) < 
\min_{1\leq i \leq d-1}(U_{z_i} + \alpha W(z_i))
\bigr)\cr
& = d\int\int 1(w\in B, u+ \alpha w \in A) h(u+ \alpha w )
e^{-u} du dF(w )\cr
& = d\int\int 1(w \in B,  y\in A) 1(w\leq y/\alpha )
h(y) e^{-y+ \alpha w } dF(w ) dy\cr
&= \int_A \psi (y) \Lambda (B| (0, y/\alpha ]) dy,}$$
\no where $\psi (y) = dh(y)e^{-y} \Lambda (y/\alpha )$.  Take 
$B=(0, \infty )$ to see that $\psi$ is the density of \break
$\min_{i \leq d} (U_{z_i} + \alpha W(z_i))$.  It follows that 
$\Lambda (\cdot \, |( 0, y/\alpha ])$ is a version of $q(\cdot \, |y)$ 
and (a) is a consequence of~(5.4). 

 \no (b) We have
$$ \eqalign{ \P(W_1 \geq t) &= d \P(1 (W (z_d)\geq t) \P(U_{z_d} + 
\alpha W(z_d ) 
 < \min_{1\leq i < d} (U_{z_i} + \alpha W(z_i)) |W (z_d)))\cr
 &\leq d\P(W(z_d) \geq t) \P(U_{z_d} + \alpha W(z_d) 
< \min_{1\leq i < d} (U_{z_i} + \alpha W(z_i)))\ \ \hbox{ (Lemma \ 
3.5)}\cr
 &= \P(W(z_d) \geq t)\ \ \hbox{ (symmetry)}\cr
 & = \P(W_0 \geq t).}$$

\no This proves $W_1 \prec W_0$.  The result now follows by 
induction and the 
fact that $p(\cdot \, |w)$ is stochastically non-decreasing in $w$. 
\qed

Since $W_n$ decreases stochastically in $n$, the laws converge 
weakly,
and we shall now show that the limit is nontrivial. 
The other goals 
for the remainder of this section are to understand the stationary 
distribution for $\{ W_n \}$ (including existence
and convergence to stationarity) and to prove an 
exponential rate of convergence through a coupling mechanism.
The first goal is achieved in the next result.

\no{\bf Theorem 5.2}.   
(a) $ \W = \{ W_n : n\in \Z_+ \}$ converges
to a stationary distribution $\pi$ on $(0, \infty )$.

\no(b) $\pi$ is the unique stationary distribution for $\W$ and the
stationary chain is ergodic.

\no(c) There are positive constants $c_{5.1}(p) (p>0)$, and $c_{5.2} $ 
such that
for all $w\geq 0$ and for all $p>0$
$$ \eqalign{
 &\pi((0, w]) \leq c_{5.1}(p) w^p \cr
 & \pi ([w , \infty )) \leq c_{5.2} e^{-dw}.} \leqno~(5.5)$$


We begin the proof with an integral 
equation that leads to an important bound~\kern-1pt(5.9) on the 
left tail of $F$.  
Let $F_1 (t) = \P(T(\infty) \leq t)$ 
and $g(x) = 1-(1-x)^d$.  $W_0$ is the minimum of $d$
independent copies of $T(\infty)$ and so
$$ F(t) = g(F_1 (t)) \in [F_1 (t), dF_1 (t)].\leqno~(5.6)$$

\no By conditioning on $U_0$ we also have
$$ F_1 (t) = \alpha e^{-t} \int^{t/\alpha}_0 e^{\alpha s} F(s) ds 
\leq F(t/\alpha ) (1-e^{-t}).\leqno~(5.7)$$

\no Equations~(5.6) and~(5.7), together with the easy fact
that $F_1 (t) > 0$ for positive $t$, imply that $F$ and $F_1$ 
are infinitely differentiable and satisfy
$$ F^\prime_1 (t) = -\alpha e^{-t} \int^{t/\alpha }_0 e^{\alpha s} F(s) 
ds 
+ F(t/\alpha ), \qquad 
F^\prime (t) = g^\prime (F_1 (t)) F^\prime_1(t), \leqno~(5.8)$$
$$ 0< F(t) / F(t/\alpha ) \leq d(1-e^{-t}) \leq td \quad 
\hbox{ for all } t>0.\leqno~(5.9)$$

\no{Next}, we let $\beta$ be the midpoint between $1$ and 
$\alpha^{-1}$, 
and show that when $W_n$ is small, the probability is at least
$1/4$ that $W_{n+1} \geq W_n \beta.$

\no{\bf Lemma 5.3}. (a) There is a $t_0>0$ such that 
$F(t\beta )/F (t/\alpha ) < {2\over 3}$ and $\Lambda (t\beta)/ 
\Lambda(t/\alpha )<~{3\over 4}$ for $t \in (0,t_0)$.

\no(b) For all $0<w \leq 1$, $\Lambda (w)/\Lambda (w/\alpha ) \leq 
dew.$

\no {\it Proof}.  Differentiate in~(5.8) to see
$$F^{\prime\prime} (t) = g^{\prime\prime} (F_1 (t)) F_1^\prime (t)^2
 + g^\prime (F_1(t)) F_1^{\prime\prime}
(t),\leqno(5.10)$$
and as $t\downarrow 0$,
$$ \eqalign{
   F^{\prime\prime}_1 (&t) = o(F (t/\alpha )) - F(t/\alpha ) 
+ F^\prime (t/\alpha ) \alpha ^{-1}\cr
 &= o(F(t/\alpha )) - F(t/\alpha ) +g^\prime (F_1 (t/\alpha )) 
F_1^\prime (t/\alpha ) \alpha^{-1}\cr
 & = o (F(t/\alpha )) - F(t/\alpha ) + g^\prime (F_1 (t/\alpha 
))\alpha^{-1} 
(o (F(t/\alpha^2))+F(t/\alpha^2))\ \ \hbox{\
(from\ (5.8))}\cr
 &= o(F(t/\alpha^2 )) - F(t/\alpha ) + d\alpha^{-1} F(t/\alpha^2).}$$
Use the latter in~(5.10) to conclude that as $t\downarrow 0$
$$ \eqalign{
F^{\prime\prime} (t)  & = o(F(t/\alpha )) + d(o (F(t/\alpha^2)) - 
F(t/\alpha ) 
+ d\alpha^{-1} F(t/\alpha^2))\ \ \hbox{ ((5.8)\  again)}\cr
 & = o(F(t/\alpha^2)) + d^2 \alpha^{-1} F(t/\alpha^2)-dF(t/\alpha )\cr
 & >0  \q \hbox{ for small $t$}. }$$
Therefore $F$ is convex near $0$ and for $t \in (0, t_0)$ 
we have
$$ \eqalignno{ (F(t/\alpha) - F(t\beta )) / F(t/\alpha ) \geq 
\thalf (F(t/\alpha ) -F (t) )/F (t/\alpha )
 &\geq \thalf (1-dt)  \q \hbox{ (by~(5.9))}\cr
 & > \third \q \hbox{ (for $t_0$ small enough)}.  }$$
For the second inequality in (a) note that for $t \in (0,t_0),$
$$\Lambda(t\beta )/\Lambda(t/\alpha )\leq F(t\beta 
)e^t/F(t/\alpha)< \threequarters, $$
where we have taken $t_0$ sufficiently small for the last inequality.
This proves (a), and (b) is a trivial consequence of~(5.9), just as 
above. \qed

\no {\it Proof of Theorem 5.2}.  (a), (c). The stochastic monotonicity 
in Theorem 5.1 shows
 that $(W_n)$ converges in distribution to a law $\pi$ on $[0, \infty 
)$.  
We must show $\pi (\{ 0\} ) = 0$.   If $t_0 >0$ is as in Lemma 5.3, 
then for 
$w \in (0, t_0 \wedge 1]$, 
$$ \eqalign{
\P(W_n \leq w)&\leq \P (W_{n+1} \leq w )\cr
 & = \P(\Lambda ((W_n/\alpha ) \wedge w )  \Lambda (W_n/\alpha 
)^{-1} )\cr
 &\leq \P(W_n \leq \beta^{-1}w ) + \P(1(\beta^{-1} w 
< W_n \leq w) \Lambda (\beta W_n ) \Lambda
(W_n / \alpha )^{-1})\cr
 &\qquad\qquad + \P (1 (W_n > w) \Lambda (w ) \Lambda (w 
/\alpha )^{-1} )\cr
 & \leq \P(W_n \leq \beta^{-1} w ) + \threequarters \P(\beta^{-1} w 
< W_n \leq w ) + dew P(W_n >w ),\cr}$$
by Lemma 5.3. Rearrange terms to conclude from the above that
$$ \P (\beta^{-1} w < W_n \leq w ) \leq 4 dew = cw \ \ \hbox{ for}\ 
w \in (0,
t_0 \wedge 1].$$
 \no
Iterating the above we find that for $k\geq k_0$,
$$ \P(0 < W_n \leq \beta^{-k} ) = \sum_{j=k}^\infty c \beta^{-j} =
 c(1- \beta^{-1})^{-1} \beta^{-k},$$
and therefore
$$ \P(0 \leq W_n \leq w) \leq c\beta (1 - \beta^{-1})^{-1} w 
= c_{5.1} w\ \q  \hbox{ for } 0<w < w_0. \leqno(5.11)$$

\no Let $n\rightarrow\infty$ to get the first inequality in (5.5) 
for $p=1$,  first for $w$ 
a continuity point of $\pi ((0, w ])$ in $(0,w_0)$, and
then (increase $c_{5.1}$ if necessary) for all $w >0$.  This proves 
the first part of $(c)$ for $p=1$ and also 
shows $\pi (\{ 0\} ) = 0$. Since $F$ is atomless, so is $\Lambda$, and 
$p(\cdot \, |w )$ is weakly continuous in $w $.  Taking limits in 
$\P(W_{n+1} 
\in \cdot ) = \P( p (\cdot \, | W_n))$, we see that $\pi$ is a 
stationary law.  
The upper bound on $\pi ([w, \infty ))$ is immediate from (5.1) and 
$\pi \prec 
F$ (again increase $c_{5.1}$ if need be). 

To prove the upper bound on the left tail of $\pi $ for general $p$ 
we consider $p=2$. The induction argument which will give the result 
for general $p$ in ${\N}$ will then be clear. The fact that 
$\pi $ is stationary implies that for $w\in (0,t_0)$,
$$
\eqalign{
\pi \left( [0,w]\right) 
  &= \int\limits_0^\infty {{\Lambda \left( {x\over \alpha }\wedge
w\right)}\over
      {\Lambda \left( {x\over\alpha }\right)}}\; d\pi (x)\cr
  &\leq \pi \left( [0,\beta^{-1}w]\right) +\int 1(\beta^{-1}w<x\leq w)\; 
      {\Lambda (w)\over \Lambda ({x\over\alpha })}\; 
      d\pi (x)\cr
  &\qquad + \int 1(w<x\leq {w\over\alpha} )\; {\Lambda (w)\over \Lambda
({x\over\alpha })}\; d\pi (x) + \int
      1({w\over\alpha }<x)\; {\Lambda (w)\over \Lambda ({x\over\alpha })}\;
d\pi (x),\cr}$$
and therefore,
$$\eqalign{
\pi \left( (\beta^{-1}w,w]\right) 
  &\leq {\Lambda (w)\over \Lambda (w\beta^{-1}\alpha^{-1})}\; \pi \left(
(\beta^{-1}w,w]\right) +
    {\Lambda (w)\over \Lambda ({w\over\alpha })}\; \pi \left(
(w,{w\over\alpha }]\right) + 
    {\Lambda (w)\over \Lambda ({w\over\alpha^2})}\cr
  &\leq \threequarters \; \pi \left(
(\beta^{-1}w,w]\right)+(dew)c_{5.1}(1){w\over\alpha }+(de)^2\;
{w^2\over\alpha }, \cr} $$
 using Lemma 5.3 and the $p = 1$ case.
The above implies
$$
\pi \left( (\beta^{-1} w,w]\right) \leq 4\left( dec_{5.1}(1)+(de)^2\right)
\alpha^{-1} w^2
$$
and the result for $p=2$ now follows as for $p=1$ in the above.
 
\no (b) (5.8) and (5.9) show that 
$F^\prime_1 (t) \geq e^{-t} F(t/\alpha ) >0$ for all $t> 0$. This and the 
second part of (5.8) imply $F^\prime (t) > 0$ 
and hence $\Lambda^\prime (t) > 0$ for all $t> 0$. Therefore $p(\cdot 
\, |w)$ 
has a strictly positive continous density on $(0,w /\alpha ]$.  It 
follows 
easily that $(W_n)$ is an indecomposable Markov chain and (b) is a 
consequence 
of (a) and Theorem 7.16 of Breiman (1968).\qed 

There is a natural coupling technique for $\W$ which, in addition to
 refining the above convergence result,
will also play an important role in the limit theorems of Section 6.  
Let 
$D=\{ (w_1, w_2 ) \ab: w_1 \geq w_2 >0\}$ and define a Markov 
kernel on 
$D$ by (set ${0/ 0} \equiv 0)$
$$ \eqalign{
\bar p (A_1 \times A_2 | (w_1, w_2)) = 
&\int^{w_2/\alpha}_0 1_{A_1\times A_2} (x, x) d\Lambda (x) 
\Lambda (w_2
/\alpha )^{-1}(\Lambda (w_2 /\alpha ) /\Lambda (w_1 /\alpha ))\cr
 & +\int^{w_1/\alpha}_{w_2/\alpha} 1_{A_1} (x_1) d \Lambda (x_1) 
(\Lambda (w_1 /\alpha ) - \Lambda (w_2/\alpha ))^{-1}\cr
 & \times \int^{w_2/\alpha}_0 1_{A_2} (x_2)  d\Lambda (x_2) 
\Lambda (w_2 /\alpha )^{-1}(1 - (\Lambda (w_2 /\alpha ) 
/\Lambda 
(w_1 /\alpha ))).}$$

\no
If $\bar p_i (\cdot \, | (w_1, w_2))$ is the $i^{th}$ marginal of 
$\bar p (\cdot \, | (w_1, w_2 ))\ (i=1, 2)$, then

$$ \eqalign{
\bar p_1 (A_1 | (w_1, w_2)) &= \int^{w_2/\alpha}_0 1_{A_1} (x) 
d\Lambda (x) \Lambda (w_1 /\alpha )^{-1} +
\int^{w_1/\alpha}_{w_2/\alpha} 1_{A_1} (x) d \Lambda (x) 
\Lambda (w_1 /\alpha )^{-1}\cr
 & = p(A_1 | w_1),}$$

\no
and similarly one sees that $\bar p_2 (A_2 | (w_1, w_2)) = p(A_2| 
w_2)$.  
We extend this coupling of $\W$ 
to a Markov kernel on $(0, \infty)^2$ by setting

$$ \bar p (A_1 \times A_2| (w_1, w_2 )) = \bar p (A_2 \times A_1 
|(w_2, w_1))\ \ \hbox{
if}\ w_2 > w_1.$$
Then $\bar p$ is a Markov kernel on $(0, \infty )^2$ with marginals 
$\bar p_i 
(\cdot \, |( w_1, w_2)) = p (\cdot \, |w_i)$ and so the induced Markov 
chain 
$((W_n^1, W^2_n), n\in \Z_+)$ will be a coupling of $(W_n)$. 
Clearly 
$W_0^1 \geq W^2_0$ (respectively, $W^1_0 \leq W^2_0$) implies 
$W^1_n \geq 
W^2_n$  (respectively $W^1_n \leq W^2_n$) for all $n\geq 0$
a.s., and if $\tcpl = \min \{ n : W^1_n=W^2_n\}$ then
$W^1_n = W^2_n$ for all $n \geq \tcpl $ a.s.

If $W^1_0 \geq W^2_0$, the chains will couple at the first time $n$ 
for which 
$W^1_n \leq W_{n-1}^2 /\alpha$ and hence it is possible for $W^2_n$ 
to jump 
onto $W^1_n$.  Unfortunately if $W^2_{n-1}$ is small the probability 
of $W^1_n 
\leq W^2_{n-1}/\alpha $ will be small, and to get a good coupling rate 
we must 
bound the time spent by $W_n^2$ in $(0, \delta ]$ for $\delta$ small. 
We start
with a stochastic lower bound on the left tail of the Markov kernel of 
$\W$.

\no{\bf Lemma 5.4}. There is a $w_0 \in (0,1]$ and a probability 
$\rho$ on $(0,\alpha^{-1}]$ such that
$$ \leqalignno{
 p ((0,w x] \,| w ) &\leq \rho ((0,x]) \quad \hbox{ for all }
 x \geq 0 \hbox{ and } w \in (0,w_0], &(5.12) \cr
 \int \log x d \rho (x) &= m >0, &(5.13) \cr
 \rho ((0,x]) &\leq w_0 dex \leq x \quad \hbox{ for all } 
x< \alpha^2.  &(5.14)\cr} $$

\no {\bf Proof}.  Let $w_1 = (2de)^{-1}$ and define $\rho$
 on $(0, w_1]\times (0,1)$ by
$$ \rho (w, x) = \prod^{n-1}_{k=0} (de\alpha^k w ) 
=  ( dew )^n \alpha^{n(n-1)/2} \quad \hbox{ if} \quad
\alpha^n \leq x < \alpha^{n-1}, n\in \N.$$

\no Then $\rho$ is non-decreasing in each variable, $\rho \leq 1/2$, 
$\rho (w , \cdot )$ is right-continuous, and

\no $\rho (w , 0+) = 0\equiv \rho (w , 0)$. Let $t_0$ be as in Lemma 5.3,
$0< w_0 \leq w_1 \wedge t_0$, and define

$$ \rho (x) = \cases{
\rho (w_0, x)&if $x<1$,\cr
{3\over 4} &if $1\leq x < \beta$,\cr
1 &if $\beta\leq x\leq \alpha^{-1}$.\cr}$$

\no
$\rho$ is the distribution function of a law (also denoted by $\rho$)
 on $(0, \alpha^{-1}]$.  If $(x, w ) \in (0,1) \times (0, w_0]$ 
and $n\in \N$ satisfies 
$\alpha^n \leq x < \alpha^{n-1}$, then
$$ \eqalign{ \Lambda (xw) \leq \Lambda (\alpha^{n-1} w)
 &\leq \prod^{n-1}_{k=0} (de \alpha^k w ) \Lambda
(\alpha^{-1} w )\ \ (\hbox{ by Lemma 5.3(b)}\ )\cr
 &= \rho (w , x) \Lambda (w /\alpha) \leq \rho (x) \Lambda (w 
/\alpha).}$$
Therefore
$$ p ((0, w x] | w ) = \Lambda (x w ) / \Lambda (w /\alpha ) \leq 
\rho (x)\ \hbox{ for}\ (x,w) \in (0,1) \times (0, w_0].$$

\no For $(x, w ) \in [1, \beta ) \times (0,w_0]$ we have, from Lemma 5.3(a),

$$ p ((0, w x] | w ) = \Lambda (x w ) /\Lambda (w/\alpha ) \leq 
\Lambda (\beta w ) / \Lambda (w /\alpha ) \leq \threequarters =\rho (x),$$
which proves (5.12).

If $\alpha^n \leq x < \alpha^{n-1}$ for some $n\geq 3$, then
$$ \eqalign{
\rho(x) = \rho (w_0, x) & = (dew_0)^n\ \alpha^{n(n-1)/2}\cr
 & \leq (dew_0) x^{(n-1)/2} \leq (dew_0)x,\cr}$$
which proves (5.14).

Finally,
$$ \eqalign{
\int^{1/\alpha}_0 \log x d\rho (x)
 &\geq \int^{\alpha^2}_0 (\log x) (dew_0) dx + (\log \alpha^2)
 \rho (w_{0}, 1-) + (\log \beta ) /4\cr
 & = (dew_0) \int^{\alpha^2}_0 (\log x) dx + (\log \alpha^2 )
 dew_0 + (\log \beta) /4,}$$
\no which is positive if we choose $w_0$ sufficiently small. \qed

For $w >0$, let $\P_w$ be the law of $(W_n, n\in \Z_+)$ starting at $W_0=w$, 
and for $(w_1,w_2) \in [0,\infty )$, let $\P_{w_1,w_2}$ be the law of 
 $((W^1_n, W^2_n), n\in \Z_+)$ starting at $(W^1_0, W^2_0) = (w_1, w_2)$.  
If $\mu$ is a law on $(0,\infty )$ or $(0, \infty )^2$, then write $\P_\mu$ 
for the law of the appropriate chain $(W$ or $(W^1, W^2))$ with initial 
distribution $\mu$. 

\no{\bf Lemma 5.5}.  For all $\theta, \varepsilon > 0$, there exist
$\delta_{5.1} , c_{5.3}, \lambda_{5.1} > 0$ 
(depending on $(\theta , \varepsilon ))$ such that
$$ \P_w (\sum^n_{j=0} 1 (W_j \leq \delta_{5.1} ) > \varepsilon n ) 
\leq c_{5.3} (w^{-\theta } + 1) e^{-\lambda_{5.1}
n}\ \ \hbox{ for all } n \in {\Z}_+,  w>0.$$

\no{\it Proof}.  Define a Markov kernel $q$ on $(0, \infty )$ by

$$ q((0,y]|w) =  \cases{
\rho((0,y/w ]) &if $w\leq w_0$\cr
p((0,y]|w_0 ) &if $w>w_0$.\cr}$$

\no Lemma 5.4 implies that 
$q(\cdot \, |w) \prec p (\cdot \, |w)$ for all $w > 0$. 
 This, the stochastic monotonicity of $p(\cdot \, |w)$, and a standard 
coupling argument (Kamae et al. (1977, Theorem 2)) show that for 
each $w
>0$ we may construct Markov chains $(W_n)$ and $(X_n)$ with 
transition kernels 
$p$ and $q$, respectively, on the same probability space such that 
$X_0 = W_0 = w$ and $W_n\geq X_n$ for all $n\geq 0$.  We abuse 
notation
slightly and let $\P_w$ denote the underlying probability. 
 Define a sequence of stopping times $(T_j, j\in \Z_+)$ by
$$ \eqalign{
T_0 &= \min \{ k \geq 0 : X_k > w_0\}\cr
T_{j+1} & = \min \{ k> T_j : X_k > w_0\}.}$$

\no
Let $\{ V_j : j\in \N\}$ be i.i.d. random variables, on some 
$(\Omega , {\cal F}, \P)$, with law $\rho$, and set \break
$S_n = \sum_{j=1}^n \log V_j$.  
Let $T^\prime_0 (w) = \min \{ k \geq 0:  w \prod^k_{j=1} V_j
> w_0 \}$, $w>0$.  The definition of $q$ implies that
$$\eqalign{ \P_w (((\log &X_{k \wedge T_0}, k \geq 0), T_0)\in 
\cdot \, )\cr
&= \P((( \log w + S_{k \wedge T^\prime_0}, k\geq 0), T^\prime_0 (w 
)) 
\in \cdot ) \hbox{ for all } w > 0.}\leqno (5.15)$$

\no
(5.14) implies $\gamma (\theta) \equiv \P (\exp (-\theta \log V_1 )) 
< \infty$
 for $\theta < 1$ and by (5.13) $\gamma^\prime
(0) = -\P (\log V_1 ) < 0$.
Choose $0<\theta <1$ sufficiently small $(0< \theta \leq \theta_0$ 
say)
 so that $\gamma (\theta )< 1$.  Then
$$\eqalign{
\P_w (T_0 > n)
 & = \P (\log w + S_k \leq \log w_0 \ \forall k\leq n)\ \ (\hbox{ by\ 
(5.15))}\cr
 & \leq \P (\exp (-\theta S_n)) (w_0/w)^\theta\cr
 &\leq w^{-\theta} \gamma (\theta )^n\ \ \ (\hbox{ recall\ } w_0 
\leq 1).\cr}\leqno (5.16)$$

\no The strong Markov property of $X$ now shows that $T_j < 
\infty$
for all $j\geq 0 $ a.s. 

If $\delta >0$ define
$$ Y_j = \sum^{T_j}_{k=T_{j-1}+1} 1(X_k \leq \delta)\ , j\in {\N};
\qq  Y_0 = \sum^{T_0} _{k=0} 1 (X_k \leq \delta).$$
For $j\in \N$, the strong Markov property of $X$ implies 
$({\cal G}_j = \sigma (X_0, \ldots , X_j))$
$$ \leqalignno{
\P_w (Y_j \in A|{\cal G}_{T_{j-1}})  & = \P_{X(T_{j-1})} 
(\P_{X(1)} (Y_0 \in A)) &(5.17) \cr
 & = \int \P_x (Y_0 \in A) p(dx |w_0 )\ \ (\hbox{ since}\ X_{T_{j-1}} > 
w_0).}$$

\no
Therefore $\{ Y_j : j\in \N\}$ are i.i.d. and
$$ \eqalign{
\P_w (Y_j)
 &= \int \P_x (Y_0) p(dx|w_0)\cr
 & \leq \int \P (\sum^\infty_{k=0} 1(\log x + S_k < \log \delta ))
 p(dx |w_0)\ \ (\hbox{ by}\ (5.15))\cr
 & \leq \sum^\infty_{k=0} \int \P (e^{-\theta S_k}) \delta^\theta 
x^{-\theta} p(dx|w_0)\cr
 & = (1 - \gamma (\theta ))^{-1} \int x^{-\theta} p(dx|w_0) 
\delta^\theta \ \ 
= c(\theta ) \delta^\theta ,\cr}\leqno (5.18)$$
where $c(\theta ) < \infty$ by (5.12) and (5.14).  Use (5.17) and 
(5.16)
 to show that for sufficiently small positive $\beta$
$$ \P_w (e^{\beta Y_j})  = \int \P_x (e^{\beta Y_0} ) p (dx |w_0)
  \leq e^\beta\int \P_x (e^{\beta T_0}) p(dx |w_0) < \infty. $$

\no
(5.18) implies ${d \over d\beta} \P_w (e^{\beta Y_j} ) |_{\beta=0}
\leq c(\theta )\delta^\theta$ and so for some
$\beta = \beta (\theta ) > 0$ we have
$$ \P_w (e^{\beta Y_j}) \leq 1 + 2c (\theta ) \beta \delta^\theta 
\leq \exp (2c (\theta ) \beta \delta^\theta ).\leqno(5.19)$$

Given $\varepsilon > 0$ choose $\delta > 0$ so that 
$c(\theta ) \beta\delta^\theta \leq \varepsilon / 8$.  Then
$$ \eqalign{
\P_w ( \sum^n_{j=0} 1 (W_j \leq\delta ) > \varepsilon n) 
 &\leq \P_w (\sum^n _{j=0} 1(X_j \leq \delta ) > \varepsilon n)\cr
 & \leq \P_w (T_0 +1 + \sum^n_{j=1} Y_j > \varepsilon n)\cr
 &\leq \P_w ( T_0 > \varepsilon n/2 -1) + \exp (-\beta \varepsilon 
n/2) 
\P_w(e^{\beta Y_1})^n\cr
 &\leq w^{-\theta} \gamma (\theta)^{[\varepsilon n/2]-1} + 
\exp (- \beta \varepsilon n/2 + 2c (\theta ) \beta \delta^\theta n)\cr
 &\qquad\qquad (\hbox{ by \  (5.16)\ and\ (5.19))}\cr
 & \leq c(w^{-\theta} + 1) e^{-\lambda n}\cr}$$

\no for some $c=c(\theta , \varepsilon )>0$ and 
$\lambda = \lambda (\theta, \varepsilon ) >0$ ( by the choice of
$\delta$).  This proves the result for $0< \theta \leq \theta_0$, 
and it follows trivially for all $\theta >
0$.\qed
\med
\no{\bf Proposition 5.6}.  For any $\theta >0$ there are 
$c_{4.4} (\theta ), \lambda_{5.2} (\theta ) >0$ such that
the coupling time $\tcpl $ satisfies
$$ \P_{w_1, w_2} (\tcpl > n) \leq c_{5.4} ((w_1 \wedge w_2) ^{-
\theta} +1)
 e^{-\lambda_{5.2}n}\ \
\hbox{ for all } n \in {\Z}_+, \, (w_1, w_2)\in (0, \infty )^2.$$
In particular, $\tcpl < \infty\ \P_\mu -$ a.s. for all laws $\mu$ on 
$(0, 
\infty)^2$.

\med {\it Proof}. Fix $0<w_2\leq w_1$, write $\P$ for
 $\P_{w_1, w_2}$, and let ${\cal G}_n =
\sigma ((W^1_j, W^2_j), j\leq n)$.  Note that $M(n)
= 1(\tcpl > n) \exp (\sum^{n-1}_{j=0} \Lambda (W^2_j /\alpha ))$
is a $({\cal G}_n)$-supermartingale because w.p. 1.
$$ \eqalign{
\P(M(n+1)|{\cal G}_n) &= M(n) \exp (\Lambda (W^2_n /\alpha ))
 \P(\tcpl >n +1| {\cal G}_n)\cr
 & = M(n) \exp (\Lambda (W^2_n/\alpha ))
 (1- \Lambda (W^2_n/\alpha ) /\Lambda (W^1_n/\alpha )) 
\hbox{ \ (by\ the\ definition of}\ \bar p )\cr
 &\leq M(n).\cr}$$
If $\theta > 0$ and $\varepsilon = 1/2$ in Lemma 5.5, then that 
result
 gives (for $\delta = \delta_{5.1}$)
$$ \eqalign{
\P(\tcpl >n)  &\leq \P(\tcpl >n, \sum^{n-1}_{j=0} 1(W^2_j \leq \delta 
) \leq n/2)
 + \P (\sum^{n-1}_{j=0} 1(W^2_j \leq \delta ) >
n/2)\cr
 &\leq \P(M(n))\exp (-\Lambda (\delta /\alpha ) n/2 ) 
+ c_{5.3} (w_2^{-\theta} + 1) e^{-\lambda_{5.1}n}\cr
 &\leq \exp (\Lambda (w_2/\alpha ) - \Lambda (\delta /\alpha ) 
n/2 ) 
+ c_{5.3} (w_2^{-\theta} + 1) e^{-\lambda_{5.1} n}.}$$

\no The result follows because $\Lambda (\delta/\alpha ) > 0$ 
(by (5.9)) and $\Lambda (w_2/\alpha ) \leq 1$.\qed

\med {\bf Corollary 5.7.} (a) There are Markov chains $(W_n, n\in 
\Z_+$)
 with law $\P_F$ and $(\overline{W}_n ,n\in \Z_+$) with law $\P_\pi$, 
defined on the same probability space, such that if
$$ \tcpl =\min \{ n\in {\Z}_+ : W_n = \overline{W}_n\}$$
then: 
\item{(i)}  $W_n = \overline{W}_n$ for all $n \geq \tcpl $,
\item{(ii)} $\P(\tcpl > n)\leq c_{5.5} e^{-\lambda_{5.3} n}$ for all 
$n \in \Z_+$ for some $c_{5.5}, \lambda_{5.3} >0$.

\no (b) If $W_n = T(X|n , \infty ) \alpha^{-n-1}$ and 
$|{\cal L} (W_n ) - \pi |$ denotes the total variation
distance between the law of $W_n$ and its weak limit $\pi$ then
$|{\cal L} (W_n) - \pi | \leq 2c_{5.5} e^{-\lambda_{5.3}n}$,  $n \geq 
0$.
 \med
{\it Proof}.  (a) Since $F \succ 
\pi$ by Theorem 5.1(b), there is a law $\mu$ on $D = \{ (w_1, w_2 ): 
w_1
\geq w_2 > 0\}$ with first and second marginals $F$ and $\pi$, 
respectively.  
Let $(W_n, \overline{W}_n)$ be a chain with
law $\P_\mu$.  Then by Proposition 5.6 with $\theta = {1\over 2}$
$$ \P_\mu (\tcpl >n) \leq c_{5.4} e^{-\lambda_{5.2} (1/2)n} 
\int w^{-1/2} + 1 d\pi (w) \leq c_{5.5} e^{-\lambda_{5.3}n}\q 
\hbox{ (by (5.5)).}$$
\no (b) is immediate because $(W_n)$ has law $\P_F$.\qed

The same argument as in (a) also gives the following corollary.

\med {\bf Corollary 5.8}. There exists $c_{5.6} > 0$ such that  
$$\P_{\pi\times \pi} (\tcpl >n) \leq c_{5.6} e^{-\lambda_{5.3} n} 
   \q \hbox{ for all } n \in \Z_+ \,. \leqno(5.20)  $$

\bigbreak
\no {\bf 6.  Limit Theorems}

In this section we use the results of the previous sections to prove a Law of 
Large Numbers (Theorem~6.1) and a Central Limit Theorem (Theorem~6.2), for 
functionals of the clusters which branch off the backbone. The next section 
then lists specific limit theorems that follow 
from these.  Recall the setting of Theorem~2.5: the
backbone is $\spine_1 , \spine_2 , \ldots$, the clusters off the backbone,
and associated percolation times to infinity are given by
$Y_{n , j}= (\CC_{n,j} , a_{n,j})$, and are conditionally independent
with distributions $\nu_{W_n}$.  Let $\Y_n$ denote the
vector $(Y_{n,1} , \ldots , Y_{n,d-1}) \in \sS_0^{d-1}$.

\no{\bf Notation}. If $\varphi : \sS_0^{d-1} \rightarrow \R$, let
 
$$\eqalign{ 
\mu_{t}(\varphi) 
   &= \int_{{\cal S}_0^{d-1}} \varphi (y_1 , \ldots , y_{d-1}) 
   \prod^{d-1}_{j=1} \nu_{t} (dy_j), \cr
\bar \mu(\varphi) &= \int^\infty_0 \mu_{t}(\varphi) 
   d\pi (t).\cr}$$
Thus $\mu_t$ is the law of $(d-1)$ independent copies of 
$\nu_t$.

\med \no {\bf Theorem 6.1}. If $\varphi : {\cal
S}_0^{d-1} \rightarrow \R$ is measurable and
$\bar \mu(\varphi)$ is finite, then
$$
\lim_{n \rightarrow \infty} n^{-1} \sum^{n-1}_{j=0} 
\varphi (\Y_j ) = \bar \mu(\varphi) \q \hbox{ a.s..}
$$
   
\no{\bf Theorem 6.2}.  Suppose that 
$$ \int^\infty_0 {\mu_t(\varphi^2)}^{2} \log^+ (\mu_t (\varphi^2))
d\pi(t) < \infty. \leqno(6.1)$$
Then 
$$  n^{-1/2} \sum_{j=0}^{n-1} (\varphi (\Y_j)- \bar \mu(\varphi))  
    \, \wconv \, N(0 , \sigma^2_{\varphi}) \quad \hbox{as } 
n \rightarrow \infty,$$
 where
$$ \sigma^2_{\varphi}= 2 \Bigl( \sum_{j=1}^\infty \P_{\pi}(\mu_{W_j}
   (\varphi) \mu_{W_0}(\varphi) - \bar \mu(\varphi)^2) \Bigr)
   + \bar \mu(\varphi^2) - \bar \mu(\varphi)^2 < \infty.
   \leqno(6.2)$$

\no{\it Proof of Theorem~6.1}.   Let $(W_j , j\in \Z_+$) and 
$(\overline{W}_j, j\in \Z_+)$ 
be the chains with laws $\P_F$ and $\P_\pi$, respectively, 
which are coupled as 
in Corollary 5.7 (a).  Let $\{ U_j : j\in \Z_+\}$ be an independent 
sequence of i.i.d. random variables which are uniformly distributed on 
$[0,1]$.  Define $G_t(x) = \mu_t (\{ (y_1, \ldots , y_{d-1} ): 
\varphi (y_1, \ldots , y_{d-1} ) \leq x\})$ and let $V_j (t) = G^{-1}_t 
(U_j)$ where $G^{-1}_t (u) = \inf \{ x: G_t (x) >u\}$.  Therefore 
$\{V_j(t):j\in \Z_+\}$ are i.i.d. and have distribution 
$\mu_t (\varphi (\cdot ) \in \cdot )$.  Moreover $V_j (t,\omega )$ is
jointly 
measurable because $G^{-1}_t (x)$ is. Theorem 2.5 shows that $(V_j (W_j), 
j\in \Z_+)$ and $(\varphi (Y^{(j)}), j\in \Z_+)$ are equal in law. 
By the coupling in Corollary 5.7 it suffices to show 
$$ \lim_{n\rightarrow\infty} n^{-1} \sum^{n-1}_{j=0} V_j (\overline{W}_j) 
= \bar \mu(\varphi) \hbox{ \ a.s.}\leqno (6.3)$$ 
\no However, $V_j(\overline{W}_j) = G^{-1}_{\overline{W}_j} (U_j)$ and
since $\{ (\overline{W}_j, U_j), j\in \Z_+\}$ 
is clearly stationary and ergodic by Theorem 5.2, $\{V_j (\overline{W}_j),
j\in \Z_+\}$ is too.  $\bar \mu(\varphi)$ is the mean of $V_j
(\overline{W}_j)$
and therefore the ergodic theorem implies (6.3). \qed 

\med {\it Proof of Theorem~6.2}.  Let $\{ (W_j, \overline{W}_j, V_j):j\in
{\Z}_+\}$ 
be as in the previous proof. As in the above argument it suffices to show 
$$
\Phi_n= n^{-1/2} \sum\limits_{j=0}^{n-1}\big( V_j(\overline{W}_j)-\overline{\mu} 
  (\varphi )\big) \wconv N(0 , \sigma_\varphi^2).
$$
Write $\Phi_n =X_n+Z_n$, where
$$
X_n=  n^{-1/2} \sum\limits_{j=0}^{n-1}
    \big( V_j(\overline{W}_j)-\mu_{\overline{W}_j} 
  (\varphi )\big)\, \hbox{ and } \,
   Z_n=   n^{-1/2} \sum\limits_{j=0}^{n-1} 
    \big( \mu_{\overline{W}_j} (\varphi )-\overline{\mu} (\varphi )\big).
$$
We will use the Lindeberg Central Limit Theorem to show that, conditional
on $\overline{\W}=(\overline{W}_j)$, $X_n \wconv $
$N(0,\sigma_1^2)$ and then use a Central Limit Theorem for stationary
ergodic processes to prove that 
$Z_n  \wconv N(0,\sigma_2^2)$.
Introduce
$$
s_n^2 = \sum\limits_{j=0}^{n-1} \P \Big( {\big(V_j (\overline{W}_j)
-\mu_{\overline{W}_j}(\varphi )\big)}^2 \Big| \overline{\W }\Big) 
      = \sum\limits_{j=0}^{n-1} \mu_{\overline{W}_j} (\varphi^2)
         -\mu_{\overline{W}_j} (\varphi )^2. 
$$
The Ergodic Theorem implies
$$
\lim\limits_{n\rightarrow\infty } s_n^2n^{-1}
  = \overline{\mu} (\varphi^2) - \int \mu_t (\varphi )^2\, 
     d\pi (t)\equiv \sigma_1^2 \,\, \hbox{ a.s.}\leqno(6.4)
$$
We claim that 
$$
\lim\limits_{n\rightarrow\infty } \P (e^{i\theta X_n} | \overline{\W} )
= e^{-\theta^2\sigma_1^2/2}\quad \hbox{ for all } \theta \in {\R }.\leqno(6.5)
$$
If $\sigma_1^2=0$ this is clear from (6.4). Assume $\sigma_1^2 >0$, let $\var 
>0$ and check the Lindeberg condition, conditional on $\overline{\W}$. For any 
fixed $K>0$, use (6.4) and the Ergodic Theorem to conclude 
 $$
\eqalign{
{\limsup\limits_{n\rightarrow \infty }} \,
& s_n^{-2} \sum\limits_{j=0}^{n-1} \P \left( {\big( V_j
(\overline{W}_j)-\mu_{\overline{W}_j} (\varphi )\big) }^2
1\big( |V_j(\overline{W}_j)-\mu_{\overline{W}_j} 
(\varphi )|>\varepsilon s_n\big) |\overline{\W}\right) \cr
&\leq {\lim\limits_{n\rightarrow \infty }}
\sigma_1^{-2} n^{-1}\sum\limits_{j=0}^{n-1}
\int {\big( \varphi (y)-\mu_{\overline{W}_j} (\varphi )\big)}^2
1 \big( |\varphi (y)-\mu_{\overline{W}_j} (\varphi )|>K\big)\, 
d\mu_{\overline{W}_j} (y)\hbox{ a.s. }\cr
&= \sigma_1^{-2} \int \int {\big( \varphi (y)-\mu_t (\varphi )\big) }^2
1\big( |\varphi (y)-\mu_t (\varphi )|>K\big)\,
d\mu_t (y)\, d\pi (t)\hbox{ a.s. }\cr  
}
$$ 
The last expression approaches zero as $K\rightarrow \infty $ because 
$\overline{\mu} (\varphi^2)<\infty $ by (6.1). This gives us the 
Lindeberg condition with respect to $\P (\cdot |\overline{\W})$
a.s., and (6.5) then follows from the Lindeberg Central Limit Theorem and
(6.4).

It is easy to use the exponentially fast coupling given by Corollary 5.8 to
see that $\{ \overline{W}_j \} $ is strongly mixing with an exponential
mixing rate (with the notation of Rio (1995)), 
$\alpha (n)\leq 2 c_{ 5.6} \exp(-\lambda_{5.3}n)$. This means that the
same is
true of the ergodic process 
$\{ \mu_{\overline{W}_j} (\varphi )-\overline{\mu }(\varphi )\} $, and
(6.1)
allows us to apply a Central Limit Theorem for strongly mixing stationary 
processes (Theorem~1 of Rio (1995) and (1.5) of that work) to conclude that
$$
\lim\limits_{n\rightarrow \infty } \P (e^{i\theta Z_n})
=e^{-\theta^2\sigma_2^2/2},\leqno(6.6)
$$
where
$$
\sigma_2^2 = \int \mu_t (\varphi )^2 \, d\pi (t)-
\overline{\mu}(\varphi )^2
+2\sum\limits_{j=1}^\infty \P 
\big(
\mu_{\overline{W}_j}(\varphi)\mu_{\overline{W}_0}(\varphi)-\overline{\mu}
(\varphi )^2\big) .
$$
The required result now follows from (6.5) and (6.6) because 
$\sigma_\varphi^2$ $=\sigma_1^2+\sigma_2^2$. \qed

\no{\bf Example 6.3}.  Define $\varphi_1: {\cal S}_0^{d-1} 
\rightarrow [0, \infty)$ by $\varphi_1 ((C_1,a_1), \ldots,
(C_{d-1} , a_{d-1})) = \sum^{d-1}_{i=1} \# C_i$.  Then 
$$ \eqalign{ \mu_t(\varphi_1^2)  &= 
(d-1) \int (\# C)^2 d\nu_t + (d-1) (d-2) (\int \# Cd\nu_t)^2\cr
 & \leq (d-1) c_{4.8}(\alpha,d) + (d-1)(d-2)c_{4.3}(\alpha,d)^2\ \ \ 
\hbox{ (Theorems\ 4.4\ and\ 4.6).}}$$
If $C_n = \cup_{i=1}^{d-1} \CC_{n,i}$, then $\varphi_1 (Y^{(n)}) = \# C_n$,
and Theorems 6.1 and 6.2 show that
$$ \lim_{n\rightarrow \infty} n^{-1} \sum^{n-1}_{j=0} \# C_j
= \bar \mu \equiv (d-1) \int\int \#C d\nu_t d\pi (t) \q 
\hbox{a.s.,} \leqno(6.7)$$
and
$$ n^{-1/2} \sum^{n-1}_{j=0} (\#C_j - \bar \mu)
\wconv N(0, \sigma^2_{\varphi_1})\ \
\hbox{ as\ }n\rightarrow \infty,\leqno (6.8)$$

\no where a (not very helpful) expression for $\sigma^2_{\varphi_1}>0$ may
be  retrieved from Theorem 6.2.  The expression for $\bar \mu$ is more 
tractable and Theorem 4.4(b) and (4.11) imply that
$$ \leqalignno{
 0 < \bar \mu = \bar \mu (\alpha , d) & 
 \leq (d-1) \sup_t \int \# C d\nu_t &(6.9) \cr
  & \leq (d-1) c_{4.3}(\alpha,d) \cr
  & \leq c_{6.1}(d) \exp (c_{4.7}(d) /(1-\alpha )).}$$

We now consider some processes which describe aspects of the growth
 dynamics of $\{ A_n\}$.  The size of the cluster at the
first time a node of height $n$ is filled is $\minsize (n) 
= \min \{ k: h(A_k) = n\}$.  The size of the cluster at the
last time a node of height $n$ or less is added is 
$\maxsize (n) = \max \{ k: A_k-A_{k-1} \in \B(0,n)\}$.  Let $\cutsize
(n) = \# (A_\infty \cap \B(0,n))$ be the number of nodes
 in $A_\infty$ of height $n$ or less. Then clearly
$$ \minsize (n) \leq \cutsize (n)\leq \maxsize (n).\leqno (6.10)$$

\no Note that each point in $A_\infty \cap \B(0,n)$
 is in a cluster $\CC_{k,i}$ for some $k< n$ or is one of the first
$n+1$ vertices along the backbone.  
If $N(n) = \sum^{n-1} _{k=0} \#(C_k)$, this gives
$$ \cutsize (n) \leq N(n) +n+1.\leqno (6.11)$$

For the next result $\sigma_{\varphi_1}^2$ is defined in Theorem 6.2
with $\varphi_1$ as in Example 6.3, and we define $$ \rate = \bar \mu
(\alpha , d) + 1\ >1.\leqno(6.12)$$

\med {\bf Theorem 6.4}.  For $\{ L(n) \} = \{ \cutsize (n)\}, 
\{ \maxsize (n)\}$ or $\{ \minsize (n)\}$ the following results hold. 

\no{(a)}\quad $\dst \lim_{n\rightarrow\infty} n^{-1} L(n) = \rate$ a.s.

\no{(b)} \quad $\dst n^{-1/2} (L(n) - n\rate ) \wconv 
N(0, \sigma^2_{\varphi_1})$ as $n\rightarrow\infty$.

\bigskip 
It is not hard to prove this by first using the bounds in Theorem 4.4
to show that $\maxsize (n) - \minsize (n)$  and $N(n)+n+1 - \cutsize (n)$ 
remain bounded in probability as $n\rightarrow\infty$, and then applying 
Example 6.3.  The interested reader may find this argument in Lemma 6.5 of the 
earlier version of this work referred to in the introduction. We will prove 
Theorem 6.4 in the next Section by showing that the above differences are 
bounded in $L^1$ by means of a dynamical decomposition of $A_\infty$ into 
independent blocks. 

Let $\ell (n) = \min \{ |x|: x\in A_\infty - A_n\}$ be the height of the 
shortest vertex which is added after step
$n$.  Hence $A_k \cap \B(0, \ell (n)-1)$ is ``frozen" for $k>n$.  
Clearly $\ell (\maxsize (k)-1) \leq k < \ell (\maxsize (k))$ and so
$$ \maxsize (k) = \min \{ n : \ell (n)>k\}.$$

\no It is a simple matter to read off limit theorems for 
$h(A_n)$ and $\ell (n)$ from the corresponding results for their
inverses, $\minsize$ and $\maxsize$, respectively. Let 
$\mu_0 = \mu_0(\alpha,d) = \rate^{-1}.$

\bigskip \no {\bf Corollary 6.5}.  
(a) $\dst \lim_{n\rightarrow\infty} n^{-1}  h (A_n) 
 = \lim_{n\rightarrow\infty} n^{-1} \ell (n)  = \mu_0$ a.s.

\no{(b)(i)} \quad $ n^{-1/2} (h (A_n) - n\mu_0) \wconv N(0, \sigma^2_{\varphi_1}
\mu_0^{3})$ as $n\rightarrow \infty$.

\no \phantom{(b)} {(ii)} \quad $ n^{-1/2} (\ell (n) - n \mu_0) 
\wconv N(0, \sigma^2_{\varphi_1}
\mu_0^{3})$ as $n\rightarrow \infty$.

\med {\it Proof}.  (a) is a trivial consequence of Theorem 6.4 (a)
 with $L(n) = \minsize (n)$ or $\maxsize (n)$.

\no (b) Let $H(n) =n^{-1/2}  (h (A_n) - n \mu_0)$ and $\Sigma (k) = 
k^{-1/2} (\minsize (k) - k \mu_0^{-1} ) $.  Fix $x \in \Bbb R$ and set
$k(n) = [n\mu_0 + x \sqrt n \,]+1$ where $[z]$ is the greatest integer 
not exceeding $z$. If 
\hfil\break $x_n = (n-k(n) \mu_0^{-1}) k(n)^{-1/2}$, then for $n$
large enough so that 
$k(n) \in \N$, 
$$ P(H_n \leq x) = P(h(A_n) < k(n)) = P(\minsize (k(n)) > n) = P (\Sigma 
(k(n)) > x_n).$$ 

\no Since $\lim x_n = -\mu_0^{-3/2} x$, Theorem 6.4(b) with $L(k) = 
\sigma (k)$ implies that
$$ \lim_{n\rightarrow\infty} P(H_n \leq x ) = P(Z> -\mu_0^{-3/2} x) = 
P(Z \leq \mu_0^{-3/2} x),$$

\no where $Z$ is a $N(0, \sigma_{\varphi_1}^2)$ random variable. (i)
follows and a similar argument proves (ii). \qed 

It is easy to translate these a.s. limit theorems as $n\rightarrow\infty$ 
into continuous time results as $t\uparrow T_\infty$.  Note that
$$ W_{k+1} \leq (T_\infty - T^{(k)} ) \alpha^{-k-1} 
\leq W (\minver_k),$$
where $\minver_k$ is the first vertex in $\B(k)$ which is added to the
cluster. Note that
$W(\minver_k)$ is equal in law to $W_0$, and $W_{k+1}\succ  \pi$
for all $k$ by Theorem 5.1(b).  
Use the above with the estimates (5.1) and (5.5) (for the 
left-hand tail of $\pi$) and a Borel-Cantelli argument to see that
$$ (k-1)^{-1} \alpha^{k-1} \leq T_\infty - T^{(k)} 
\leq (\log k)\alpha^k \ \hbox{ for large $k$ a.s.}$$

\no By considering $T^{(k)} \leq t< T^{(k+1)}$ this gives
$$ h(C_t)^{-1} \alpha^{h(C_t)} < T_\infty - t < \log (h(C_t)) 
\alpha^{h(C_t)} \ \hbox{ for\ } 0< T_\infty-t \ \hbox{
sufficiently small, a.s.}\leqno (6.13)$$

\med{\bf Theorem 6.6}.  (a) $\dst \lim_{t\uparrow T_\infty} (h(C(t))
 (\log_\alpha (T_\infty - t))^{-1}= 1$ a.s.

\no {(b)} $\dst \lim_{t\uparrow T_\infty} \# (C(t))
 (\log_\alpha (T_\infty - t))^{-1} = \rate$  a.s.

\no{\it Proof}. (a) is immediate from (6.13).
Corollary 6.5 (a) and a trivial interpolation (recall $A_n =
 C(\grow_n))$ shows that \no $\lim_{t\uparrow T_\infty} h(C(t))
/\# C(t) = \mu_0$  a.s.  (b) follows from this and (a).  \qed

\med{\bf Theorem 6.7}.  (a) $\dst \lim_{n\rightarrow \infty}
n^{-1} (\log_\alpha (T_\infty - \grow_n))  = \mu_0$ a.s.
\item {(b)} $\dst n^{-1/2} (\log_\alpha (T_\infty - \grow_n) - n \mu_0)  
\wconv N(0,\sigma^2_{\varphi_1} \mu_0^{3})$ as $n\rightarrow\infty$.

\no{\it Proof}.  For (a) set $t= \grow_n$ in Theorem 6.6(b).  (b) follows from 
Corollary 6.5(b)(i) and 
$$ \lim_{n\rightarrow\infty} n^{-1/2} (h(C(\grow_n)) - \log_\alpha (T_\infty - 
\grow_n))  = 0 \hbox{ a.s. } \leqno (6.14)$$ 
\no (6.14) is an easy consequence of (6.13) and the trivial bound
 $h(C(\grow_n )) \leq n.$\qed

\bigb {\bf 7.  A Decomposition of the Infinite Cluster into i.i.d. Blocks.}
\medskip
The drawback of the decompositions of $A_\infty$ into clusters off the
backbone (Theorem 2.5) is that it is not a dynamical decomposition.
This means that some additional work is needed before the results on
the growth rate of $h(A_n)$ such as Corollary 6.5 can be derived from the
limit theorems for the clusters $(\CC_{n,i})$ given in Example 6.3.  We now
establish a dynamical decomposition of $A_\infty$ into i.i.d. pieces,
which will lead to a proof of Theorem 6.4, and 
will also be used in Section 8 to establish properties of the shape
of the cluster ``as viewed from the tip".

Let ${\cal C}_n = \sigma (Y_{k,i}, W_k: k< n, i<d) \vee \sigma (1(W_n \leq 
\delta_0)) \vee \sigma (\spine_n)$.  Here $\delta_0$ is a sufficiently small 
positive number whose precise value will be prescribed below. Let $\bar
{\cal 
C}_n = {\cal C}_n \vee \sigma (W_n)$.  Clearly both $({\cal C}_\cdot)$ and 
$(\bar {\cal C}_\cdot)$ are filtrations. 

Let $M_n = \min \{ a_{n,i} (x) : x \in \CC_{n,i}, i<d\}$ and recall $C_n =
\cup_{i<d} \CC_{n,i}$.  Inductively define regeneration times
$\{ \reg_j: j\in \Z_+\}$ by 
$\reg_0 =0$ and
$$\eqalign{
\reg_{j+1}= &\min \{ k> \reg_j:\ \min_{\reg_j\leq n<k} M_n \alpha^n >
\delta_0 \alpha^k, C_n \subset \B(0,k-n-2) \cr
   &\hbox{ for all } n\in [\reg_j,k),
  W_{k-1} > \alpha \delta_0, W_k \leq \delta_0\}.}$$

\no Here $\delta_0>0$ will be chosen below.  Clearly $\reg_j$ is a 
$({\cal C}_n)$-stopping time for all $j\in \Z_+$.
The first two and last conditions in the above inductive
definition will imply that all the points in $A_\infty \cap \B(0,\reg_j)$
are added to the cluster before all the points in $A_\infty \cap \B(\reg_j,
\infty)$ and that $\{ \spine_{\reg_j} \}$ is the only point in $A_\infty$
of generation $\reg_j$.  The next to last condition will ensure the
independence of the blocks $B(j)$, defined by
 $$B(j)=A_\infty \cap \B (\reg_j,\reg_{j+1}-1).$$ 

Our goal is to prove the following theorem (for a sufficiently small
$\delta_0 >0)$. 

\no{\bf Theorem 7.1}. The sequence of $({\cal C}_n)$-stopping
times, $\{\reg_j:j\in \Z_+ \},$ is a.s. finite and satisfies

\item {(a)}  $A_\infty \cap \B (\reg_j) = \{ \spine_{\reg_j}\}$ and 
$\minsize(\reg_j)
= \maxsize (\reg_j) \hbox{ for all } j\in \Z_+.$
\item{(b)} $\P((\# B(j))^2) \leq c_{7.1} 
(\alpha , d)  \hbox{ for all } j\in \Z_+.$
\item {(c)} By (a) we may define $D_j \in S_0$ by
$$ B(j) = \{ (\spine_{\reg_j}) \oplus x:x\in D_j\},\ j\in \Z_+.$$
\no Also let $T(j;x) = T(\spine_{\reg_j}, 
(\spine_{\reg_j})\oplus x)\alpha^{-\reg_j}$ for
$x\in D_j$.  Then $\{ (D_j, T(j) ) : j\in \Z_+ \}$ 
are independent ${\cal S}_0$-valued random vectors and are
 identically distributed for $j\geq 1$.  Moreover for $j\geq 1$,
$$ \P((D_j, T(j)) \in \cdot \, | {\cal C}_{\reg_j}) 
= \int^{\delta_0}_0 P((D_0, \TT^{(0)})\in \cdot \, | W_0 = w )
d\Lambda (w ) \Lambda ([0,\delta_0])^{-1}.$$
\med
\no{\bf Remark}.  We abuse notation slightly, and extend the definition
of $T(j;x)$ to all $x \in \B$. Clearly 
$T(0; x) = \TT^{(0)} (x) \hbox{ for all } x\in \B.$ 

\med {\bf Lemma 7.2}. $\bar {\cal C}_n\subset {\cal W}_n 
\hbox { for all } n \in \Z_+.$

\med {\it Proof}.  By Theorem 5.1 $W_n$ is ${\cal W}_n$-measurable, and
 $\spine_n$ is trivially ${\cal W}_n$-measurable.  Fix $k<
n$, $i<d,$ and $x,x^\prime$ in $\B$.  Then for a Borel set $B$
$$ \{ a_{k,i} (x^\prime ) \in B, \spine_n = x \} 
= \{ W(x|k) - T (x|k , e_i (x|k+1) \oplus x^\prime ) \alpha^{-k-1} \in B,\ 
\spine_n=x \} \in \bar {\cal F}_x.$$

\no Therefore $Y_{k,i}$ is ${\cal W}_n$-measurable.\qed

\no{\it Proof of Theorem 7.1(a)}.   
Assume $\reg_j < \infty$ for a fixed $j\in \N$.  If $n< \reg_j, i<d$
and 
$x\in \CC_{n,i}$, then
$$ \leqalignno{ T(e^n_i \oplus x) = T(\infty) - \alpha^{n+1} a_{n,i} (x) 
&< T(\infty) -\alpha \delta_0 \alpha^{\reg_j} &(7.1)\cr
 & \leq T(\infty) - \alpha^{\reg_j+1}W_{\reg_j} = T(\spine_{\reg_j}).\cr}$$

\no Note we have used the fact that if $\reg_j < \infty$ 
and $j\in \N$, then $M_n\alpha^n > \delta_0 \alpha^{\reg_j}$ for
all $n<\reg_j$ (and not just $\reg_{j- 1}\leq n< \reg_j)$.  Similarly
we have $C_n \subset
\B(0,\reg_j - n-2)$ for all $n<\reg_j$ and this clearly implies that
all the clusters which break off the backbone before $\reg_j$ are in
$\B(0,\reg_j-1)$ and therefore $$ T(\reg_j) = T(\spine_{\reg_j}),\q \minsize
(\reg_j) = \# C(T(\spine_{\reg_j})).\leqno (7.2)$$

\no By (7.1) no new points are added to the first $\reg_j$ clusters 
$(C_0, \ldots C_{\reg_j-1})$ after time $T(\spine_{\reg_j})$.  
Hence $\spine_{\reg_j}$ 
is the last 
node in $\B(0,\reg_j)$ added to $A_\infty$.  This means 
 $$ \maxsize (\reg_j) = \# C (T (\spine_{\reg_j})).\leqno (7.3)$$

\no Clearly $A_\infty \cap \B(\reg_j) = \{ \spine_{\reg_j}\}$ since 
$\spine_{\reg_j}$ is the 
last point added to $A_\infty$ in $\B(0,\reg_j)$ and the first point added
to 
$A_\infty$ in $\B (\reg_j)$.  (7.2) and (7.3) also give 
$\minsize (\reg_j) = 
\maxsize (\reg_j)$ and hence (a) is proved once we show 
$\reg_j < \infty$ for all
$j \in \Z_+$.

The proof of Theorem 7.1(b) requires several preparatory lemmas.

\med {\bf Lemma 7.3}.  (a) $\nu_t (\min_{x\in C} a(x) \leq \delta) 
\leq c_{4.3} (\log^+ (1/\delta) +1) (\delta \wedge 1) \equiv g(\delta )$
for all $t$, $\delta \geq 0$.

\item {(b)} There exists $\beta \in (\alpha , 1)$ such that
 $\nu_t (C\not\subset \B(0,j-1)) \leq \beta^{j+1}$ for all
$j\in \Z_+$, $t\geq 0$.

\no{\it Proof}.  (a) We have 
$$ \eqalign{
\nu_t (\min_{x\in C}a(x) \leq \delta) & \leq \P(\max_{x\in C(t)} T(x) \geq
t-\delta|T(\infty) > t)\cr
 & \leq \P(\# (C(t) - C((t-\delta)-)) | T(\infty) > t)\cr
 &\leq g(\delta)\ \ \ \hbox{ (by\ Theorem\ 4.4(b)).}}$$

\item {(b)} Using Theorem 4.4(a) we have 
$$ \eqalign{ \nu_t (C\not\subset \B(0,j-1)) 
& \leq \P(\sum_{n\geq j} \# (C(t) \cap \B (n)) |T(\infty) >t)\cr
& \leq \sum_{n\geq j} c_{4.1}^d c_{4.2} \alpha^{n+1} 
  = c_{4.3} \alpha^{j+1},}$$ 
while
$$ \nu_t (C = \emptyset ) = \P( U_0 > t|T(\infty) > t) 
= (e^t \P(T(\infty) > t))^{-1} \geq c_{4.1}^{-1} \hbox{ (by Lemma 4.1).}
$$
The result follows trivially from the above two inequalities. \qed

\med{\bf Lemma 7.4}.  If $0<\delta_0 \leq \delta_{5.1} (1/2, 1/2)$
 ($\delta_{5.1}$ as in Lemma 5.5) there exist $c_{7.2},
c_{7.3} > 0$ such that $N_m = \sum^m_{k=1} 1 (W_{k-1} 
> \alpha \delta_0 , W_k \leq \delta_0)$ satisfies
$$ \P(\exp (-N_m) |W_0)\leq c_{7.2} (W_0^{-1/2} + 1) e^{-c_{7.3}m} 
\q \hbox{ for all } m\in \N.$$

\no{\it Proof}.  For $\delta_0$ as above and a fixed $q>0$
 let $M_m = \sum^m_{k=1} 1(W_{k-1} > \alpha
\delta_0 )(q-1(W_k \leq \delta_0))$.  Then
$$ \eqalign{
\P(e^{M_m} |{\cal W}_{m-1}) &= \exp (M_{m-1} + q 1(W_{m-1} 
> \alpha \delta_0))\cr
 &\qquad \times [ \exp (-1 (W_{m-1} > \alpha \delta_0)) \Lambda 
([0,\delta_0 \wedge W_{m-1} /\alpha )]) \Lambda
([0,W_{m-1} /\alpha])^{-1}\cr
 &\qquad + \Lambda ((\delta_0 , W_{m-1} /\alpha ]) 
\Lambda ([0, W_{m-1}/\alpha ))^{-1} ] \ \ \hbox{ (Theorem 5.1)}\cr
 &\leq \exp (M_{m-1} ) [1(W_{m-1} \leq \alpha \delta_0)\cr
 &\qquad + 1(W_{m-1} > \alpha \delta_0 ) e^q \bigl(e^{-1} 
\Lambda ([0,\delta_0 ]) + \Lambda ((\delta_0, \infty ))\bigl) 
\Lambda([0, \infty ))^{-1}].}$$

\no As $\Lambda ((0, \delta_0])>0$ (see (5.9)), we may choose 
$q=q(\delta_0)> 0$  small enough such that $e^{M_m}$ is a supermartingale. 
If $N^\prime_m = \sum^m_{k=1} 1(W_k \leq \alpha \delta_0)$, then
$$ e^{M_m}  = \exp (q m - q N_m^\prime - N_m)
 \geq 1(N^\prime_m \leq m/2) \exp (qm/2 - N_m) $$
and therefore
$$ \eqalign{ \P(\exp (-N_m ) |W_0) & \leq \P(N^\prime_m > m/2 |W_0 ) 
+ e^{-qm/2} \P(e^{M_m} |W_0)\cr
 & \leq c_{5.3} (W_0^{-1/2} + 1) e^{-\lambda_{5.1}m} + e^{-qm/2}.}$$

\no In the last line we used Lemma 5.5 and the supermartingale property 
of $e^{M_m}.$\qed
\med
\no{\bf Lemma 7.5}.  $\P(W_{\reg_j} \in B | {\cal C}_{\reg_j}) = 
\Lambda (B| \, [0, \delta_0])$ a.s. on $\{\reg_j < \infty \} $
for all $j\in \N.$

\no{\it Proof}.  If $A\in {\cal C}_n$ we claim that
$$ \P(A| {\bf W}) = \P (A|W_0, \ldots ,W_{n-1} , 1(W_n \leq
\delta_0)).\leqno(7.4)$$

\no To see this consider
$$ A= \bigcap^{n-1}_{k=0} \bigcap_{i<d} \{ Y_{k,i} \in A_{k,i} \} 
\bigcap \{ \spine_n = x\} \cap \{ (W_0, \ldots , W_{n-1} )
\in B, W_n \in D\}$$

\no where $D = [0, \delta_0]$ or $(\delta_0, \infty)$.  (7.4) is then an
easy 
consequence of Theorem 2.5 and the independence of $\spine_n$ and $\bf{W}$
(which 
holds by symmetry).  (7.4) follows for general $A \in {\cal C}_n$ by a 
monotone class argument. 

Assume $A\in {\cal C}_n$ is a subset of 
$\{ W_{n-1} > \alpha \delta_0, W_n \leq \delta_0\}$.  Then
$$ \eqalign{
\P(A,W_n \in B) & = \int \P(A|W_0,\ldots , W_{n-1} , 1(W_n \leq \delta_0)) 
1 (W_n \in B) d\P \ \hbox{ (by\ } (7.4)) \cr
 &= \int \P(A| W_{n-1} , 1(W_n \leq \delta_0)) 
\P(W_n \in B |W_{n-1} , 1(W_n \leq \delta_0))d\P, \cr}$$
 by the Markov property of ${\bf W}$. 
Use the form of the transition kernel found in Theorem 5.1 to see that
$$ \eqalign{ \P(W_n \in B |W_{n-1}, 1 (W_n \leq \delta_0)) 
= \Lambda (&B |[0, \min (\delta_0, W_{n-1} /\alpha)]) 1 (W_n \leq
\delta_0)\cr
&+ \Lambda (B|(\delta_0, W_{n-1}/\alpha ]) 1 (W_n > 
\delta_0 , W_{n-1} /\alpha > \delta_0).}$$

\no Therefore we have (by our assumption on $A$)
$$ \leqalignno{
\P(A,W_n \in B) &= \int \P(A|W_{n-1}, 1(W_n \leq \delta_0)) 
\Lambda (B |[0, \delta_0]) d\P &(7.5)\cr
 &= \P(A) \Lambda (B|[0, \delta_0]).}$$
Let $j\in \N$,  $A \in {\cal C}_{\reg_j}$ and 
$A_n = A\cap \{ \reg_j = n\} \in {\cal C}_n.$  Clearly
$A_n \subset \{ W_{n-1} > \alpha \delta_0, W_n \leq 
\delta_0\}$ and so by (7.5)
$$\eqalignno{
\P(A, \reg_j < \infty , W_{\reg_j} \in B) &= \sum^\infty_{n=1} 
\P(A_n, W_n \in B)\cr
 & = \P(A, \reg_j < \infty ) \Lambda (B|[0, \delta_0]). &\square}$$

The next result is an easy consequence of Theorem 2.5.

\med{\bf Lemma 7.6}.  Let $K$ be a $({\cal C}_n)$-stopping time 
(possibly infinite).  Then conditional on

\no ${\cal C}_{K} \vee \sigma ({\bf W})$ and on 
$\{ K < \infty \}$, $\{ Y_{K+n,i} : n\in \Z_+, i<d\}$
are independent ${\cal S}_0$-valued random vectors such that
$$ \P(Y_{K+n, i} \in A|{\cal C}_{K} \vee \sigma ({\bf W})) 
=\nu_{W_{K+n}} (A) \ \hbox { a.s. on } \{K <\infty\}.$$

\no{\it Proof of Theorem 7.1(b)}.  Choose $\theta < \delta_0 \leq
\delta_{5.1} 
({1\over 2}, {1\over 2}) \wedge 1$ small enough so that
$$ p_1 (\delta_0)= 1 - \prod^\infty_{j=1} (1-g (\delta_0 \alpha^j))^{d-1} 
+ 1-\prod^\infty_{j=1} (1-\beta^j)^{d-1} <1.\leqno (7.6)$$
Here $g$ and $\beta$ are as in Lemma 7.3.  Choose $n_0 \in \N$
 sufficiently large so that if \break
$ h(\delta_0) =\ab c_{4.3} \delta_0 ( 1 + \log 1/\delta_0 + \log 1/\alpha
)$,
 then (recall $\beta >\alpha )$
$$ (d-1) h (\delta_0) j \alpha^j \leq \beta^j \hbox{ \ for\ }j>
  n_0 \hbox{ \ and\ } \gamma = \beta^{n_0}
(p_1 (\delta_0) - \beta^{n_0} ) d (1-\beta )^{-1} + p_1 (\delta_0)<1.\leqno
(7.7)$$
\no Fix $j\in \Z_+$ and let 
$ k(0) \equiv  \reg_j$, and  
$ k(0) < k(1) <  \ldots <  k(n ) < \ldots $ 
denote the successive times for which 
$W_{k (n)} \leq \delta_0$ and $W_{k(n)-1} > 
\alpha \delta_0\ (n \in \N)$.  If $\reg_j = \infty$ set 
$ k(n) = \infty$ for all $n$.  Each $ k(n)$ is a
 $(\sC_n)$-stopping time.  Lemma 7.4 and the strong Markov
property of $W$ show that each $k(n)$ is finite if 
$\reg_j$ is.  Let $K(i) = k (in_0)$ and define 
$$ 
\eqalign{  &B_i = \{ \min \{M_n\alpha^n : \reg_j \leq n < K (i) \} \leq 
\delta_0 \alpha^{K (i)}\},\cr
 & D_i = \{ C_n \not\subset \B(0, K (i) - n-2) 
\hbox{\ for\ some\ }n \hbox{\ in\ } [\reg_j, K (i))\},\cr
 & A_N = \cap^N_{i=1} (B_i \cup D_i).}$$
\no Note that $\reg_{j+1} - \reg_j>m$ and $K(N) \leq m + \reg_j$ 
together imply $A_N$ (since each $K(i)$ for $i\leq N$ must
violate one of the defining conditions for $\reg_{j+1})$.  If 
$$ N^\prime_m = \sum^m_{k=1} 1(W_{\reg_j + k} \leq \delta_0,
W_{\reg_j +k-1} > \alpha \delta_0),$$
 
\no this easily gives for $\varepsilon > 0$
$$ \{ \reg_{j+1} - \reg_j >m\} \subset \{ N^\prime_m \leq \varepsilon m\}
\cup A_{[\varepsilon m / n_0]}.\leqno (7.8)$$

\no If $\omega \in A_N$ and $\reg_j(\omega) < \infty,$ 
then $\omega \in B_N \cup D_N,$ and so for some $n$ in $[\reg_j, K(N)),$
either $C_n$ is not contained in $\B(0,K(N)-n-2)$ or $M_n \alpha^n \leq
\delta_0 \alpha^{K(N)}.$  Choose $i=i(\omega ) \in \{ 1, \dots , N\}$
such that $K(i-1) \leq n < K(i)$ and note that $\omega \in A_{i-1}$
as well ($\omega \in A_j$ for all $j \leq N$). This shows that
on $\{\reg_j < \infty \}$,
$$ \eqalign{
  \P(&A_N | {\cal C}_{\reg_j} \vee \sigma ({\bf W}))\cr
 & \leq \sum^N_{i=1} \P (A_{i-1} \cap \{ \min (M_n \alpha^n : K (i-1) 
\leq n < K (i)) \leq \delta_0 \alpha^{K (N)}\} |
{\cal C}_{\reg_j} \vee \sigma ({\bf W}))\cr
 & \qquad + \P (A_{i-1} \cap \{ C_n \not\subset \B(0, K (N) - n-2)
\hbox{\ for\ some\ } n \in [ K (i-1) , K (i))\}
|{\cal C}_{\reg_j} \vee \sigma ({\bf W}))\cr
 & = \sum^N_{i=1} \P (1 (A_{i-1}) \{ 1 - \prod^{K (i)-1}_{n=K (i-1)}
 \nu_{W_n} ( \min_{x\in C} a(x) >
\delta_0 \alpha^{K (N) - n} ) ^{d-1} +1 \cr 
 &\qquad - \prod^{K (i) - 1}_{n=K  (i-1)} \nu_{W_n} 
( C \subset \B(0, K (N) - n-2))^{d-1} \} | {\cal
C}_{\reg_j} \vee \sigma ({\bf W})) \q \hbox{ (by\ Lemma\ 7.6)} }$$
$$\eqalign{
 &\leq \sum^N_{i=1} \P(A_{i-1} | {\cal C}_{\reg_j} \vee \sigma ({\bf W})) 
\Big\{ 1-\prod_{n=K (i-1)}^{K(i)-1} (1-g (\delta_0 \alpha^{K (N)-n}))^{d-1}\cr
 &\qquad\qquad +1 - \prod^{K(i) -1}_{n=K (i-1)} (1-\beta^{K(N)-n})^{d-1} 
\Big\} \q \hbox{ (Lemma 7.3)}\cr
 &\leq \P(A_{N-1} | {\cal C}_{\reg_j} \vee \sigma ({\bf W})) p_1
(\delta_0)\cr
 &\qquad + \sum^{N-1}_{i=1} \P(A_{i-1} | {\cal C}_{\reg_j} 
\vee \sigma ({\bf W}))
 (d-1) (\sum^{K (i) - 1}_{n=K (i-1)} g (\delta_0
\alpha^{K (N) -n} ) +\beta^{K (N)-n})\cr
 & \leq \P(A_{N-1} | {\cal C}_{\reg_j} \vee \sigma ({\bf W})) p_1
(\delta_0)+
 \sum^{N-1}_{i=1} \P (A_{i-1} | {\cal C}_{\reg_j} \vee \sigma ({\bf W}))
\sum^\infty_{m=n_0 (N-i)+1} {\kern -2 pt} d \beta^m.}$$

\no In the last line we used the first part of (7.7) and the inequality
 $K (N) - K (i) \geq (N-i) n_0$.  A simple induction
argument using the above and the definition of $\gamma$ in (7.7) gives
$$ \P(A_N | {\cal C}_{\reg_j} \vee \sigma ({\bf W})) \leq \gamma^N\ \ 
\hbox {for all } N\in \N \hbox{ a.s. on\ }\{ \reg_j <
\infty\}.\leqno (7.9)$$

\no Use (7.8) with $\varepsilon = c_{7.3} / 2$ and the strong Markov
 property of ${\bf W}$ with respect to $(\bar {\cal C}_n)$ (recall
Lemma 7.2 and Theorem 5.1) to see that if $m\in \N$ (recall $\P_w$ 
is the law of ${\bf W}$ starting at $w$)
$$ \eqalign{
\P(\reg_{j+1} & - \reg_j > m |\bar {\cal C}_{\reg_j})\cr
 & \leq \P(N^\prime_m \leq \varepsilon m |\bar {\cal C}_{\reg_j} ) + 
\P(A_{[\varepsilon m/n_0]} | \bar {\cal C}_{\reg_j})\cr
 & \leq \exp (\varepsilon m) \P_{W(\reg_j)} (e^{-N_m}) + \gamma^{-1} 
(\gamma^{\varepsilon /n_0 })^m\ \ \ \hbox{ (by \ (7.9))}\cr
 & \leq c (W(\reg_j)^{-1/2} + 1) e^{-c^\prime m } \q \hbox{ (by\
 Lemma\ 7.4\ and\ the\ choice\ of\ } \varepsilon ).}$$

\no Now condition on ${\cal C}_{\reg_j}$, use $F(t) \leq dt$ (see (5.9)) if
 $j=0$, and this together with Lemma 7.5 if $j\geq 1$, to
derive
$$ \P(\reg_{j+1} - \reg_j > m | \sC_{\reg_j} ) 
\leq c_{7.4} e^{-c_{7.5} m}\ \ 
\hbox{ for all } j \in \Z_+, m \in \N.\leqno (7.10)$$
Let
$$ M^\prime_n = \sum^{n-1}_{i=0} 1 (\reg_j \leq i < \reg_{j+1})
 (\# C_i - (d-1) \nu_{W_i} (\# C)), \ n \in \Z_+.$$
\no Lemma 7.6 shows $( M_n^\prime , {\cal C}_n \vee \sigma (\bf{W}))$ is 
a martingale and Theorems 2.5 and 4.6, and (7.10) readily show it is
$L^2$-bounded.  By (a) each point in $B(j)$ 
either belongs to a cluster which branched off the backbone
at generation $i\in [\reg_j, \reg_{j+1})$ or to the backbone itself.  
Therefore
$$ \eqalign{
\P(\# (B(j))^2 )
& = \P ((\sum_i 1 (\reg_j \leq i < \reg_{j+1} ) (\# (C_i) +1))^2)\cr
 &\leq 2(\P(M^\prime_\infty)^2 ) + \P((\sum_i 1 (\reg_j \leq i <
\reg_{j+1}) 
(\nu_{W_i} (\# (C)) (d-1) + 1))^2)\cr
 &\leq 2\P((M^\prime_\infty)^2) + c \P((\reg_{j+1} - \reg_j)^2)\hskip
1.0in \hbox{ (Theorem 4.4)}\cr
 &\leq c_{7.6} }$$

\no by the above and (7.10).  This proves (b), and also shows that
$\reg_j < \infty$ for all $j$ a.s.
\med
For (c), one more lemma is required.

\no{\bf Lemma 7.7}. For all $j \in \N$ and measurable 
$A \in [0, \infty )^{\B}$ we have 
$$\P((T(j;x), x\in \B) \in A|{\cal C}_{\reg_j})
 = \int^{\delta_0}_0 \P((T_0 (x), x\in \B ) \in A| W_0 = w )
d \Lambda(w) \Lambda ([0, \delta_0 ])^{-1}.$$

\no{\it Proof}.   Let $p((U_x, x\in \B - \{0\} ) \in \cdot \, |w)$ be a
regular conditional probability for
\no $(U_x, x\in \B - \{ 0\})$
given $W_0 = w$.  Let $U(j;x) = U((\spine_{\reg_j})\oplus x)$ and choose 
$F\subset \B - \{0\}$
 finite and $B$ a measurable subset of $[0, \infty )^F$.  We will prove
$$ \P((U(j;x), x\in F) \in B |{\cal W}_{\reg_j} ) (\omega ) = p ((U_x, x
\in F) 
\in B |W_{\reg_j} (\omega )) \hbox{ a.s. } \leqno (7.11)$$

\no Since $T(j; x) = \sum_{0\not= y\leq x} U(j; y) \alpha^{|y|}$, we can
 condition on ${\cal C}_{\reg_j} \subset {\cal W}_{\reg_j}$ and
use Lemma 7.5 to obtain the desired result.  Turning to (7.11), note that
for
 $A\in {\cal W}_{\reg_j}$ and $x_0 \in \B$,
$$ \P((U(j;x), x\in F) \in B, A, \spine_{\reg_j} = x_0) 
= \P((U_{x_0\oplus x} , x\in F) 
\in B, A, \spine_{\reg_j} = x_0),\leqno (7.12)$$
that $A\cap \{ \spine_{\reg_j} = x_0 \}$ is in $\bar {\cal F}_{x_0} = \sE_{x_0}
 \vee \sigma (W(x_0))$, and that $\{(U_{x_0\oplus x},
x\in F) \in B\} \in {\cal F}_{(x_0 , \infty)}$. The independence of 
$\sE_{x_0}$ and ${\cal F}_{(x_0, \infty )}$  and the inclusion 
$\sigma (W(x_0)) \subset {\cal F}_{(x_0, \infty )}$ therefore shows that 
(7.12) equals

$$ \eqalign{
\P(\P((U_{x_0\oplus x}, x &\in F) 
\in B|W (x_0))\ 1 (A, \spine_{\reg_j} = x_0))\cr
 & = \P( p ((U_x, x\in F) \in B|W_{\reg_j})\ 1 (A, \spine_{\reg_j} = x_0)).}$$

\no Sum over $x_0$ to obtain (7.11).\qed

\medskip Since ${\bf W}$ and $(Y_{n,i})$ are measurable functions of 
$\TT^{(0)}(\cdot ),$ 
there is a measurable map \break
$ r : [0, \infty )^{\B}
\rightarrow \N$ such that $\reg_1 = r(\TT^{(0)})$.  It is 
straightforward to check that

\no $\reg_{j+1} - \reg_j = r (T(j; \cdot ))$ for all 
$j\in \Z_+$ (use (a)).  Define $\Phi : [0, \infty )^{\B}
\rightarrow {\cal S}_0$ by
$$ \eqalign{
 & \Phi (\TT^{(0)} (\cdot )) = (\overline {D_0} (\TT^{(0)} (\cdot )), 
\TT^{(0)} (\cdot )|_{\overline {D_0}}),\ \hbox{ where}\cr
 & \overline {D_0} (\TT^{(0)} (\cdot )) = \{ x\in \B : \TT^{(0)} (x) 
< \lim_{n \rightarrow \infty} \inf \{ \TT^{(0)} (y) : |y| = n\}, |x| < r
(\TT^{(0)})\} .}$$

\no Then $(D_j , T(j; \cdot )) = \Phi (T(j; \cdot ))$ for all $j\in \Z_+$ 
and so Lemma 7.7 shows that if $A \subset {\cal S}_0$ is measurable, then
$$ \leqalignno{
\P((D_j, T(j ; \cdot )) \in A|{\cal C}_{\reg_j}) & 
= \int^{\delta_0}_0 \P (\Phi (\TT^{(0)}) \in A | W_0 = w ) 
d \Lambda (w) \Lambda ([0,
\delta_0])^{-1}\cr &=\int^{\delta_0}_0 \P((D_0, \TT^{(0)}) \in A |W_0 = w) 
d \Lambda (w) \Lambda ([0, \delta_0])^{-1}.}$$

To complete the proof of Theorem 7.1(c) it suffices to show that
$(D_j , T(j;\cdot))$ is ${\cal C}_{\reg_{j+1}}$-measurable for all 
$j\in \Z_+$. This reduces to showing that for a fixed $x \in \B$,
$$\leqalignno{
\TT^{(0)} ((\spine_{\reg_j})\oplus x) 1 (\TT^{(0)} ((\spine_{\reg_j}) 
  &\oplus x) < \alpha W_0, |x| < \reg_{j+1} - \reg_j) \cr
 &\hbox{ is ${\cal C}_{\reg_{j+1}}$-measurable.}&(7.13) }$$

\no If $n\in \Z_+$ and $x_0 \in \B$ satisfy $n+|x|< |x_0|$ then on 
$\{ \spine_{\reg_{j+1}} = x_0, \reg_j = n\}$ we consider the following two
cases:

\no{\it Case 1}.  $(x_0 |n ) \oplus x \leq x_0$.

Then $(\spine_{\reg_j}) \oplus x =x_0 | (n +|x|)$, $\TT^{(0)} (\spine_{\reg_j} \oplus x)
=
 \alpha W_0 - \alpha^{n+|x|} W_{n+|x|}$, and so
$$ \eqalign{ \{  \TT^{(0)} ((&\spine_{\reg_j}) \oplus x) \in B, 
\TT^{(0)} ((\spine_{\reg_j}) \oplus x)< \alpha W_0, \spine_{\reg_{j+1}} = x_0, 
\reg_j = n\}\cr
 & = \{ \alpha W_0 - \alpha^{n + |x|} W_{n+|x|} \in B, \reg_{j+1} 
= |x_0|, \spine_{|x_0|} = x_0 , \reg_j = n\}\cr
 & \in {\cal C}_{|x_0|} \q  \hbox{ (because\ } n+ |x| < |x_0|).}$$

\no{\it Case 2.}  Case 1 fails.

Then  there exists $m \in [n, |x_0| - 1) \cap \Z_+$, $i<d$ and 
$x^\prime$ in $\B$ (depending on $(x, x_0))$ such that 
$(\spine_{\reg_j}) \oplus x = e_i (x_0 |m+1) \oplus x^\prime$. Now
$\TT^{(0)} (\spine_{\reg_j}) \oplus x) < \alpha W_0$ if and 
only if $x^\prime \in C_{m,i}$ in which case
$\TT^{(0)} ((\spine_{\reg_j} ) \oplus x) = \alpha W_0 - \alpha^{m+1} a_{m,i} 
(x^\prime )$.  
Therefore
$$ \eqalign{ \{ \TT^{(0)} ((&\spine_{\reg_j} \oplus x ) \in B, 
\TT^{(0)} ((\spine_{\reg_j} \oplus x ) < \alpha W_0, \spine_{\reg_{j+1}} = x_0 
, \reg_j = n\}\cr
 &= \{ \alpha W_0 - \alpha^{m+1} a_{m,i} (x^\prime) \in B, 
x^\prime \in C_{m,i} , \spine_{|x_0|} = x_0 , \reg_{j+1} = |x_0|, \reg_j = n
\}\cr
 & \in {\cal C}_{|x_0|} \q (\hbox{ because\ } n \leq m < |x_0|).\cr}$$

\no Taking the union over $n< k-|x|$ and $|x_0|=k$ in the above cases we
have
$$\eqalign{ 
\{ \TT^{(0)} ((\spine_{\reg_j})\oplus x) &\in B, 
\TT^{(0)} ((\spine_{\reg_j} \oplus x)< \alpha W_0, \cr
 &\hbox{ for } |x| < \reg_{j+1} - \reg_j , 
   \reg_{j+1} = k\} \in {\cal C}_k  \cr}$$
 and so (7.13) follows.  \qed

We will now use Theorem 7.1 to complete the proof of the main limit
theorem, Theorem 6.4. Let $I(n)=[\reg_{k-1}, \reg_k )\cap \Z_+$
iff $\reg_{k-1} \leq n < \reg_k$ and let
$$Z_n=\sum\limits_{j=0}^\infty \# C_j 1\bigl( j\in I(n)\bigr) , n\in
\Z_+,$$
be the size of the ``regeneration block" spanning generation $n$. 
Theorem 7.1(a) implies that
$$\maxsize (n) \leq \minsize (n) + Z_n\leqno(7.14)$$
and (recall that $N(n) = \sum^{n-1}_{k=0} \#(C_k)$)
$$N(n)+n+1\leq \minsize (n) + Z_n. \leqno(7.15)$$

\no{\bf Lemma 7.8}. The sequence $\{ Z_n \} $ is bounded in $L^1$. 

\no{\it Proof}.  It follows easily from the exponential estimate (7.10)
on the tail of $\reg_k - \reg_{k-1}$  and the Renewal Theorem (see (4.16) 
in Ch. XI of Feller (1971)) that
$$\sup\limits_n \P \left( {\bigl( \# I(n)\bigr)}^q\right) 
< \infty \hbox{ for all } q>0. \leqno(7.16)$$  
Note that Theorems 2.5 and 4.4 together with the definition
of $\overline{\sC}_j$ imply
$$\P (\# C_j|\overline{\sC}_j) \leq (d-1) c_{4.3},\leqno(7.17)$$
and that
$$Z_n=\sum\limits_{j=0}^\infty \sum\limits_{k=1}^\infty \# C_j
1 \bigl( j,n \in [\reg_{k-1}, \reg_k)\bigr) .$$
Take means in the above to conclude that if $p>1$,
$$ \eqalign{
\P (Z_n) 
&\leq \sum\limits_{j=0}^{n-1} \sum\limits_{k=1}^\infty 
   \P \Bigl( 1 \bigl( j\in [\reg_{k-1}, \reg_k )\bigr) 
   \# C_j {(\reg_k -j)}^p {(n-j)}^{-p} \Bigr)\cr
&\qquad  +\sum\limits_{j=n}^\infty \sum\limits_{k=1}^\infty  
    \P \Bigl( 1 \bigl( \reg_{k-1}\leq n\leq j < \reg_k \bigr) 
    \P ( \# C_j|\overline{\sC }_j)\Bigr) \cr
&\leq \sum\limits_{j=0}^{n-1} \P \Bigl( \bigl( \# C_j\bigr)
    {\bigl( \# I(j)\bigr)}^p\Bigr) (n-j)^{-p}
     + (d-1) c_{4.3} \P \Bigl( \sum\limits_{j=n}^\infty
     1\bigl( j\in I(n)\bigr) \Bigr) \cr
&\leq \sum\limits_{j=0}^{n-1} {\P {\bigl( (\# C_j)^2\bigr) }^{1/2}
     \P {\bigl( \# I(j)}^{2p}\bigr) }^{1/2} 
     (n-j)^{-p} +(d-1)c_{4.3} \P \bigl( \# I (n)\bigr) .\cr}
$$
Use (7.16) and Theorems 2.5 and 4.6 to see that the final expression above
is uniformly bounded in $n$. \qed  

\no{\it Proof of Theorem 6.4}. For each choice of $L(n)$, writing 
$$L(n)= N(n)+n+1+\Delta_n,$$
the inequalities (6.10), (6.11), (7.14) and (7.15) imply that 
$|\Delta_n| \leq Z_n$ and so is bounded in $L^1$ by the previous Lemma. 
The Borel-Cantelli Lemma implies $n^{-2} \Delta_{n^2}  \rightarrow 0$ a.s. and 
so Example 6.3 gives a.s. convergence along the subsequence $\{ n^2 \} $. A 
standard interpolation argument completes the proof of (a). Since 
 $n^{-1/2} \Delta_n  \rightarrow 0$ in probability, (b) is now immediate from 
the Central Limit Theorem for $N(n)$ (Example 6.3). \qed 


\bigb {\bf 8. Some Concluding Remarks }

While equations (4.1) and (4.2) specify the law of $T(\infty)$
(see the Remark following \break  Lemma~4.1),
precise estimates on its distribution seem quite difficult to obtain.
It is however possible to derive some asymptotic results as $\alpha 
\uparrow 
1$. Our starting point is the following 

\med {\bf Proposition 8.1.} (Kingman (1975)). If $\alpha=1$ then
$n^{-1} T(n) \rightarrow  c_{8.1}(d)$ a.s. as $n \rightarrow \infty$, 
where
$c=c_{8.1}(d)$ is the unique root in $(0,1)$ of $dc e^{1-c}-1 =0$.

As $\alpha \in (0,1)$ will vary in the following, we will use notation 
such
as $T_\alpha(x)$ or $T_{\alpha}(\infty)$ to denote dependence on $\alpha$. 
Bear in mind 
that
the times $T_\alpha (x)$ are all defined on a common probability
space as sums of the same variables $U_x$ with different weights.

\med {\bf Theorem 8.2.} $(1-\alpha) T_\alpha (\infty)
\rightarrow c_{8.1}(d)$ as $\alpha \uparrow 1$ a.s. and in $L^1.$

\no {\it Proof.} Fix $\ee > 0$, and use Proposition 8.1  to choose
$K(\ee,\omega)$ such that $n^{-1} T_1 (n) > c_{8.1}-\ee$ for 
$n \geq K(\ee,\omega)$.
If $x \in \B$, then
$$ \eqalign { T_{\alpha}(0,x) &= \sum_{i=1}^{|x|} \alpha^i 
(T_1(0,x|i) -  T_1(0,x|i-1) ) \cr
&\geq \sum_{i=K(\ee)}^{|x|} (\alpha^i - \alpha^{i+1}) T_1(0,x|i) +
 \alpha^{|x|+1} T_1(0,x) \cr
&\geq (c_{8.1} - \ee) \sum_{i=K(\ee)}^{|x|} (1-\alpha) i \alpha^i  
 \cr
&\geq (c_{8.1} - \ee) (1-\alpha) {\alpha}^{K(\ee)}\sum_{i=0}^
{|x|-K(\ee)} i \alpha^i.}$$
Take the minimum over $|x|=M$ and let $M \rightarrow \infty$ to see that

$$ \liminf_{\alpha \uparrow 1} (1-\alpha) T_\alpha (\infty) 
\geq c_{8.1}(d) \q \hbox{ a.s. } \leqno (8.1) $$


For the other direction, note that $T_1 (n) / n$ is uniformly
$L^2$-bounded (being stochastically smaller than the
average of $n$ i.i.d. exponentials) and so the convergence in
Proposition~8.1 holds in $L^1$ as well.  Choose $N$ large enough 
so that $\P (T_1 (N)) \leq (c_{8.1} + \ee) N$.  Now choose random 
vertices $x_k$ inductively $N$ generations apart so that
$x_0 = 0$ and $x_{k+1}$ minimizes $T_1 (x_k , x_{k+1})$ among
all descendants of $x_k.$  The times $T_1(x_k,x_{k+1}) = 
T_1 (x_{k+1}) - T_1 (x_k)$ will be i.i.d.  Use the crude bound 
$$T_\alpha (\infty) \leq \sum_{k=0}^\infty \alpha^{kN} 
   T_1 (x_k , x_{k+1}) $$
together with summation by parts and the Strong Law (as above)
to see that \hfil\break $\limsup_{\alpha \uparrow 1} (1-\alpha) T_\alpha 
(\infty)
\leq c_{8.1} + \ee$ almost surely, and that
$$\P [(1 - \alpha) T_\alpha (\infty)] \leq \sum_{k=0}^\infty
   (1 - \alpha) \alpha^{kN} (c_{8.1} + \ee) N \rightarrow 
   c_{8.1} + \ee \leqno (8.2) $$
as $\alpha \uparrow 1$ for fixed $N$.  The first conclusion and (8.1) give
the required almost sure convergence.  Fatou's Lemma and (8.2)
then show that the mean value of $(1-\alpha)T_\alpha (\infty)$
approaches $c_{8.1}$ as $\alpha \uparrow 1$.  Convergence in
$L^1$ now follows.
\qed

\med {\bf Remark 8.3}. Recall from Lemma 4.1 that 
$e^t\P(T_{\alpha} 
(\infty) >t)$ increases to a finite limit $c_{4.1}(\alpha,d)$ as 
${t \rightarrow \infty}$. It is possible to show that
$$ \lim_{\alpha \uparrow 1} (1-\alpha) \log c_{4.1}(\alpha,d) 
= c_{8.1}(d). \leqno(8.3)$$
The lower bound is an easy consequence of (8.1) and the 
monotonicity of
$e^t \P(T_\alpha (\infty) > t)$  in $t$. The upper bound is more 
involved and we will not give a proof. It uses the anticipating 
equation~(4.2). Equation~(8.3) 
shows that the bound on $c_{4.1}$ in Lemma 4.1 is far from optimal. 
For 
example if $d=2$, then $c_{8.1} (d) \approx .23$ and~(8.3)~implies 
$c_{4.1} (\alpha , 2) \leq \exp ((.23 + \varepsilon ) / (1-\alpha ))$ 
for $\alpha $ close to $1$; Lemma 4.1 gives the same kind of 
bound but with $.23 + \varepsilon$ replaced by $\log 2 \approx .69$.

\bigbreak
We conclude this paper by mentioning an associated particle system. 
Set
$$ u_t(n) = \#( \partial C_t \cap \B(n)), \q t \geq 0, \, n \in \Z_+.$$
The \llq particles'' (i.e. sites in $\B$ on the boundary of the cluster 
$C_t$)
evolve independently: each $u$-particle at a site $n \in \Z_+$ dies at 
rate
$\alpha^n$, and is replaced by $d$ particles at $n+1$. The process 
$U=(u_t(.) : \, t \geq 0)$ captures the essential features of the DLA 
processes
$C_t$ and $A_n$: only the labels of the branches are lost. For various 
limit
theorems on the process $U$ in the case $\alpha >1$ see Aldous and 
Shields
(1988). 

To study the evolving cluster for $\alpha < 1$ it is more helpful to 
consider 
the following modification of $U$.  Define a random time change
$\sigma_t$ so that $C_{\sigma_t}$ always adds neighbours to its
deepest vertices at a constant rate: 
$$ \eqalign{ M_t &= \max\{n: \, u_t(n)>0 \},\cr
L_t &= \int_0^t \alpha^{-M_s} ds, \q 0<t<T(\infty), \cr
\sigma_t &= \inf \{ s: \, L_s > t \},\cr
V_t(n) &= u_{\sigma_t}(M_{\sigma_t}-n), \q t \geq 0, \, n \in \Z_+, 
\cr}$$
where $u_t(n)$ is taken to be zero for $n < 0$.
At each time $t$, $V_t$ is a function on $\Z_+$ counting how 
many vertices are in $\partial C_{\sigma_t}$ at each level 
below the highest one.
Note that $L_t < \infty$ if $t < T(\infty)$, and that   
$\lim_{t \uparrow T(\infty)} L_t = \infty$; thus $\sigma_t < 
T(\infty)$ for all $t \geq 0$. Straightforward calculations show that the 
process $\V = (V_t(.) : t \geq 0)$ evolves as follows: 

\item{(i)} Particles at site $n$, $n\geq 1$, die at rate $\alpha^n$ and 
are
replaced by $d$ particles at site $n-1$.
\item{(ii)} Particles at 0 die at rate 1, and are replaced by $d$ 
particles
at site $-1$. The whole configuration is then immediately shifted to the
right by 1 step.

We call $\V$ the ``tip process'': it describes the form of the cluster when 
viewed backwards from the tip. The size of the process $\V$ near $0$ arises 
from the interaction of two effects: first the strongly supercritical 
branching (at an accelerating rate as particles approach $0$), and secondly 
the right shifts, which move particles away from $0$, and so slow down their 
branching. Since the process $\V$ is a functional of $C$, it should be 
possible to deduce many properties of $\V$ from our results on $C$. However, 
just as some work was needed to obtain results such as Corollary 6.5 (giving 
the growth rate of 
 $A_n$) from the cluster decomposition, so also passing from $C$ to $\V$ is 
not completely straightforward. A further study of $\V$ may be the subject of 
a future paper: here we will just give a sketch proof that $\V$ is (in a 
certain sense) recurrent. 

For $f: \Z_+ \rightarrow \R$ and $\lambda \in (0,1)$, set 
$ ||f||_\lambda = \sum_{n=0}^\infty \lambda^n |f(n)|.$
  
\bigb {\bf Theorem 8.4.} Let $\lambda \in (0,1)$. Then
the process $||V_t||_\lambda $ is recurrent in the sense that
there exists $c_{8.2}< \infty$ such that
$$ \{t: ||V_t||_\lambda < c_{8.2} \} \q \hbox{ is unbounded.} $$

\no {\it Proof.} 
Recall from Section 7 the definition of the regeneration 
times $\reg_j$, and set 
$$ \eqalign{  \xi_j &=\#( B(j-1)), \cr
 S_j & = \inf \{ t \geq 0: M_{\sigma_t} = \reg_j \}. }$$
It is clear that $S_j < \infty$ for all $j$ and that 
$\lim_{j \rightarrow \infty} S_j = \infty$. Then

$$\eqalign{
||V_{S_j}||_\lambda &= \sum_{n=0}^{\reg_j} \lambda^{\reg_j-n} 
u_{\sigma_{S_j}}(n) 
\cr
 &= \sum_{i=1}^{j} \sum_n 1( \reg_{i-1} \leq n < \reg_i) \lambda^{\reg_j-n}
u_{\sigma_{S_j}}(n) + u_{\sigma_{S_j}}(\reg_j) 
 \leq \sum_{i=1}^j \lambda^{j-i} \xi_{i} + 1. \cr }$$

\no Here we have used the facts that $\reg_j\geq \reg_{j-1}+1$, and 
that the cluster
$A_n \cap \B(0,\reg_i)$ is frozen for $n \geq \maxsize(\reg_i)$. 
By Theorem 7.1
$\xi_j, j\geq 1$ are i.i.d. with $P(\xi_j^2) < \infty$. Thus (b)  follows
by comparison with the interval recurrent Markov chain
$Z_n = \sum_{i=1}^n \lambda^{n-i} \xi_{i}$.
\qed

\bigbreak 

\everypar{\hangindent=40pt}
\centerline {{\bf References}}

\no Aldous, D.J. (1991).  The continuum random tree, I.  {\it Ann. Probab. } 
{\bf 19}, 1 - 28. 

\no Aldous, D.J. and Shields, P. (1980).  A diffusion limit for a class of 
randomly-growing binary trees.  {\it Prob. Th. Rel. Fields} {\bf 79}, 509-542. 

\no Athreya, K.B. (1985) Discounted branching random walks. {\it Adv. Appl. 
Prob.} {\bf 17}, 53-66. 

\no Barlow, M.T. (1993) Fractals and diffusion limited-aggregation. {\it Bull. 
Sci. Math. } {\bf 117 }, 161-169. 

\no Bramson, M. (1978).  Minimal displacement of branching random walk. {\it 
Zeit. Wahr.} {\bf 45}, 89 - 108. 

\no Breiman, Leo (1968). {\it Probability}.  Addison-Wesley, Reading, Mass. 

\no Brennan, M.D. and Durrett, R. (1986).  Splitting intervals.  {\it Ann. 
Probab.} {\bf 14}, 1023-1036. 

\no Devroye, L. (1986).  A note on the height of binary search trees.  {\it 
Journal ACM } {\bf 33}, 489-498. 

\no Durrett, R. (1991).  Probability: theory and examples.  Wadsworth, 
Belmont, CA.  

\no Evans, W., Kenyon, C., Peres, Y. and Schulman, L. (1995). A critical 
phenomenon in a broadcast process.  {\it Preprint.} 

\no Feller, W. (1971). {\it An Introduction to Probability Theory and its 
Applications Volume II, 2nd edition}. John Wiley and Sons, New York. 

\no Kamae, T., Krengel, U. and O'Brien, G.L. (1977).  Stochastic inequalities 
on partially ordered spaces.  {\it Ann. Probab.} {\bf 5}, 899-912. 

\no Kesten, H. (1987).  How long are the arms in DLA?  {\it J. Phys. A} {\bf 
20}, L29-L33. 

\no Kesten, H. (1990).  Upper bounds for the growth rate of DLA.  {\it Physica 
A} {\bf 168}, 529 - 535. 

\no Kingman, J.F.C. (1975).   The first birth problem for an age-dependent 
branching process.  {\it Ann. Probab.} {\bf 3}, 790-801. 

\no Lawler, G.F. (1991).  {\it Intersections of Random Walks}.  Birkh\"auser, 
Boston. 

\no Liggett, T.M. (1985).  {\it Interacting Particle Systems}.  Springer-
Verlag, New York. 

\no Lorden, G. (1970).  On excess over the boundary.  {\it Ann. Math. Stat.} 
{\bf 41}, 521-527. 

\no Lyons, R. (1990).  Random walks and percolation on a tree.  {\it Ann. 
Probab.} {\bf 18}, 931 - 958. 

\no Mahmoud, H. (1992).  {\it Evolution of Random Search Trees}.  John Wiley 
and Sons, New York. 

\no Newman, C.M. (1980).  Normal fluctuations and the FKG inequalities.  {\it 
Commun. Math. Phys.} {\bf 74}, 119-128. 

\no Pemantle, R. and Peres, Y. (1994). Domination between trees and 
application to an explosion problem.  {\it Ann. Probab.} {\bf 22}, 180-194. 

\no Pittel, B. (1984).  On growing random binary trees. {\it J. Math. Anal. 
Appl. } {\bf 103}, 461-480. 

\no Pittel, B. (1985).  Asymptotical growth of a class of random trees, {\it 
Ann. Probab.} {\bf 13}, 414-427. 

\no Rio, E. (1995).  The functional law of the iterated logarithm for 
stationary strongly mixing sequences.  {\it Ann. Probab.} {\bf 23}, 1188-1203. 
 
\no Vannimenus, J., Nickel, B. and Hakim, V. (1984).  Models of cluster growth 
on the Cayley tree.  {\it Phys. Rev. B.} {\bf 30}, 391-399. 

\no Vicsek, T. (1989).  {\it Fractal Growth Phenomena. } World Scientific, 
Singapore. 

\no Witten, T.A. and Sander, L.M. (1981).  Diffusion limited aggregation, a 
kinetic critical approach. {\it Phys. Rev. Lett.} {\bf 47}, 1400-1403. 

\no Ziv, J. (1978).  Coding theorems for individual sequences.  {\it IEEE 
Trans. Inf. Theory } {\bf 24}, 405-412. 

\bigskip
\settabs 3 \columns
\+ Address for M. T. Barlow and E. A. Perkins: &&    Address for R. 
Pemantle:\cr
\+  Department of Mathematics                  &&     Department of 
Mathematics\cr
\+  University of British Columbia             &&     University of 
Wisconsin\cr
\+  Vancouver, B.C. V6T 1Z2                    &&     480 Lincoln Drive \cr

\+  Canada                                     &&     Madison, WI 53706 
U.S.A.\cr
\end